\documentclass[12pt]{article}

\usepackage{epic} 
\usepackage{eepic} 
\usepackage{latexsym}
\usepackage{amsthm}
\usepackage{amssymb}
\usepackage[all]{xy}
\usepackage{graphicx}
\usepackage[usenames]{color}
\usepackage{makeidx}
\usepackage{amsmath} 
\usepackage{rotating}
\makeindex

\providecommand{\noglossaryignore}[1]{}
\newcommand{\globalglossaryentry}[3]{\makebox[1.5in][l]{\tt $\backslash${#1}} 
\makebox[1.1in][l]{{$#2$}} \makebox[2.5in][l]{{#3}}\newline} 
\newcommand{\newcommandabbreviation}[3]{\newcommand{#1}{#2}%
\noglossaryignore{\globalglossaryentry{#3}{#2}{}}}
\newcommand{\renewcommandabbreviation}[3]{\renewcommand{#1}{#2}%
\noglossaryignore{\globalglossaryentry{#3}{#2}{}}}
\newcommand{\newcommandmacro}[4]{\newcommand{#1}{#2}%
\noglossaryignore{\globalglossaryentry{#3}{#2}{#4}}}
\newcommand{\gge}[3]{\noglossaryignore{\globalglossaryentry{#1}{#2}{#3}}}
\newcommand{\myaddress}%
{\parbox{3in}{\footnotesize \begin{center} 
Mathematics Department, City University, \\  
Northampton Square, London EC1V 0HB, UK.\end{center}}}
    
\newcounter{minidef}[section]
\renewcommand{\theminidef}{\thesection.\arabic{minidef}}
\newcommand{\mdef}{\refstepcounter{minidef} 
\medskip \noindent ({\bf \theminidef}) }
\newcounter{minicapt}

\newcounter{perchapter}[section]
\renewcommand{\theperchapter}{\thechapter.\arabic{perchapter}}

\newenvironment{pr}{\refstepcounter{perchapter}
\par \noindent
{\bf Proposition \theperchapter.}\it}

\newenvironment{theo}{\refstepcounter{perchapter}
\par \noindent
{\bf Theorem \theperchapter.}\it}

\newenvironment{lem}{\refstepcounter{perchapter}
\par \noindent
{\bf Lemma \theperchapter.}\it}

\newenvironment{cor}{\refstepcounter{perchapter}
\par \noindent
{\bf Corollary \theperchapter.}\it}

\newenvironment{de}{\refstepcounter{perchapter}
\par \noindent
{\bf Definition \theperchapter.}}

\newenvironment{rem}{\refstepcounter{perchapter}
\par \noindent
{\bf Remark \theperchapter.}}

\newenvironment{exa}{\refstepcounter{perchapter}
\par \noindent
{\bf Example \theperchapter.}}

\noglossaryignore{GREEK ETC.\newline}
\newcommandabbreviation{\e}{\epsilon}{e}        
\newcommandabbreviation{\lam}{\lambda}{lam}  
\newcommandabbreviation{\la}{\langle}{la}        
\newcommandabbreviation{\ran}{\rangle}{ran}
\newcommandabbreviation{\ha}{\#}{ha}             
\newcommandabbreviation{\rmap}{\rightarrow}{rmap}
\newcommandabbreviation{\aaa}{\alpha}{aaa}        
\newcommandabbreviation{\ab}{\alpha,\beta}{ab}
\newcommandabbreviation{\aab}{a(\ab )}{aab}       
\noglossaryignore{\newline RINGS\newline}
\newcommandabbreviation{\HH}{H \!\!\! I}{HH}              
\newcommandabbreviation{\C}{\mathbb C}{C}
\newcommandabbreviation{\N}{\mathbb N}{N}  
\newcommandabbreviation{\Z}{\mathbb Z}{Z}     
\renewcommandabbreviation{\Re}{\mathbb R}{Re}
\newcommandabbreviation{\R}{{\mathbb R}}{R}
\newcommandabbreviation{\Q}{\mathbb Q }{Q}
\renewcommandabbreviation{\H}{\mathbb H }{H}
\noglossaryignore{\newline SYMMETRIC GROUP\newline}
\def\Sym(#1){\Sigma(#1)}                  
\gge{Sym(-)}{\Sym(-)}{symmetric group on - objects}
\def\Sy(#1){\Sigma_{#1}}                  
\gge{Sy(-)}{\Sy(-)}{symmetric group irreducible -}
\def\sym(#1){\mbox{\LARGE s}(#1)}       
\gge{sym(-)}{\sym(-)}{symmetric group on - objects (variant)}
\def\sy(#1){\mbox{\LARGE s}({#1})}       
\gge{sy(-)}{\sy(-)}{symmetric group irreducible - (variant)}
\newcommandmacro{\cs}{\C \, \sy(n)}{cs}{symmetric group algebra over $\C$}
\noglossaryignore{\newline PARTITIONS/SETS\newline}
\newcommand{\Nset}[1]{\underline{#1}}
\gge{Nset\{-\}}{\Nset{-}}{set of natural numbers to -}
\def\nset(#1){ \{ #1 \}_{ \underline{n} }}
\gge{nset(-)}{\nset(-)}{a set $-\times\Nset$}
\def\ul(#1){_{\underline{#1}}}            
\gge{ul(-)}{{}\ul(-)}{subscript underline -}
\def\Ee(#1){{\bf E}_{#1}}                 
\gge{Ee(-)}{\Ee(-)}{set of equivalence relations on set -}
\def\Eee(#1){{\bf E}_{\{ #1 \}_{\underline{n}}}}  
\gge{Eee(-)}{\Eee(-)}{ditto for nset}
\def\Een(#1,#2){{\bf E}_{\{ #1 \}_{\underline{#2}}}}  
\def\Ssn(#1,#2){{\bf S}_{\{ #1 \}_{\underline{#2}}}}  
\def\Ss(#1){{\bf S}_{#1}}                 
\def\Sss(#1){{\bf S}_{\{ #1 \}_{\underline{n}}}}  
\def\bbc(#1){((\beta_1)(\beta_2)...(\beta_{#1}))}     
\newcommandmacro{\Ln}{{\Gamma}^{n}}{Ln}{large index set}
\newcommandmacro{\LnQ}{{\Gamma}^{n}_Q}{LnQ}{index set}
\newcommandmacro{\Zz}{\zeta}{Zz}{`shape' function}
        
\noglossaryignore{\newline PARTITION ALGEBRA\newline}
\def\ka(#1){\kappa_{#1}}                  
\def\Sm(#1){\Sigma_{#1}}                  
\newcommandmacro{\com}{\bullet}{com}{bullet composition}
\newcommandmacro{\enm}{\; e^n(\! m\! ) \;}{enm}{product of idempotents}
\def\Ai(#1){ A^{ #1 \cdot } }             
\def\Aij(#1,#2){ A^{ #1  #2 } }           
\newcommandmacro{\One}{\mbox{\bf $1 \!\!\! 1$}}{One}{algebra unit 1}

\newcommandmacro{\Bp}{B_p}{Bp}{partition basis}
\def\Bb(#1){B_p[#1]}                      
\def\Pp(#1){P_n[#1]}                      
\def\Ps(#1){P_n[#1] \! /}                 
\newcommandmacro{\Ph}{\hat{P}}{Ph}{P hat  algebra}
\def\Is(#1){\sim^{#1}}                    
\noglossaryignore{\newline MODULES\newline}
\def\Wm(#1){{\cal S}_{#1}}                
\gge{Wm(-)}{\Wm(-)}{Weyl module with index -}
\def\wm(#1,#2){{}_{#1}{\cal S}_{#2}}      
\gge{wm(-1,-)}{\wm(-1,-)}{Weyl module with index -}
\def\Ind(#1,#2,#3){\mbox{Ind}_{#1}^{#2}#3}
\gge{Ind(-1,-2,-)}{\Ind(-1,-2,-)}{induction}
\def\Res(#1,#2,#3){\mbox{Res}_{#1}^{#2}#3}
\gge{Res(-1,-2,-)}{\Res(-1,-2,-)}{restriction}
\newcommandabbreviation{\weyl}{standard}{weyl}
\newcommandabbreviation{\head}{\mbox{head }}{head}
\newcommandabbreviation{\Weyl}{Weyl}{Weyl}
\def\SS(#1){{\cal S}_{#1}}                
\gge{SS(-)}{\SS(-)}{Specht/Weyl module index -}
\def\LL(#1){{\cal L}_{#1}}                
\gge{LL(-)}{\LL(-)}{simple module index -}
\newcommand{\modl}[1]{\mbox{$#1$-}\mod}    
\gge{modl\{-\}}{\modl{-}}{set of left modules for ring -}
\newcommand{\modr}[1]{\mod\mbox{-$#1$}}    
\gge{modr\{-\}}{\modr{-}}{set of right modules for ring -}
\noglossaryignore{\newline FUNCTORS/MAPS\newline}
\newcommandmacro{\Gg}{{\cal G}}{Gg}{G Functor}
\newcommandmacro{\Fg}{{\cal F}}{Fg}{F Functor}
\newcommandmacro{\ra}{\rightarrow}{ra}{}
\def\ses(#1,#2,#3){0\ra #1 \ra #2 \ra #3 \ra 0}  
\gge{ses(1,2,-)}{\ses(1,2,-)}{\hspace{.5in} short exact sequence}
\def\starr(#1){ \stackrel{ #1 }{\longrightarrow} }
\gge{starr(-)}{\starr(-)}{}
\newcommandmacro{\doublerightarrow}{\; -\!\!\! -\!\!\!\!\!\! \gg \;}
{doublerightarrow}{}
\noglossaryignore{\newline PARTITION ALGEBRA MAPS\newline}
\newcommandmacro{\smap}{s}{smap}{`inclusion' map}
\newcommandmacro{\tmap}{t}{tmap}{$ P_n -> S_n$}
\newcommandmacro{\pmap}{\psi}{pmap}{$ S_n -> P_n $}
\noglossaryignore{\newline MISC.\newline}
\def\Amap(#1){{\cal A}_{#1}}              
\gge{Amap(-)}{\Amap(-)}{}
\def\Rr(#1){R_{#1}}                       
\gge{Rr(-)}{\Rr(-)}{restriction of E}
\def\Cr(#1){C_{#1}}                       
\gge{Cr(-)}{\Cr(-)}{restriction of E to N}
\newcommandmacro{\Tm}{{\cal T}}{Tm}{Transfer Matrix}
\def\On(#1){{\cal I}_{#1}}
\gge{On(-)}{\On(-)}{}
\newcommandmacro{\UU}{\underline{\sqcup}}{UU}{}  
\newcommandmacro{\UUU}{\sqcup}{UUU}{}  
\newcommandmacro{\Vq}{V_Q^{\otimes n}}{Vq}{Potts config. space}
\def\bs(#1,#2){\mbox{{\Large $\ast$}}^{#1}_{#2}} 
\gge{bs(-,-)}{\bs(-,-)}{general plumbing multiplier}
\newcommand{\ignore}[1]{}
\def\choo(#1,#2){ \left( \begin{array}{c} #1 \\ #2 \end{array} \right) }
\gge{choo(-1,-)}{\choo(-1,-)}{choose}
\newcommand{\Qed}{$\Box$}
\gge{Qed}{\mbox{\Qed}}{QED}
\def\staq(#1){\stackrel{#1}{=}}           
\gge{staq(-)}{\staq(-)}{}
\def\stam(#1){\stackrel{#1}{\rightarrow}} 
\gge{stam(-)}{\stam(-)}{}
\def\mat{ \left( \begin{array} }    
\def\tam{ \end{array}  \right) }
\gge{mat/tam}{...}{matrix delimiters}
\newcommand{\beq}{\begin{equation} }
\def\eql(#1){ \begin{equation} \label{#1} 
}
\newcommand{\eq}{\end{equation} }
\def\eqal(#1){\begin{eqnarray} \label{#1} }
\def\eqa{\end{eqnarray} }
\def\lab(#1){\label{#1}
}
\def\prl(#1){ \begin{pr} \label{#1} 
}

\def\theol(#1){ \begin{theo} \label{#1} 
}

\def\leml(#1){ \begin{lem} \label{#1} 
}

\def\corl(#1){ \begin{cor} \label{#1} 
}

\def\del(#1){ \begin{de} \label{#1} 
}

\def\reml(#1){ \begin{rem} \label{#1} 
}

\def\exal(#1){ \begin{exa} \label{#1} 
}

\gge{smeq\{-\}}{...}{small equation}
  
\gge{fneq\{-\}}{...}{very small equation}
\newcommand{\beqa}{\begin{eqnarray}}%
\newcommand{\eeqa}{\end{eqnarray}}%

\noglossaryignore{\newline HECKE/BLOB\newline}
\newcommandmacro{\Hnq}{H_n(q)}{Hnq}{ * freestanding symbol}
\newcommandmacro{\Hn}{H_n}{Hn}{      *-mod etc.}
\newcommandmacro{\A}{{\cal A}}{A}{}
\newcommandmacro{\Cwts}{C}{Cwts}{}
\newcommandmacro{\CA}{{\cal A}}{CA}{}

\newcommandmacro{\calA}{{\cal A}}{calA}{}
\newcommandmacro{\modi}{\mbox{Mod} }{modi}{was mod not modi!}
\newcommandmacro{\Wgen}{{\Bbb S}}{Wgen}{}
\def\ol(#1){\overline{#1}}
\newcommandmacro{\st}{\mbox{St}}{st}{}
  
\def\CMult(#1,#2){(#1:#2)}
\def\CM(#1,#2){( #1 : #2 )}
\def\FMult#1,#2{(#1:#2)}
\def\CF#1,#2{(#1:#2)}

\newcommandmacro{\Top}{\mbox{Top}}{Top}{}
\newcommandmacro{\Soc}{\mbox{Soc}}{Soc}{}
\newcommandmacro{\Head}{\mbox{Head}}{Head}{}
\newcommandmacro{\Filt}{{\cal F}}{Filt}{}
\newcommandmacro{\Mod}{\mbox{mod}}{Mod}{}
\newcommandmacro{\Resi}{\mbox{Res }}{Resi}{was without i!}
\newcommandmacro{\Indi}{\mbox{Ind }}{Indi}{was without i!}
  
\def\RR(#1,#2){R^{#1}_{#2}}  
\def\TT(#1,#2){T^{#1}_{#2}}

\def\bigplus{\mbox{\LARGE $+$}}

\newcommandmacro{\Ann}{\mbox{Ann}}{Ann}{}
\newcommandmacro{\Cen}{\mbox{Cen}}{Cen}{}
\newcommandmacro{\End}{\mbox{End}}{End}{}
\newcommandabbreviation{\semisimple}{semisimple}{semisimple}
\newcommandabbreviation{\Bratteli}{Bratteli}{Bratteli}
\newcommandabbreviation{\JBC}{Jones Basic Construction}{JBC}
\newcommandabbreviation{\pa}{partition algebra}{pa}
\newcommandabbreviation{\TM}{transfer matrix}{TM}
\newcommandabbreviation{\PM}{Potts model}{PM}
\newcommandabbreviation{\QSC}{quantum spin chain}{QSC}
\newcommandabbreviation{\Hamiltonian}{Hamiltonian}{Hamiltonian}
\newcommandabbreviation{\YS}{Young symmetrizer}{YS}

\newcommand{\mystufffont}{\textsc}

\newtheoremstyle{pu}
{7pt}%
{7pt}%
{\it}
{}
{}
{.}
{ }
{\thmnumber{({\bf #2}) }\thmname{\textsc{#1}}\thmnote{#3}}
\newtheoremstyle{puu}
{3pt}%
{3pt}%
{\rm}
{}
{}
{.}
{ }
{\thmnumber{({\bf #2}) }\thmname{\textsc{#1}}\thmnote{#3}}

\theoremstyle{pu}
\newtheorem{mmpr}[minidef]{Proposition}
\newtheorem{mlem}[minidef]{Lemma}
\newtheorem{mth}[minidef]{Theorem}
\newcommand{\mpr}[1]{\begin{mmpr} #1 \end{mmpr}}

\theoremstyle{puu}

\newtheorem{mex}[minidef]{Example}

\newcommand{\mupr}{\smallskip\noindent{\mystufffont{Proposition.}} }

\newcommand{\glossaryentry}[5]%
{\makebox[1in][l]{{$#2#5$}} \makebox[3.5in][l]{{#3}} 
 \makebox[1in][l]{{\tt $\backslash$#4 }}\pageref{#4}\newline}
\newcommand{\metacommand}[5]{\newcommand{#1}{#2} 
}
\metacommand{\Set}{{\bf Set}}{Category of sets}{Set}{} 
\metacommand{\Ab}{{\bf Ab}}{Category of abelian groups}{Ab}{} 
\metacommand{\Tail}{{\bf Soc}}{tail of module (aka Socle)}{Tail}{}
\metacommand{\rad}{{\bf rad}}{radical}{rad}{}
\metacommand{\iparts}{\Lambda}{set of young diagrams of degree $n$}{iparts}{_n}
\metacommand{\tabset}{\Lambda}{set of young tableau of shape $\nu$}{tabset}{_{\nu}}
\metacommand{\Simple}{{\tt L}}{Simple module}{Simple}{_{\nu}}
\metacommand{\PerM}{{\tt M}}{Permutation module}{PerM}{_{\nu}}
\metacommand{\rank}{\mbox{rank }}{rank of matrix/module $M$}{rank}{M}

\newcommand{\nome}[2]{\nomenclature{#1}{#2}}
\newcommand{\rB}{partial Brauer} 
\newcommand{\eP}{even partition} 
 
\newcommand{\el}{tonal partition} 
\newcommand{\El}{Tonal partition} 
\newcommand{\eli}{$l$-tone} 
\newcommand{\tone}{tone} 
\newcommand{\red}[1]{{\textcolor{red}{#1}}}
\newcommand{\redx}[1]{#1}
\newcommand{\redxx}[1]{}
\newcommand{\blu}[1]{{\textcolor{blue}{#1}}}
\newcommand{\blux}[1]{{{#1}}}

\title{
\El\ algebras: fundamental and geometrical aspects of representation theory 
}
\author{Chwas Ahmed${}^{*\ddagger}$, 
  Paul Martin${}^*$ and
  Volodymyr Mazorchuk${}^\dag$
  \\ ${}^\ddagger$ {\small Department of Mathematics, College of
    Science, University of Sulaimani,} \\ {\small Kurdistan Region, Iraq}
\\ ${}^*$ {\small School of Mathematics, University of Leeds, Leeds, UK}
\\ ${}^\dag$ {\small Department of Mathematics, 
                   Uppsala University, Uppsala, Sweden}}
\date{}
\begin{document}
\maketitle


\newcommand{\Tri}{\mbox{Tri}} 
\newcommand{\Dgroup}{{\mathcal D}}
\newcommand{\Aregular}{$A$-regular}
\newcommand{\nngraph}{\mbox{-graph}}
\newcommand{\hash}{\#}
\newcommand{\ex}{\sigma} 
\newcommand{\chiy}{\chi} 
\newcommand{\prop}{propagating}
\newcommand{\propno}{\hash^p}
\newcommand{\MM}{{\mathcal M}}
\newcommand{\chiyn}{\chiy_{}}

\newcommand{\rb}[2]{\raisebox{#1}{#2}}
\newcommand{\ig}{\includegraphics}


\renewcommand{\nome}[2]{}
\newcommand{\href}[2]{#2}
\newcommand{\highdetail}[1]{}
\newcommand{\complete}[1]{\check{#1}}
\newcommand{\nullseq}{{\mathfrak n}}
\newcommand{\RHS}{RHS}
\newcommand{\RTP}{RTS}

\newcommand{\mdeff}[1]{{\sc #1}}
\newcommand{\mdefff}[2]{{\sc #2}}
\newcommand{\mdeft}[4]{{ #2} & #3 $\;$ #4}


\newcommand{\Trace}{\mbox{Tr}}
\newcommand{\Power}{P} 
\newcommand{\Part}{P} 
\newcommand{\Pset}{{\mathsf P}}
\newcommand{\Pair}{J} 
\newcommand{\pcirc}{\circ} 
\newcommand{\dcirc}{\oslash} 
\newcommand{\ER}{E} 
\newcommand{\udot}{u}  
\newcommand{\bdot}{u^*}
\newcommand{\RB}{R}  
\newcommand{\PPart}{\RB}
\newcommand{\PPPart}{\RB'}
\newcommand{\PParto}{\RB^{o}}
\newcommand{\PParte}{\RB^{e}}
\newcommand{\pPart}{B}

\newcommand{\CC}{{\mathcal C}}
\newcommand{\rook}{partial}
\newcommand{\Rook}{Partial}

\newcommand{\specht}{{\mathcal S}} 

\newcommand{\muth}{{\textsc{Theorem}.}\ } 
\newcommand{\mulem}{{\textsc{Lemma}.}\ }

\newcommand{\myceil}[1]{\left \lceil #1 \right \rceil }
\newcommand{\floor}[1]{\left \lfloor #1 \right \rfloor }
\newcommand{\flip}{\star} 
\newcommand{\mE}{\mathtt{E}}
\newcommand{\mP}{\mathtt{P}}
\newcommand{\mD}{\mathtt{D}}
\newcommand{\kk}{K}  
\newcommand{\edge}{ .\newline --- WAITING FOR UPDATE BELOW HERE --- \newline}
\newcommand{\mm}{{\mathbf m}}
\newcommand{\iss}{S_l \times S_m}
\newcommand{\kiss}{k\Smm}
\newcommand{\kissin}{k(S_{i} \times S_{n-2i})}
\newcommand{\Wl}{W^l}
\newcommand{\Wlb}{W^l_b}
\newcommand{\Wll}{W} 
\newcommand{\Pl}{P^l}
\newcommand{\Al}{A^l}
\newcommand{\gammal}{\gamma^{l,n}}   

\newcommand{\lmu}{\underline{\mu}}

\begin{abstract}
For $l,n \in \N$
we  define  tonal partition algebra $P^l_n$ over $\Z[\delta]$.
We construct modules $\{ \Delta_{\lmu} \}_{\lmu}$ for $P^l_n$ over $\Z[\delta]$,
and hence over any integral domain
containing $\Z[\delta]$ 
(such as $\C[\delta]$),
that pass to a complete set of irreducible modules over the field of
fractions.
We show that $P^l_n$ is semisimple there.
That is,
  we construct for the tonal partition algebras
  a modular system in the sense of Brauer \cite{Brauer39}.
(The aim is to
investigate the non-semisimple structure of the tonal partition
algebras over suitable quotient fields of the natural ground ring,
from a geometric perspective.)
Using a `geometrical' index set for the $\Delta$-modules, 
we give an order with respect to
which the decomposition matrix over $\C$
(with $\delta \in \C^{\times}$) is upper-unitriangular. 
We establish several crucial properties of the $\Delta$-modules.
These include a tower property, with respect to $n$,
in the sense of Green \cite[\S6]{Green80}
and Cox {\it et al} \cite{CoxMartinParkerXi06};
contravariant forms with respect to a natural involutive
antiautomorphism; 
a highest weight category property;
and branching rules. 
\end{abstract}

\section{Introduction} \label{ss:intro}

\newcommand{\monotone}{monotone}
\newcommand{\mtone}{tone}
\newcommand{\catP}{{\mathcal P}}  
\newcommand{\CatP}{{\mathfrak P}} 
\newcommand{\catT}{{\mathcal T}}
\newcommand{\mPaa}[2]{\mP_{#1 \sqcup #2}} 

\begin{figure}
  \newlength{\iin}
  \setlength{\iin}{3cm}
\hspace{.51cm}
  \raisebox{-.76531in}{
\includegraphics[width=1.52\iin]{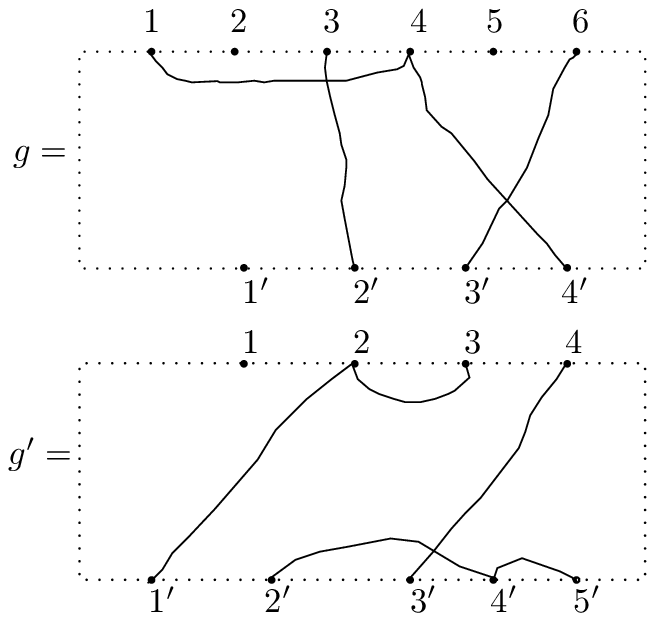}}
\hspace{.2593052in}
\hspace{.93052in}
$g|g' =\! $
\raisebox{-.631in}{
  \includegraphics[width=1.262\iin]{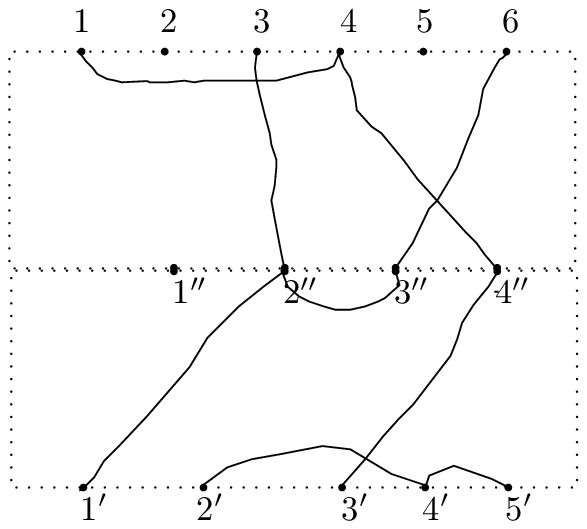}}
\\ \vspace{.51in} \\
\mbox{ } \hspace{.259063in}
\hspace{1.2in}
$g * g' \!\!= \delta  \!\!\!$
\raisebox{-.31in}{
 \includegraphics[width=1.352\iin]{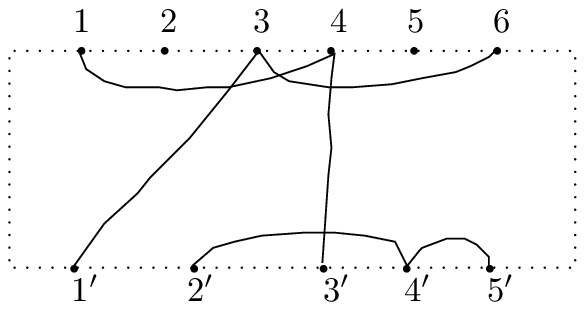}}
\caption{Schematic of composition in category $\catP$
  using
  graphs to represent partitions
  (see Appendix~\ref{ss:pdefa} for details). \label{fig:compx}}
\end{figure}

Fix $\kk$ a commutative ring and $\delta \in \kk$.
Let $\mP_S$ denote the set of set partitions of a set $S$.
The partition category $\CatP$ 
(as defined in \cite{Martin94})
is a  $\kk$-linear category: 
the objects of $\CatP$ 
are finite sets; and the hom-set $Hom_{\CatP}(S,T)$ has $\kk$-basis  
$\mPaa{S }{T}$.
(We recall the
$\delta$-dependent
category composition in \S\ref{ss:pdefa}.
A schematic is in Figure~\ref{fig:compx}.)   
The full subcategory on objects of form 
$\underline{n} = \{ 1,2,...,n \}$ is a skeleton in $\CatP$, denoted
$\catP$.
\newcommand{\calC}{{\mathcal C}}
In the notation for a category
in which 
$\calC = (objects,arrows,composition)$
we summarize the category $\catP$ as
$\catP = (\N_0 , K \mP_{\underline{n}\sqcup \underline{m}} , *)$. 


The category $\catP$ has many interesting subcategories, such as the
Brauer and Temperley--Lieb categories
\cite{Brauer37,TemperleyLieb71,Martin94,Martin08a,LehrerZhang14}. 
Here we study  subcategories chosen to serve as  testbeds for questions
in geometrical representation theory
(in the sense arising for example in
\cite{CoxDeVisscherMartin0609,Martin0915,EhrigStroppel13,EhrigStroppel14b}).

\medskip

Consider a part $p_i \in p \in \mPaa{S}{T}$. 
Define
\beq \label{eq:kerpi}
\ker(p_i) = |p_i \cap S | - | p_i \cap T |   .
\eq
For $l \in \N$
part $p_i$ is said 
  to be {\em $l$-\mtone} if $\ker(p_i) \equiv 0$
mod.$l$.
A partition is said to be $l$-\mtone\ if every part is $l$-\mtone.
It is routine to show
(see e.g. (\ref{eq:kerp}) in \S\ref{ss:catalg})
that the category composition in $\catP$ closes on the span of 
the subsets
$ \mP^l_{\underline{n}\sqcup\underline{m}} $
of $l$-\mtone\ partitions. 
Thus for each $l \in \N$ we have a subcategory of $\catP$, denoted
$\catP^l$:
\[
\catP^l = (\N_0 , K\mP^l_{\underline{n}\sqcup\underline{m}}  , *)
\]
There is also a corresponding subcategory of non-crossing partitions,
denoted $\catT^l$. Together we call these \el\ categories. 

\medskip

\newcommand{\rr}{{\mathsf r}}

For $\delta = m \in\N$
the partition algebra itself is Schur--Weyl dual to the `left'
symmetric group $S_\delta$ \cite{Martin91,Jones94,Martin2000}.
That is, the Potts/tensor action of 
$\catP$ \cite[\S 8.2]{Martin91}
is dual to the diagonal action of symmetric group $S_\delta$ on $Y^{\otimes n}$,
where $Y$ is the $\lambda=(\delta-1,1)$ Young module.
Neither action is faithful in general.
Let $\Lambda$ denote the set of all integer partitions,
$\Lambda^* = \Lambda\setminus \{ \emptyset \}$ and 
\[
\Lambda_i \; = \; \{ \lambda \vdash i \}
\]
denote the set of integer partitions of $i$. 
The natural  index sets for simple
modules over $\C$ are
$
\Lambda(P_n) = \sqcup_{i = 0,1,..,n} \Lambda_i
$
for the partition algebra $P_n$
and
$
\Lambda(\C S_m) \; = \; \Lambda_m
$
for $S_m$.
Let $\rr : \Lambda^* \rightarrow \Lambda$ denote the map that removes
the first part (i.e. removes the first row in the corresponding Young
diagram). 
Fixing $m$, the duality-induced connection between the index sets is
\beq \label{eq:r}
\rr : \Lambda( \C S_m ) \hookrightarrow \Lambda(P_n ) 
\eq
for sufficiently large $n$ (that $\C S_m$ acts faithfully)
\cite{GoodmanWallach98,Martin2000}.
This has a useful geometrical realisation --- see
\cite{MartinWoodcock98}. 


\newcommand{\mk}{m}

In case $\kk =\C$
the $\kk$-algebras $Hom_{\catP^l}(S,S)$ 
are isomorphic to subalgebras of the partition algebra 
studied by Tanabe \cite{Tanabe97}, 
Kosuda
\cite{Kosudax,
  Kosuda06,Kosuda08}
and Orellana \cite{Orellana07}.
Tanabe showed that the Schur-Weyl duality between the symmetric group 
$S_\mk$ and
the partition algebra $P_n(\delta)$ with $\delta =\mk \in \N$ 
generalises to a duality between various
reflection groups and partition algebras. 
Kosuda then studied the complex semisimple representation theory of these
algebras in the generic case (of $\delta \in \C$) and in certain cases
relevant for duality
\cite{Kosudax,
  Kosuda06}.  
Orellana also studied the representation theory from the duality
perspective \cite{Orellana07} (together with an elegant parallel study
of the `coloured' partition algebras).  


In the general case of the original $S_\mk/P_n$ duality
(just as for classical $Gl_\mk/S_n$ duality)
the partition algebra does not act faithfully on
tensor space. Indeed it clearly acts semisimply when $K=\C$, 
but it is
not generally a semisimple algebra. 
The way that the tensor space
action `sits inside' the full algebra $P_n$ is 
(representation theoretically) rather interesting
\cite{MartinSaleur94b,Martin}, and
relates nicely to the geometric-linkage approach to 
geometric representation
theory \cite{Jantzen,MartinWoodcock98}.
Here the aim is to 
investigate the lift of this geometric approach to the 
$l$-\mtone\ cases.
To this end we construct a corresponding tower of 
$\pi$-modular systems, in the sense of 
\cite{Brauer39,CurtisReiner90,Green80,Benson95,
  CoxMartinParkerXi06,
  Martin0915}.


For our modular system we need first a construction for ordinary
irreducible representations over a suitable `ordinary' ground field.
In fact we construct modules directly over an integral ground ring
--- we do this over $\Z[\delta]$,
but the 
domain of complex polynomials over the indeterminate $\delta$ will
be adequate for our immediate purposes --- and show
that they pass by base change to ordinary irreducibles (over the field
of fractions).
To do this we construct  contravariant forms with respect to a
natural involutive antiautomorphism; and 
determine cases where they are non-degenerate.
We show that the algebra is semisimple in these cases.
We then show that the corresponding decomposition matrix has a unitriangular
property with respect to a suitable (partial) order.
(To achieve this we must {\em establish} a suitable partial order.)
To verify that our order has the required properties we
proceed by showing
that a certain quotient algebra $A^l_n$ is semisimple over $\C$.
NB, This last step is an addition to the steps needed in the classical
$P_n$ and $B_n$ (Brauer algebra) cases. 
It fulfills our requirements, but it also
presents some interesting new features in the representation theory,
as we shall elucidate in \S\ref{ss:int}.

\medskip

\noindent
{\bf Overview:}
The integral part of the modular
tower representation theory of $P^l_n$ follows the same steps as for
$P_n$ in \cite{Martin94,Martin96}.
However it is more complex in the detail. 
The general representation theoretic machinery
is collected in \S\ref{ss:gen}.
In \S\ref{ss:catalg} we define the algebra. 
In \S\ref{ss:2}, \S\ref{ss:ideal2}
we construct a poset of ideals (in $P_n$ this is a
chain) with relatively small sections `controlled' by symmetric
groups.
In \S\ref{ss:pol} we give a polar decomposition of partitions in an algebra basis
that facilitates construction of standard module bases.
In \S\ref{ss:int} we construct our `standard' modules and cv forms. 
In \S\ref{ss:cvf}-\ref{ss:uut} 
we study the algebra that is the top section in a natural
tower structure, and hence derive the unitriangularity theorem.
In \S\ref{ss:branch} we give restriction rules for our standard
modules.
Relating these to induction rules, one has potential analogues of the
powerful translation functors of Lie theory \cite{Jantzen}
(to complete this picture
and hence connect to Kazhdan-Lusztig Theory, cf. \cite{Martin0915}, 
we need a linkage principle - this will be
discussed elsewhere). 

The main Theorems here are as follows.
\\
Theorem~\ref{th:UEU}, which shows that $P^l_{n-l} \cong W^l P^l_n W^l$
for suitable $W^l \in P^l_n$. 
This tells us that the module category of $P^l_{n-l}$ fully embeds in
that of $P^l_n$, by a `globalisation' functor. 
This in turn tells us that we can determine the
structure of these module categories iteratively on $n$.
To this end we  
construct canonical modules in each $n$
that are well-behaved under globalisation
--- `standard' modules. 
\\
Theorem~\ref{th:ssss}, Theorem~\ref{th:ssssc}
which show that each $P^l_n$ gives a modular system.
 \\
Theorem~\ref{de:remcx}: upper-triangularity of the standard module
decomposition matrix.
\\
Theorem~\ref{th:qh5}: that if $k=\C$ and $\delta\neq 0$ then
$P^l_n -\!\!\mod$ is a highest weight category.
\\
Theorem~\ref{th:brat1}: branching rules.


\medskip

Remark: 
Here we use the term {\em tone} as the generalisation to general $l$
of the notion of parity for congruence mod.2.
(Kosuda's name of `modular party algebra' also sounds harmonious, but does not
quite fit our purpose.)


\section{\El\ categories and algebras} \label{ss:catalg}


\begin{figure}
\[
  \includegraphics[width=.6042in]{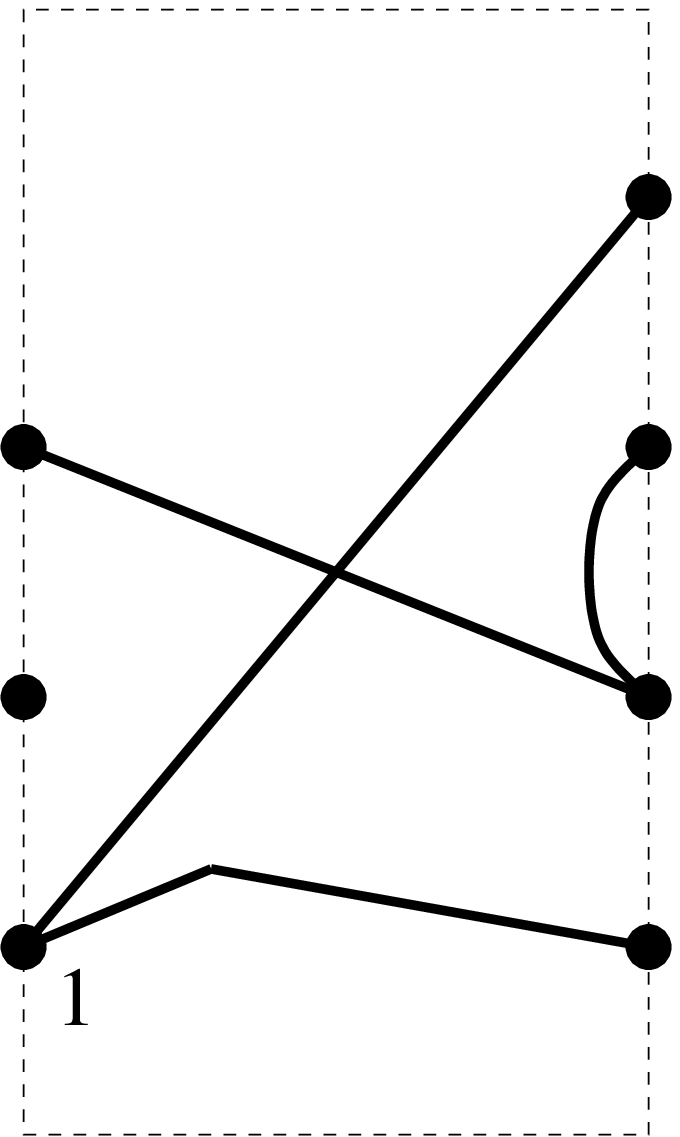}
\mbox{$\;\;\;\;$ rotate: $\;\;$ }
\includegraphics[width=1.02in]{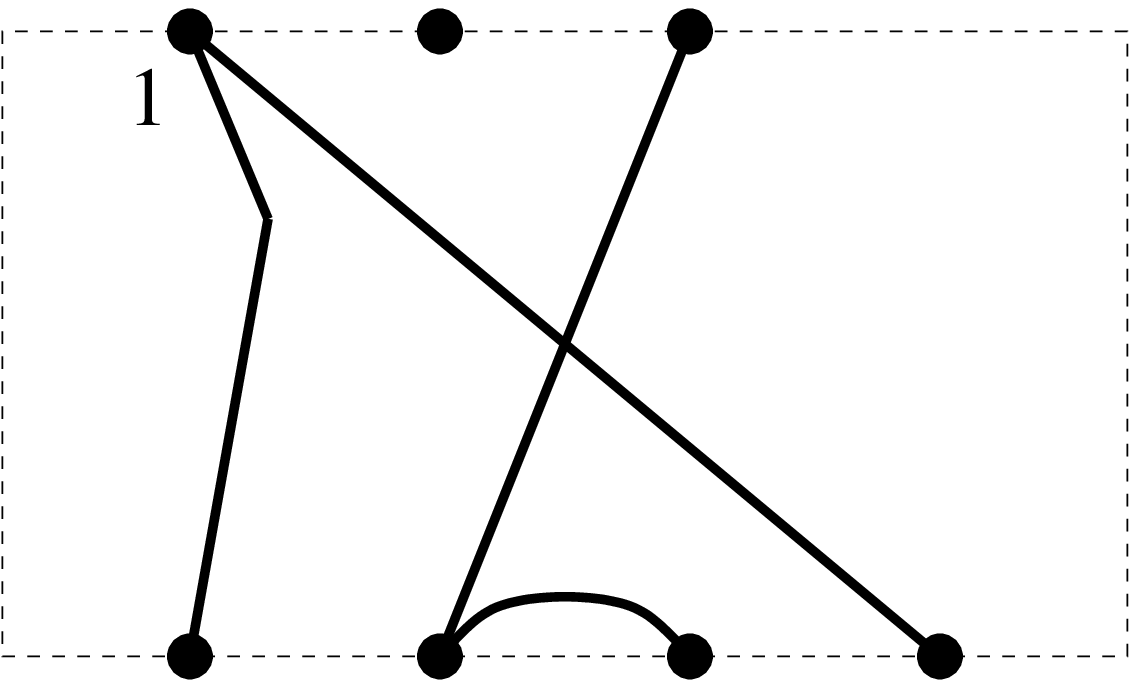}
\]
\caption{
  A partition picture
  (for partition
  \redxx{now I tried to make this consistent with appendix etc
  --- need to add a prime in figure!}
  $\{ \{ 1',4',1 \}, \{ 2',3',3 \}, \{ 2 \} \}$)
  drawn in different orientations.
  \label{fig:uppsalalondon}}
\end{figure}

\begin{figure}
\[
\includegraphics[width=1.2in]{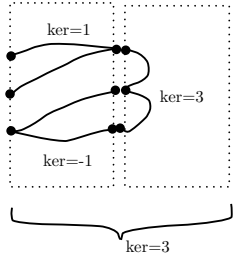}
\]
\caption{ Examples of partitions
(here drawn left-to-right)
  and parts with kernel numbers.
  \label{fig:kercount1}}
\end{figure}




\newcommand{\corange}{corange}
\newcommand{\range}{range}
\newcommand{\cora}{{\rm cora}}
\newcommand{\rang}{{\rm ran}}

Let $\catP = (\N_0 , \catP(n,m), *)$ denote the usual partition category
over a given commutative ring $\kk$, with parameter $\delta\in \kk$
\cite{Martin96}.
This is a skeleton in $\CatP$ obtained by restricting to objects 
$\underline{n} = \{ 1,2,..., n \}$.  
Thus $\catP(n,m) = \kk \mP_{n,m}$ where $\mP_{n,m}$ is the set of set
partitions of 
$\underline{n}         
\cup \underline{m'} $  
(with $\underline{n} = \{ 1,2,...,n \}$, 
$\underline{m'} = \{ 1',2',...,m' \}$).

We draw pictures of partitions in $\mP_{n,m}$
as for example in Fig.\ref{fig:uppsalalondon}
and \ref{fig:kercount1}.
Fig.\ref{fig:kercount1} also shows the kernel count as in
(\ref{eq:kerpi}).
Composition is as in Fig.\ref{fig:compx}. 

Write $\otimes$ for the usual monoidal composition in 
the category $\catP$. 

Write $P_n$ for the usual partition algebra $P_n= \catP(n,n)$,
and $\mP_n = \mP_{n,n}$ for the basis of partitions.



Write 
$\mP^{l}_{n,m} \subset \mP_{n,m}$ for the subset of
$l$-\tone\ 
partitions.


{\mth{ {\rm (cf. \cite{Tanabe97,Kosuda08})}
Fix $l \in \N$. The restriction of category $\catP$ to 
the span of 
$l$-\el s defines a
monoidal 
    subcategory, the {\em \el\ category} $\catP^l$. 
(Hence defining the \el\ algebras 
$
P^l_n = \catP^l(n,n) . 
$)
}}
\proof
Consider the product of composable partitions $p,p'$ in $\catP$.
Note (e.g. from \S\ref{ss:pdefa})
that in   
the definition of the product $p p'$
one first forms the concatenation $p|p'$,
then discards the `middle' vertices to form $pp'$.
Thus whenever a part $(pp')_i$ is formed in
composition the process is that (in some number of instances)
two vertices, one in some 
$\rang(p_i)=p_i \cap \underline{m}'$ and one in some 
$\cora(p'_j) = p'_j \cap \underline{m}$, 
are identified and then discarded from some union of parts. 
Thus 
\beq \label{eq:kerp}
\ker ((pp')_i) \; = \; \sum_{\pi} \ker( \pi )
\eq
--- sum over parts from $p,p'$ involved in $(pp')_i$
(cf. Fig.\ref{fig:kercount1}).
Thus if the incoming parts are all \eli\ ($\ker$ divisible by $l$), 
then the new part is again \eli.

For the monoidal property note that if $a,b$ are $l$-tone then so is
$a \otimes b$. 
\qed



\medskip

\newcommand{\bb}{b}
\newcommand{\ee}[2]{{\mathtt e}_{#1}^{(#2)}}
\newcommand{\AAA}[1]{{\mathsf A}^{#1}}
\newcommand{\Onne}[1]{1_{#1}} 

\begin{figure}
$1_n=$
\raisebox{-.431in}{
  \includegraphics[width=1.82in]{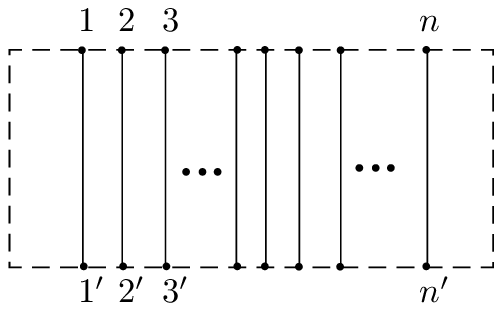}}
, \hspace{.421in}
$\AAA{i \; i\! +\! 1} = $
\raisebox{-.431in}{
  \includegraphics[width=1.82in]{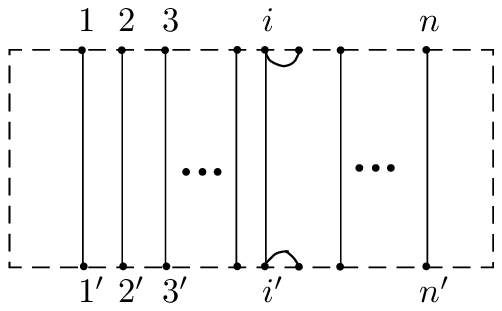}}
,
\\ \vspace{.21in} \\
\mbox{ } \hspace{.421in}
$\varepsilon_i^{(n)} = $
\raisebox{-.431in}{
  \includegraphics[width=1.82in]{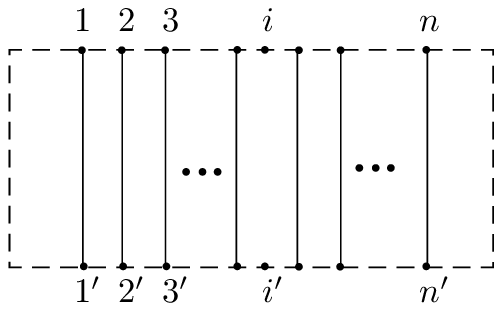}}
    ,\hspace{.421in}
$\ee{i}{n} = $
\raisebox{-.431in}{
  \includegraphics[width=1.82in]{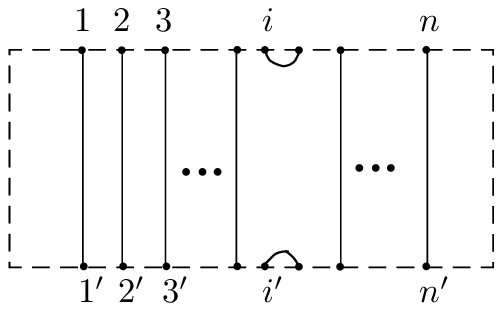}}
\caption{Special elements in the partition algebra $P_n$. \label{fig5}}
\end{figure}

\mdef 
For $l \in \N$ define
\[
b^l \; = \; \{\{ 1,2,..,l,1',2,',..l' \}\} \in \mP_l
\]
Consider Fig.\ref{fig5}.
Define
$u = \ee{1}{2}$,
$a=\AAA{12} \in \mP_2$.
Note that $b^1 = 1_1$, $b^2 = a$, and for $l>1$  we have
$
\bb^l  = \; \AAA{12} \AAA{23} ... \AAA{l\! -\! 1 \; l} \;\; \in \mP_l .
$
Fixing $n$, 
define
$U = \ee{1}{n}$ 
as the  
partition depicted in Fig.\ref{fig5}. 
We have 
\[
U = u \otimes \Onne{n-2}
\]
For $p \in \mP_{n,m}$ write
$
p^\star
$
for the `flip' image in
$\mP_{m,n}$
(given by $i \leftrightarrow i'$ in $p$) \cite{Martin94}.



\mdef \label{de:even1}
A set partition is {\em even} if all its parts are of even order. 
For example, all the partitions depicted in Fig.\ref{fig5} are even except for
$\varepsilon_i^{(n)}$.  
Write $\mE_n$ for the set of \eP s in $P_n$.
Note that $\mE_n = \mP^2_n$. 
Write $E_n$ for $P^2_n$. 


\mdef
For given $l$ let $w$ denote the unique 
partition
in  $\catP^l(l,0)$; and $w^\flip \in \catP^l(0,l)$.
Let $v$ denote specifically the unique partition in $\catP^2(2,0)$ and
$v^*$ the unique partition in $\catP^2(0,2)$.
Let $\sigma_1$ denote the unique elementary transposition in
$\catP(2,2)$. 
Then $E := \catP^2$ is generated as a linear monoidal category
by $1_1$, $v$, $v^\flip$, $\sigma_1$ and $a$.


\mdef
\label{de:flip}
Recall the `flip' antiisomorphism of $P_n$, denoted 
$p \mapsto p^{op}$ (given by $p \mapsto p^\star$ on partitions). 
This fixes the subset $P_n^l$, so that 
$P^l_n$ is isomorphic to its opposite. 

Note that the symmetric group $S_n$ is a subgroup of $P_n$ that is
also in $P^l_n$, and that the restriction of the flip antiautomorphism
to this subgroup is the usual $g \mapsto g^{-1}$ antiautomorphism.

\mdef
\label{de:latflip}
Let $w_n \in S_n$ denote the
order reversing (Coxeter longest word) element --- an
involution. Recall the lateral-flip endofunctor on $\catP$ given on
partitions $p \in \mP_{m,n}$ by
$
p \mapsto \bar{p} = w_m p w_n .
$
This takes a partition to its mirror image.
Note that the endofunctor fixes $\catP^l$, and indeed $P^l_n$.


\newcommand{\SLondon}{South London} 
\newcommand{\London}{London}


\newcommand{\tanabereview}{{
\subsection{Review of Theorems on $P^l_n$ from 
   Tanabe, Kosuda and Orellana}


{\em 
This section is temporary --- just organising what is already known,
to get referencing etc correct.}

\medskip

In \cite{Tanabe97} Tanabe introduces $P^l_n$ as the centralizer
algebra of the natural action of a unitary reflection group on tensor
space (generalising the $l=1$ case where the group is the symmetric group). 

In \cite{Kosuda06,Kosuda08} Kosuda first introduces the party algebra
--- the regular case of $P^l_n$; then gives a presentation for
$P^l_n$. 
In \cite{Kosudax} Kosuda computes the generic irreducible
representations in the general case. 

In \cite{Orellana07} Orellana ...

...



}}

\section{Basic
         properties of the algebra $\Pl_n$} \label{ss:2}

Here
 we develop an analogue of the propagating ideals of
$P_n$ as in \cite[\S6.1]{Martin94}.

\subsection{Set partitions: Co-$i$ parts and propagating numbers}

\newcommand{\has}{\hash^{s}}
\newcommand{\hashi}{\hash^{-}}

\mdef Recall that, for any $n$,  we write 
$\hash : \mP_n \rightarrow \N_0$ for the map taking a partition
$p \in \mP_n$ to the number of propagating parts in $p$.
Note that  a product of partitions in $\catP$ is a scalar times a partition.
The map $\hash$ extends to apply to a product of partitions 
in the obvious way.
Recall the `bottleneck principle':

{\mlem{{\rm{ \cite[\S6.1]{Martin94}}}
For $p,p' \in \mP_n$ we have $\hash( p p') \leq \hash(p)$. \qed
}}

\newcommand{\side}{side} 

\mdef
Fix $l$. 
For $p \in \mP^{l}_{n}$ a propagating part $s \in p$ is {\em co-$i$}
if the restriction of $s$ to one
`\side' 
of the underlying set
(the set $\underline{n}$ say) 
is of  order congruent to $i$ mod.$l$.
(Note that the definition is independent of the choice of \side\ 
in case $p$ in $\mP^{l}_n$ but not in general in $\mP_n$.)  

Define $\hash^i : \mP^{l}_n \rightarrow \N_0 $ so that 
$\hash^i (p)$ is 
the number of co-$i$ propagating parts. 
Define {\em propagating vector}
\[
\hashi (p) \; = \; (\hash^{1} (p) , \; \hash^{2} (p) , \hash^{3} (p) ,
... , \hash^{l} (p))
\]
Example: For $l=3$, 
$
\hashi(1_n) = (n,0,0) .
$


\label{de:pA12}
Consider $p \in \mP^{}_n$. Then $p \AAA{12}$ is a partition similar to $p$
but with the parts containing vertices
    {$1', 2'$}
combined.
Thus if vertices
{$1', 2'$} intersect at most
one propagating part in $p\in \mP^{l}_n$ then
$\hashi(p \AAA{12}) = \hashi(p)$.


\subsection{The index set $\gammal$ and corresponding partitions}

For $l \in \N$ and
$\mm  = (m_1, m_2, ..., m_l ) \in \Z^l$
define $r_\mm = \sum_{i=1}^l i m_i$.
For given $l$ and $n$, define 
\[
\gammal = \{ \mm \in \N_0^l \; : \; (n- 
r_\mm )/l \in \N_0 \}
\]


For $\mm \in \gammal$ define a set partition 
\beq \label{eq:amm}
 a^{\mm} \; 
 := \; a^{\mm}_n \; 
= \; \left( \otimes_{i=l}^1 (b^i)^{\otimes m_i} \right) \otimes
        (ww^\flip )^{(n-r_\mm )/l}
\eq
Note that $a^\mm \in \mP_n$. 
For example 
$$
a^{(4,4)}_{16} \;  = \;\;\; 
\raisebox{-.21in}{\includegraphics[width=4.9cm]{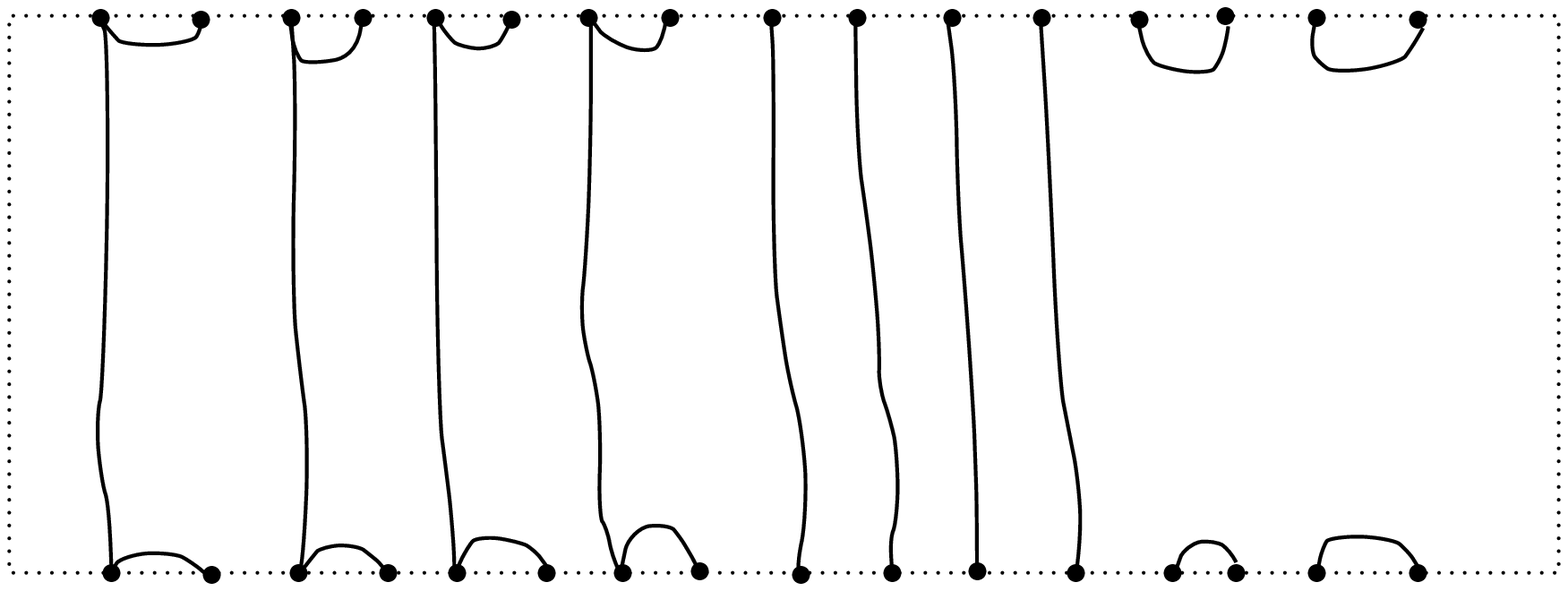}} .
$$
$$
 \; a^{(4,4,2)}_{24} \; 
  = \;\;\; 
\raisebox{-.21in}{\includegraphics[width=6.4cm]{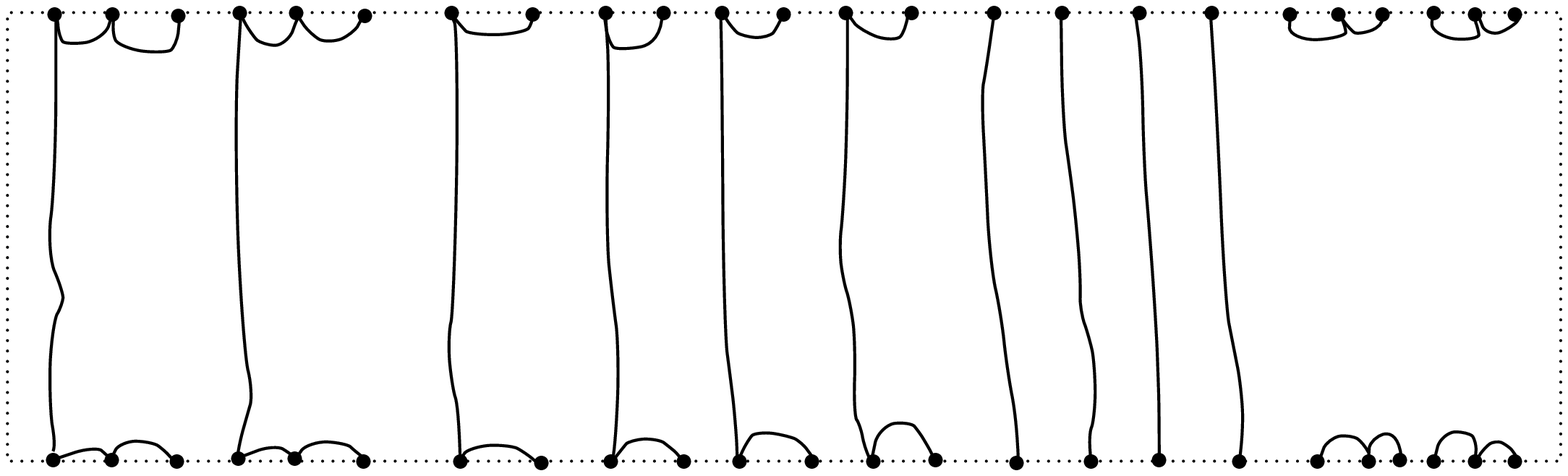}} .
$$
Thus
if $\delta$ is invertible then $a^{\mm}$ is a
(not necessarily normalised) idempotent in $P_n$.


If $\mm \neq 0$, so that $a^\mm_n$ has at least one propagating part,
we also define for each $a^{\mm}_n$ a partition $b^{\mm}_n$,
obtained from $a^{\mm}_n$ by combining the last (`rightmost')
propagating part with all the non-propagating parts.
Thus:
\[
\; b^{(4,4,2)}_{24} \;
  = \;\;\; 
\raisebox{-.21in}{\includegraphics[width=6.4cm]{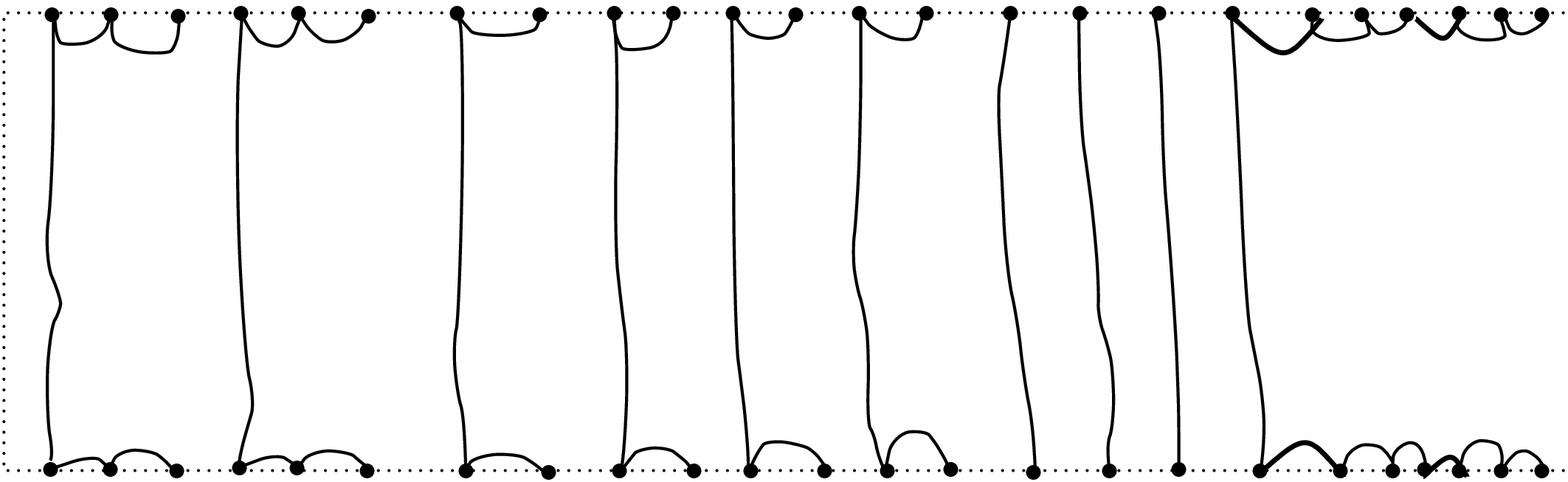}} .
\]
Note that $b^{\mm}_n$ is idempotent in $P_n$ (for any $\delta$).
Also
\beq \label{eq:abaa}
b^\mm a^\mm b^\mm = b^\mm
\qquad \mbox{ and } \qquad   a^\mm b^\mm a^\mm = a^\mm
\eq
For example
\[
\; b^{(4,4,2)}_{24} a^{(4,4,2)}_{24} b^{(4,4,2)}_{24} \;
  = \;\;\; 
\raisebox{-.21in}{\includegraphics[width=6.4cm]{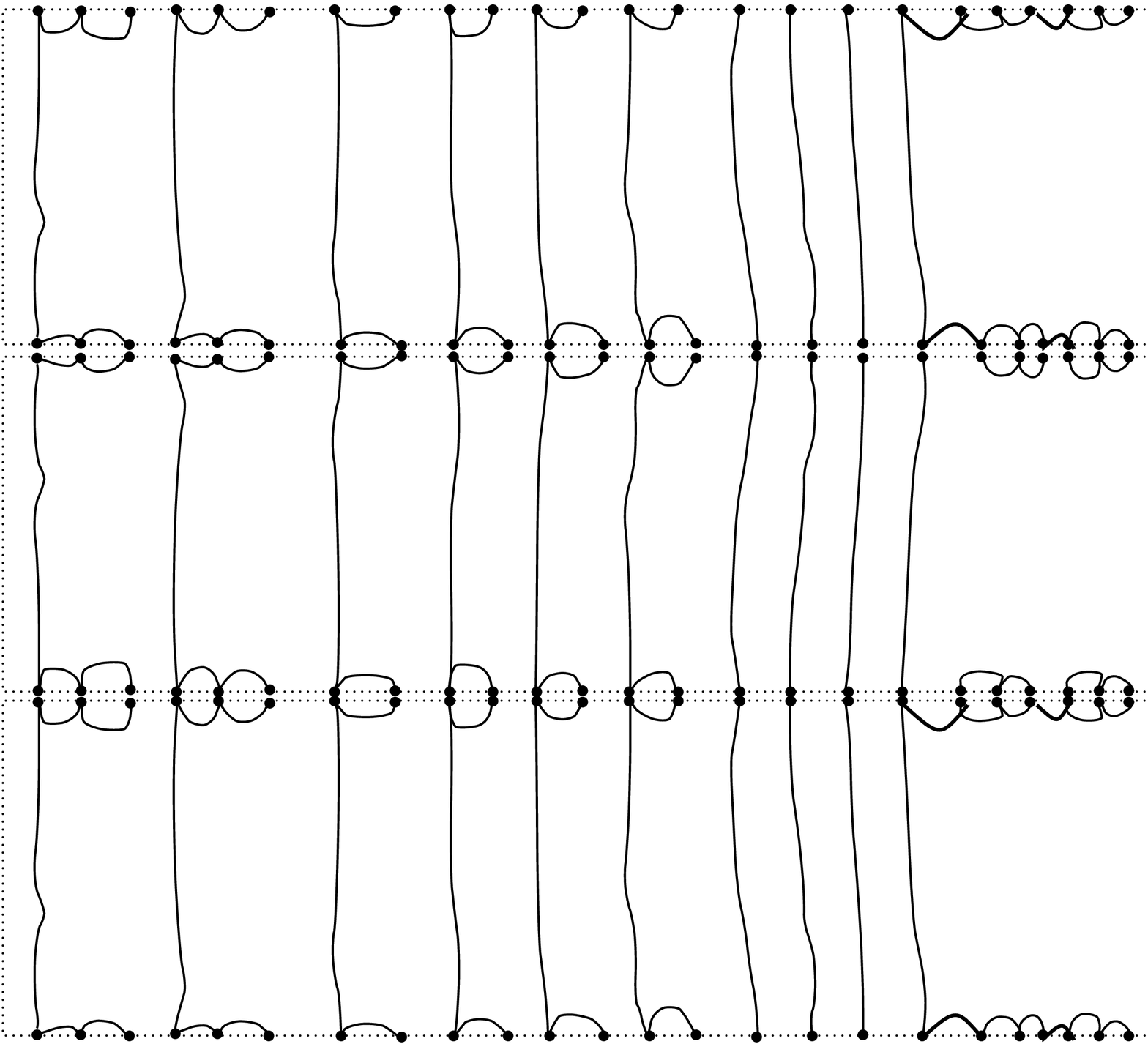}} .
\]

\subsection{Poset structure on $\gammal$} \label{ss:gammalposet}


\mdef \label{de:gammal}
Fix $l$. Define  $V \subset \Z^l$ as follows. 
Define $v_{ij}$ for $1 \leq i \leq j \leq l$ by 
\[
v_{ij} = (0,0,...,0,
   \underbrace{1}_{i+j},0,...,0,\underbrace{-1}_j,0,...,0,\underbrace{-1}_i,0,...,0)
\;\;\in \Z^l
\]
where the index $i+j$ is understood mod.$l$.
Then $V=\{ v_{ij} \}_{i,j}$.
In particular
\[
v_{ii} = (0,0,...,0,
   \underbrace{1}_{2i},0,...,0,\underbrace{-2}_i,0,...,0)
\qquad \mbox{ and } \qquad 
v_{ll} = (0,0,..,0,-1)
\]
Note that there are a total of $\frac{l(l-1)}{2} +1$ of these vectors
in $V$.

Define a poset structure on $\gammal$ by $\mm \geq \mm'$ if 
$\mm' - \mm$  
lies in the nonnegative integral span of $V$. 
For example 
$(9,0,0) > (7,1,0)$ since $-(9,0,0)+(7,1,0) = (-2,1,0) = v_{11}$.

\begin{figure}
\includegraphics[width=4.12cm]{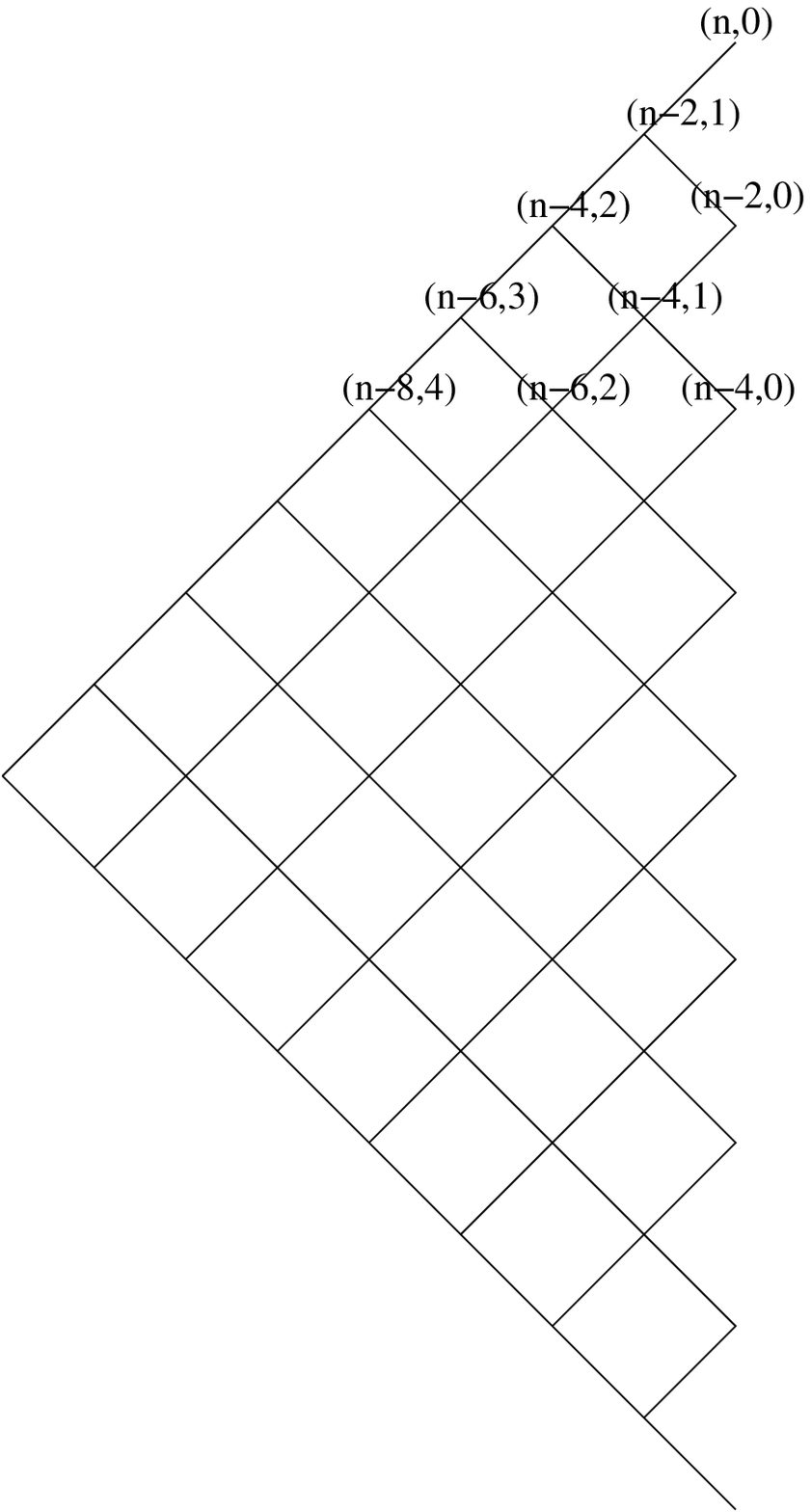}
\qquad
\includegraphics[width=4.2cm]{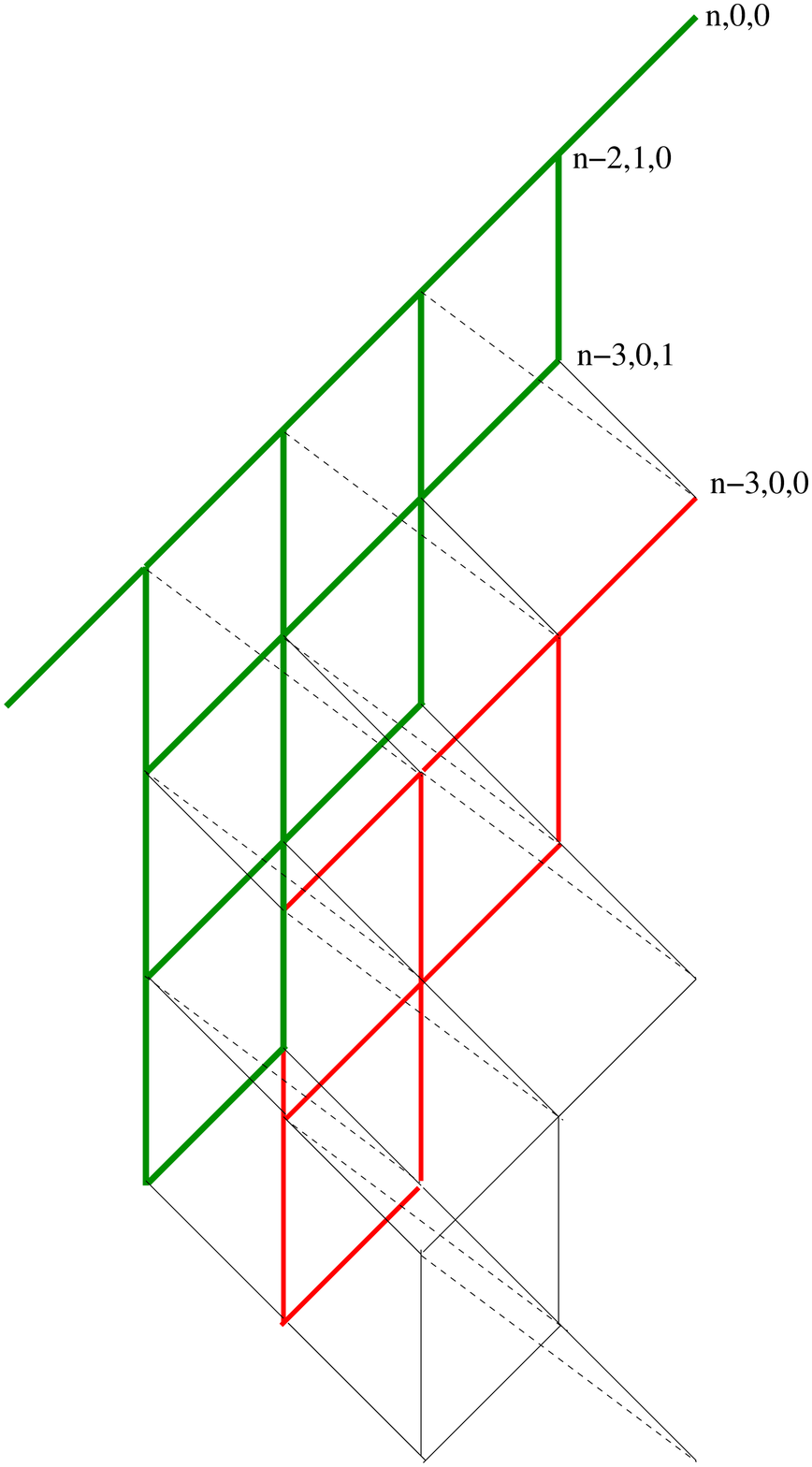}
\qquad
\includegraphics[width=4.2cm]{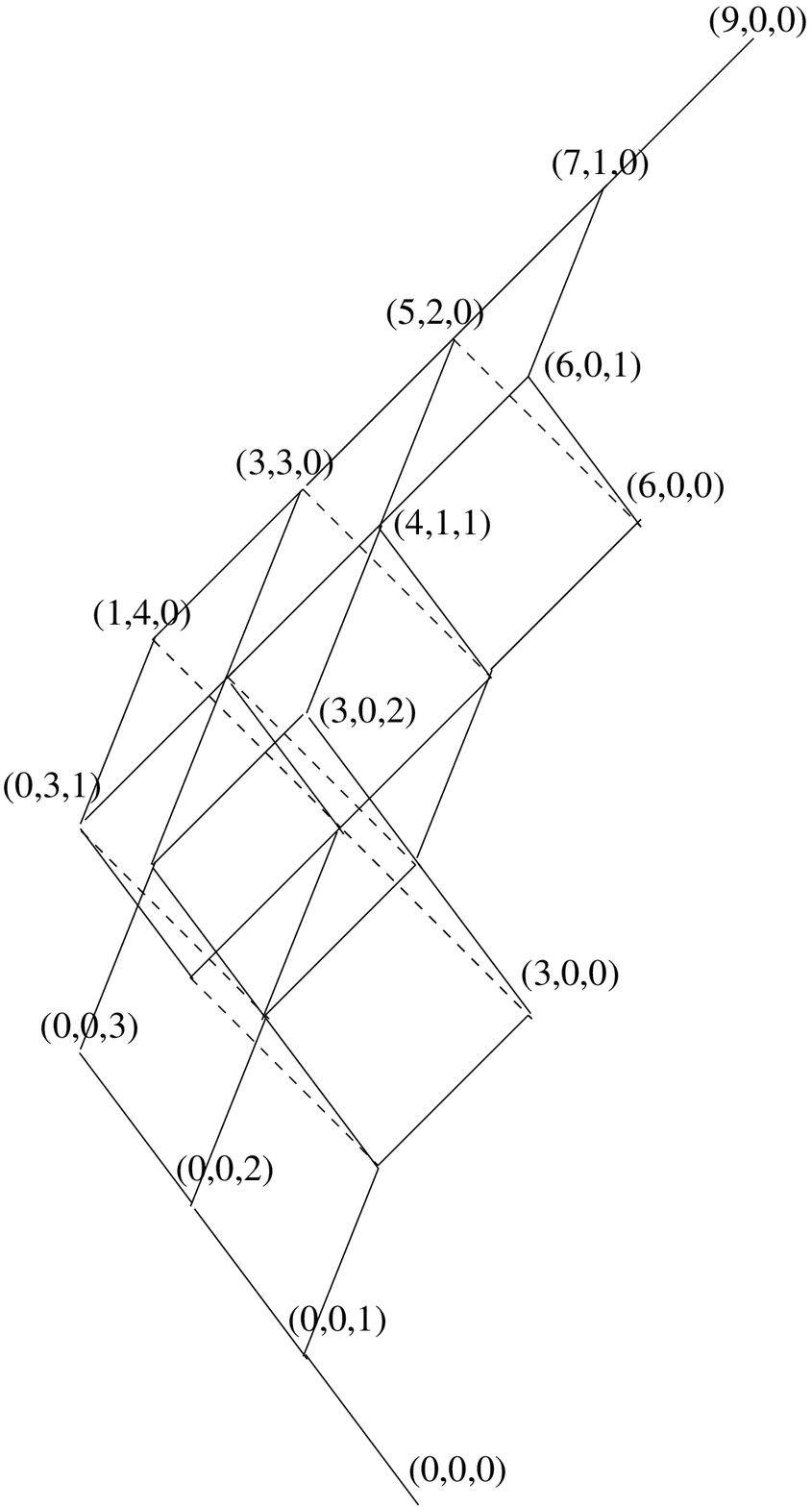}

\caption{Hasse diagram for  inclusion poset 
  for ideals $\Pl_n a^{\mm} \Pl_n$ in cases $l=2$ (left)
  and $l=3$ (right). 
  Vertex $\mm$ corresponds to ideal  $\Pl_n a^{\mm} \Pl_n$.
\label{fig:eventriPlattice1}}
\end{figure}

Note that $(n,0,0,...,0)$ is the unique top element in $\gammal$ for
any $l$. The Hasse diagrams in the cases $l=2,3$ are indicated in
Fig.\ref{fig:eventriPlattice1}.
See also Fig.\ref{fig:eventriPlattice12}. 


\mdef
Note that every element of $\gammal$ is in the positive cone of 
$(n,0,0,...,0)$ with respect to the subset 
$V' = \{ v_{11}, v_{12}, ...,v_{1 l} \}$. 
Thus $\gammal$ includes an $l$-dimensional lattice (in the crystal
lattice sense). 
The subset $V'$ is manifestly a basis for the underlying $\R^l$ containing
$\Z^l$. It follows that none of the remaining vectors in $V$ are
$\R$-linearly independent of $V'$. 
However we claim they are positive-integrally independent. 
Typically for every three sides in a cube in $\gammal$ then there is
an element in $V \setminus V'$ that is the main diagonal in this cube.
For example with $l=3$
\[
v_{11} + v_{22} = v_{12} + v_{13}
\]
That is to say, $v_{22} = v_{12} +v_{13} -v_{11}$, a non-positive
combination of basis elements. 
This example of $v_{22}$ corresponds to the dashed lines in
Fig.\ref{fig:eventriPlattice1}. 


\newlength{\ccm}
\setlength{\ccm}{0.79cm}

\begin{figure}
$n=4$ \includegraphics[width=1.9013\ccm]{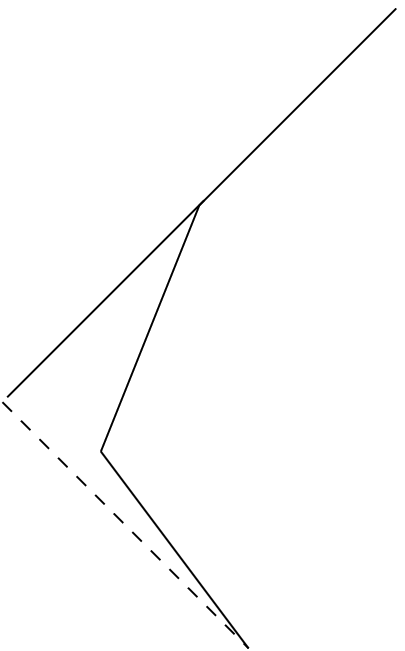} \hspace{.43in}
$n=5$ \includegraphics[width=2.13\ccm]{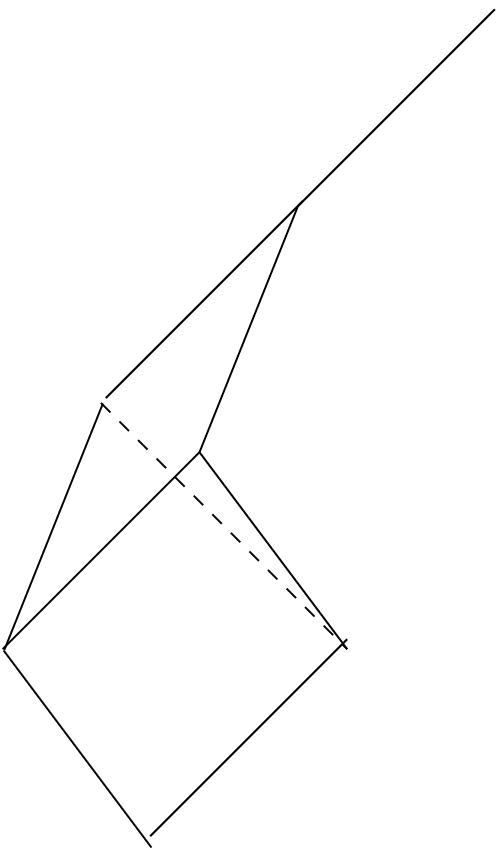} \hspace{.43in}
$n=6$ \includegraphics[width=2.463\ccm]{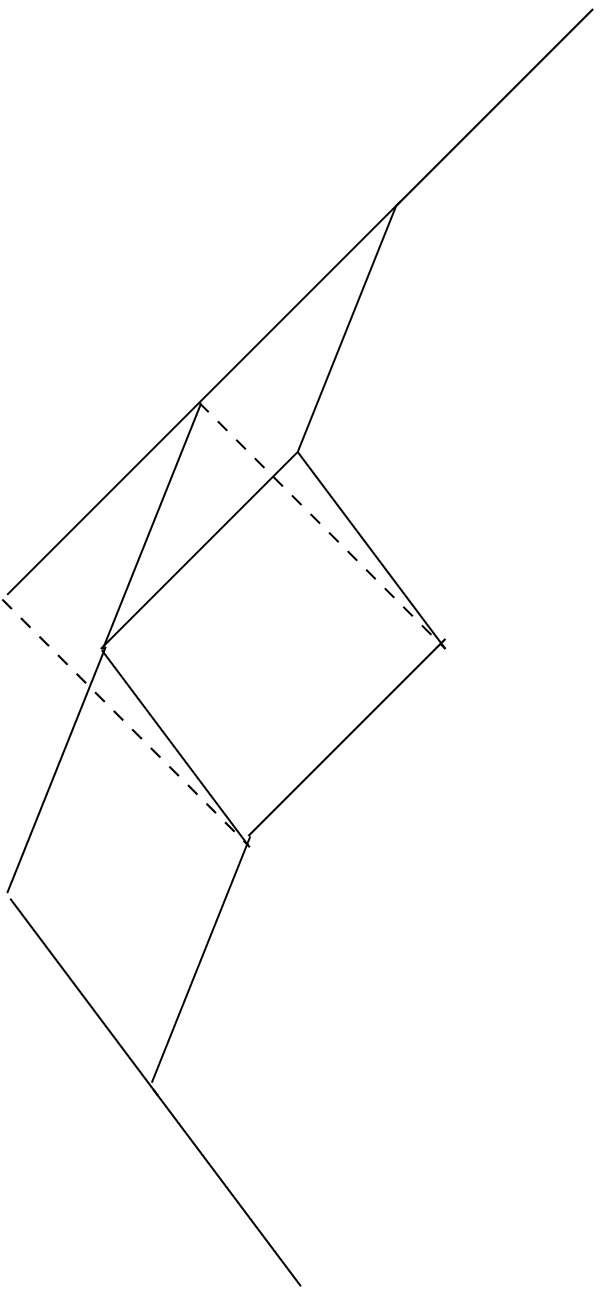}

$n=7$ \includegraphics[width=3\ccm]{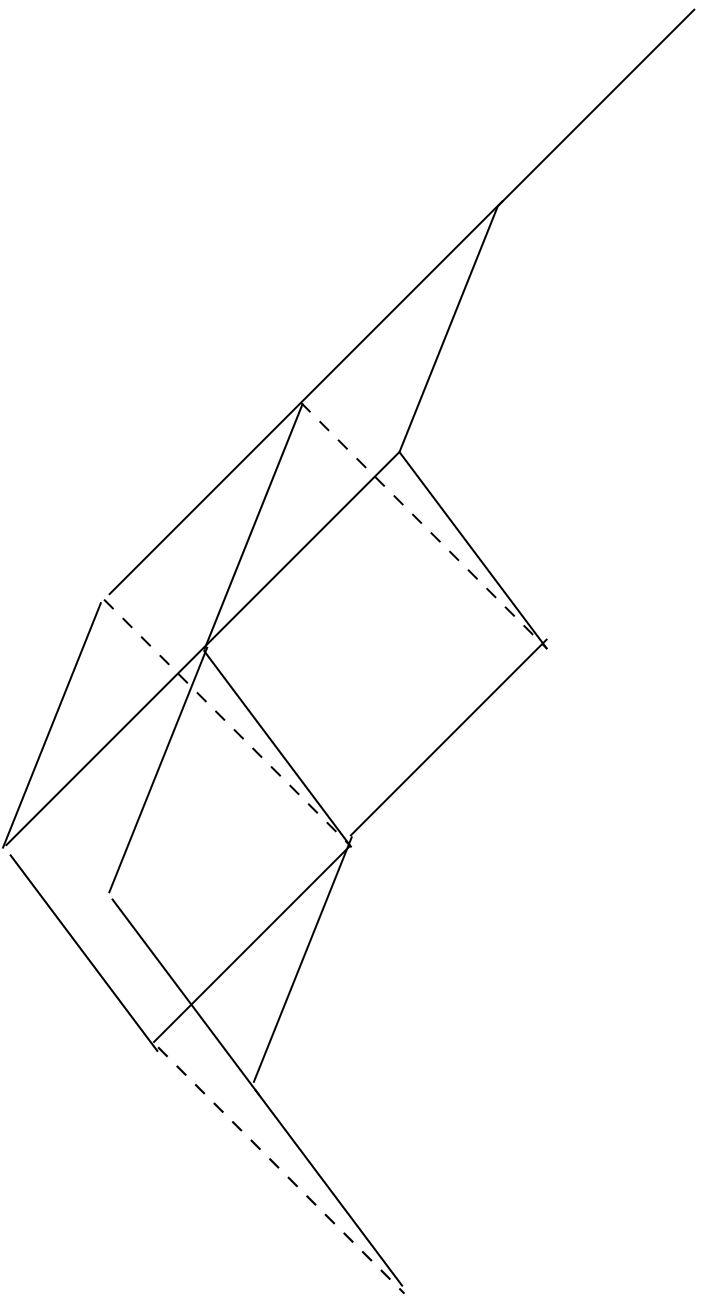}
$n=8$ \includegraphics[width=3.3\ccm]{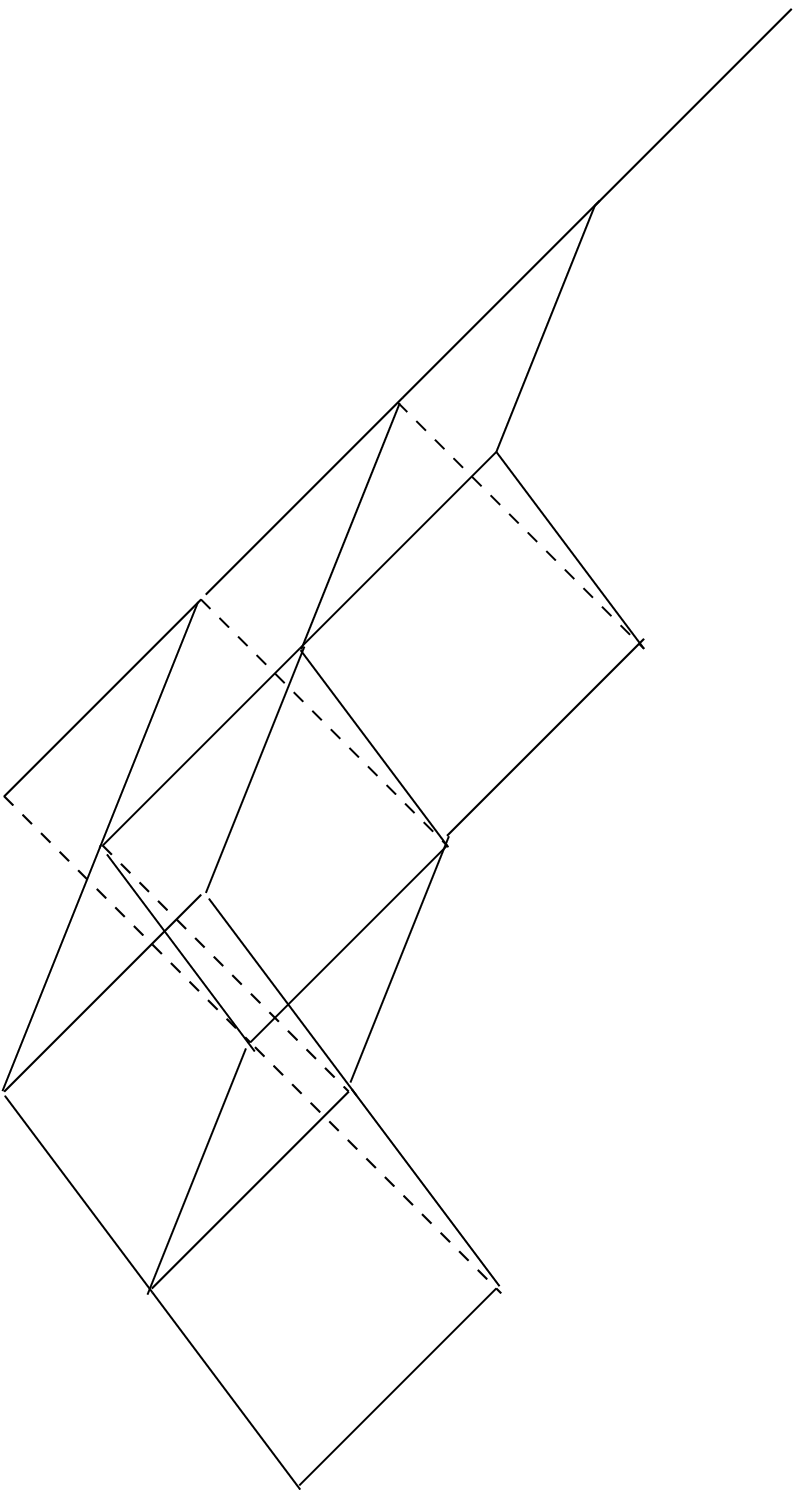}
$n=9$ \includegraphics[width=3.63\ccm]{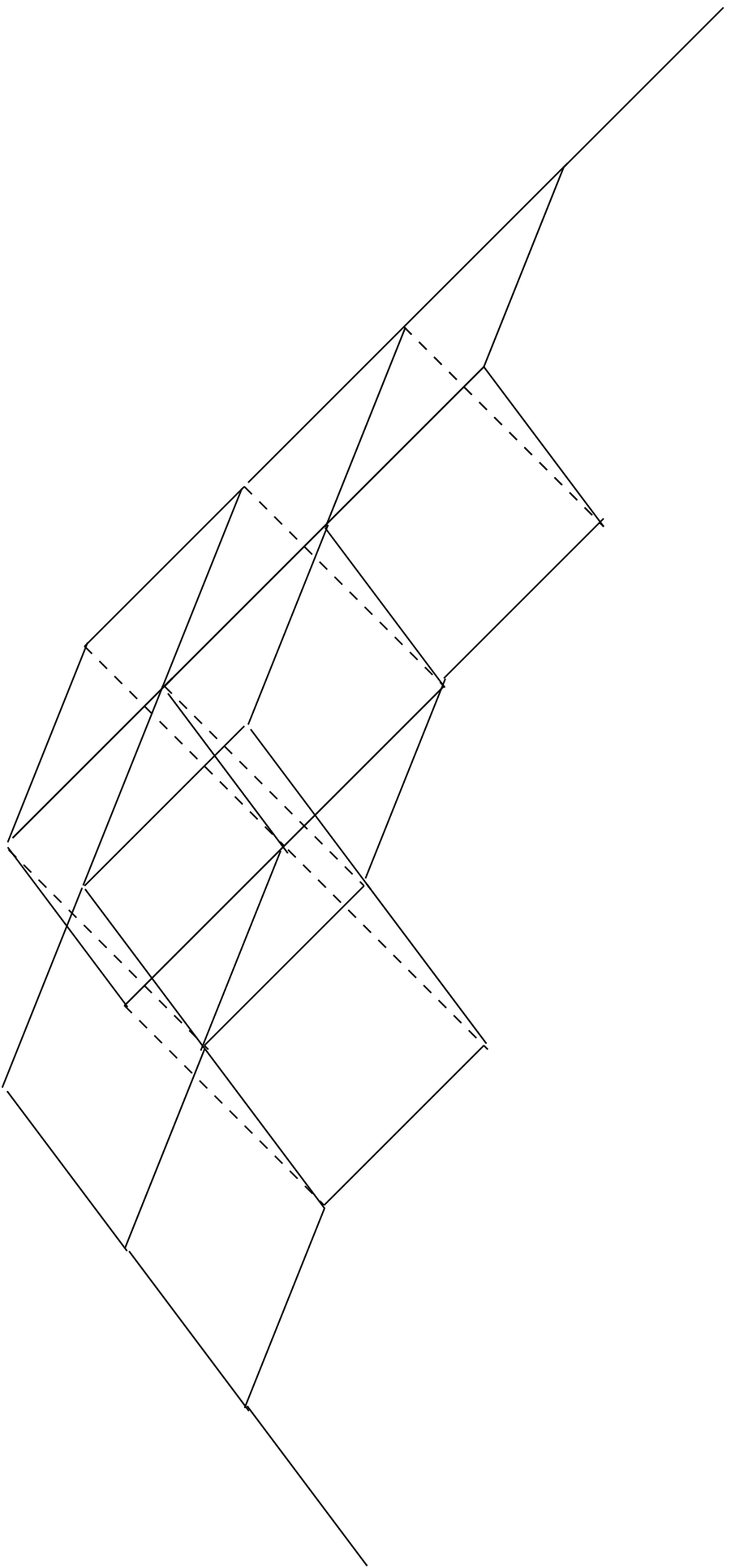}

\caption{Hasse diagram for  inclusion poset 
  for ideals $\Pl_n a^{\mm} \Pl_n$ in case $l=3$, for $n=4,5,...,9$. 
Vertex $\mm$ corresponds to ideal  $\Pl_n a^{\mm} \Pl_n$.
\label{fig:eventriPlattice12}}
\end{figure}


\begin{figure}
\[
\includegraphics[width=3.5cm]{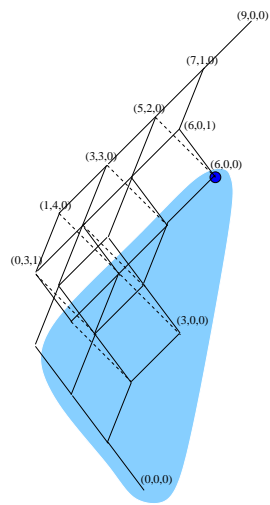} 
\]
\caption{Schematic of the set $h^l_n$, i.e. $\gammal$
  with the poset ideal generated by $(n-l,0,0,..)$ removed. (This case is $h^3_9$.)
  \label{fig:h/h}}
\end{figure}



\newcommand{\hl}{h^l} 

\mdef \label{de:h1}
Note from (\ref{de:gammal}) that 
$\gamma^{l,n-l} \hookrightarrow \gammal$.
As in Fig.\ref{fig:h/h}, 
consider the subset of $\gammal$ 
\beq \label{eq:de:hln}
\hl_n \; := \; \{ \mm\in\gammal \; | \; \mm \not\leq (n-l,0,0,..,0) \}
\eq
We have $\gammal = \hl_n \sqcup \gamma^{l,n-l}$.
Note that
no element of $ \gamma^{l,n-l} $ lies above any element of $\hl_n$ in
the poset $(\gammal,<)$.



\newcommand{\floattwo}{{
\subsection{Combinatorics of the poset $\gammal$}

\begin{figure}
\includegraphics[width=4in]{xfig/fs13o.eps}
\caption{Foreshortenned $sl_3$ view of $\gammal$ in case $l=3$. 
For $n=0,1,2,...$ we include only the part on and below the corresponding
dashed line.  
\label{fig:fs13o}}
\end{figure}


\begin{figure}
\includegraphics[width=4in]{xfig/fs13x.eps}
\caption{Foreshortenned $sl_3$ view of $\gammal$ in case $l=3$,
showing numbers of ideals.  
\label{fig:fs13x}}
\end{figure}

\mdef For later use it will be convenient to have a view of the underlying
$l$-dimensional hypercubical lattice of $\gammal$ from the $(0,0,...,0,1)$
direction. That is, in a projection in which $(0,0,...,0,n)$ and
$(0,0,...,0,n-l)$ coincide. 
(This is what might be called the `Lie theory projection' --- see
later.) 
The poset then takes on the foreshortenned appearance as in
Fig.\ref{fig:fs13o}.  


\mdef
Using the view of $\gammal$ in Fig.\ref{fig:fs13o}
we can compute the number of elements, as in 
Fig.\ref{fig:fs13x}. Looking at the columns of numbers on the left:

The rightmost column of numbers is the level $n$, labelling the 
dashed line. 

The next column gives the number of vertices in each line, ignoring
the fact that the picture is a projection.
The pattern here will be clear, and is elementary to verify:
\[
1,0,1,1,1,1, \;\; 2,1,2,2,2,2, \;\; 3,2,3,3,3,3, \;\; ...
\]

The next column gives the number of vertices in the line $n$ and the
lines $n-3$, $n-6$ and so on. That is, the number of vertices taking
account of the projection.

The final (leftmost) column gives the cumulant sum of vertices up to
this level. 

}}

\newcommand{\floatone}{{


\ignore{{
\newcommand{\edge}{ .\newline --- WAITING FOR UPDATE BELOW HERE --- \newline}
\newcommand{\nend}{ \[ \] MORE TO FOLLOW. \end{document} } 
\newcommand{\mm}{{\mathbf m}} 
\newcommand{\iss}{S_l \times S_m}
\newcommand{\kiss}{k\Smm}
\newcommand{\kissin}{k(S_{i} \times S_{n-2i})}
\newcommand{\Wl}{W^l} 
\newcommand{\Wll}{W}  
\newcommand{\Pl}{P^l} 
\newcommand{\gammal}{\gamma^{l,n}}    
}}


\subsection{Combinatorics of the ideal poset $\gammal$}

See \cite{AMM1} for now.
In \cite{AMM3} the poset $\gammal$ for $l=3$ and fixed $n$ is 
denoted $I(n {\bf e}(1))$, and
\[
C^{(3)}_n \; := \; | I(n {\bf e}(1))  |
\] 
This sequence begins 1,1,2,4,5,7,11,13,17,....


\begin{figure}
\includegraphics[width=4in]{xfig/fs13o.eps}
\caption{Foreshortenned $sl_3$ view of $\gammal$ in case $l=3$. 
For $n=0,1,2,...$ we include only the part on and below the corresponding
dashed line.  
\label{fig:fs13o}}
\end{figure}

\ignore{{

\begin{figure}
\includegraphics[width=4in]{xfig/fs13x.eps}
\caption{Foreshortenned $sl_3$ view of $\gammal$ in case $l=3$,
showing numbers of ideals.  
\label{fig:fs13x}}
\end{figure}

\mdef For later use it will be convenient to have a view of the underlying
$l$-dimensional hypercubical lattice of $\gammal$ from the $(0,0,...,0,1)$
direction. That is, in a projection in which $(0,0,...,0,n)$ and
$(0,0,...,0,n-l)$ coincide. 
(This is what might be called the `Lie theory projection' --- see
later.) 
The poset then takes on the foreshortenned appearance as in
Fig.\ref{fig:fs13o}.  


\mdef
Using the view of $\gammal$ in Fig.\ref{fig:fs13o}
we can compute the number of elements, as in 
Fig.\ref{fig:fs13x}. Looking at the columns of numbers on the left:

The rightmost column of numbers is the level $n$, labelling the 
dashed line. 

The next column gives the number of vertices in each line, ignoring
the fact that the picture is a projection.
The pattern here will be clear, and is elementary to verify:
\[
1,0,1,1,1,1, \;\; 2,1,2,2,2,2, \;\; 3,2,3,3,3,3, \;\; ...
\]

The next column gives the number of vertices in the line $n$ and the
lines $n-3$, $n-6$ and so on. That is, the number of vertices taking
account of the projection.

The final (leftmost) column gives the cumulant sum of vertices up to
this level. 

}}

\section{Bijection with nonisomorphic hollow hexagons}

In structural chemistry, 
the structural classes of polycyclic conjugated hydrocarbons are
generally of necessity enumerated by computer \cite{,Cyvin95}.
One class for which there is a closed form is the nonisomorphic hollow
hexagons \cite{}.
We may define these mathematically as follows.

We start with some alcove geometry, which is not strictly necessary,
but helps to place the problem in a  wider setting useful for later
generalisations. 

\subsection{Basic alcove geometry}

\mdef Consider the Euclidean plane. 
Reflection in a line passing through the origin defines a copy
of the group $S_2$ acting on the plane \cite{Humphreys90}. 
Now add reflection in another
1d subspace that makes an angle $\pi/3$ with the first. 
Together they generate a copy of $S_3$. 
In particular we get a new simple reflection 
line which is the image of the
first in the second. 
Note that these three lines may be considered to subdivide the plane
into various parts, including six open unbounded `faces' or chambers.
 
Adding an affine line at unit distance from the first then
generates a copy of the $A_3^{\tilde{}}$ affine reflection group
in an analogous way. 
The corresponding closed collection of lines gives a triangulation of
the plane 
into `alcoves'.


\mdef Consider a given tiling $\Delta$ of the plane by equilateral
triangles as above.
Associated to $\Delta$ we consider
 the following sets: the set $\Delta_0$ of vertices; 
$\Delta_1$ of lines; and $\Delta_2$ of faces. 
A tile is the closure of a face.

Consider also the geometric dual tiling by regular hexagonal tiles. 
In each case consider the set of partial tilings --- 
the power set of the set of tiles. 
A partial tiling is bounded if
it lies within a finite interval of the plane.
Hereafter we only consider bounded partial tilings.

\mdef
The triangular convex hull of a set of points in $\Delta$ is the
collection of all points, open lines and open faces that intersect the
ordinary geometric convex hull of this set.  
A set of points is minimal if it is minimal among sets with its convex
hull. A minimal set of points is {\em good} if its convex hull and
triangular convex hull coincide. 


\begin{figure}
\[
\includegraphics[width=2.5in]{xfig/hydro1.eps}
\]
\caption{All hollow hexagons up to rank 14
and their `duals'. 
\label{fig:hydro1}}
\end{figure}


\mdef
{\em A good set $P$ contains 6, 4, 3, 2 or 1 points.}
(This claim is obvious, but we write out a proof just to `warm up the
engine'.) 

\proof 
The convex hull is a convex polygon, with the points as
vertices. Thus the points may be ordered clockwise from some root point.
By the good condition two successive points lie on a good line, in the
obvious sense.
By the minimal condition no three points are colinear, so make an
angle of at most $2\pi/3$. \qed

\mdef We have shown that every good polygon is a hexagon or some
degeneration thereof. 

\subsection{Hollow hexagons}

\mdef
Consider a 
partial tiling by hexagonal tiles, $H$, that is
connected and simply connected.
Note that such a partial tiling has an external boundary. 
We consider that the `dual' $d(H)$   
consists only of the dual points of the faces and the 
{\em internal} edges and vertices --- thus possibly only points and lines, not
necessarily complete triangles.
The are some examples in Fig.\ref{fig:hydro1}.


\begin{figure}
\includegraphics[width=2.53in]{xfig/hydro2.eps}
\includegraphics[width=1.92013in]{xfig/hydro4.eps}
\caption{ (a) Intersection of six half-planes making a 
(degenerate) hexagon; and two corresponding triangles. 
Here $t=6$ and $(a,b,c) = (3,3,0)$ so $h=30$ and $\psi=(0,3,0)$.
(b) Intersection of two triangles making a 
(maximally degenerate) hexagon.
Here $t=6$ and $(a,b,c)=(6,0,0)$ so $h=30$ and  $\psi=(6,0,0)$. 
\label{fig:hydro2}}
\end{figure}

 
A {\em hollow hexagon} is a partial tiling of the plane by regular
hexagonal tiles whose dual is a good hexagon. That is, the dual
 is the  convex hull of (at most) six points. 

\mdef
Note that a hollow hexagon $H$ is not in general strictly a hexagon
(its sides, as a hexagon, are zig-zag).
However its dual hexagon $d(H)$ is a hexagon 
(in general irregular, and possibly degenerate). 

\mdef
In practice one considers only the boundary of a hollow hexagon (hence
`hollow'), since this represents the chain of carbons in a molecule.
The interior is only used for the mathematical construction above. 

The rank of a hollow hexagon $H$ is the number of vertices in the
boundary. We have
\[
h(H) \; = \; 6 + 2c_{d(H)}
\]
where $c_{d(H)}$ is the circumference of the dual hexagon.
For example, the dual hexagon in one of the examples in
Fig.\ref{fig:hydro1} is simply a chain of length 2 edges. The circumference
in this case is 4 since, viewing the chain as a degenerate hexagon,
one must start at some point, go around the `hexagon' and return to
the start.

\mdef Two hollow hexagons are isomorphic if they are related by
translation, rotation and/or reflection in the plane
(or the latter can be considered as an orientation reversing 
`flip' of the plane).

Let $S(h)$ denote the set and let 
$N(h)$ denote the number of isomorphism classes of hollow hexagons
of rank $h$. 
Let 
\[
S \; := \; \bigcup_h S(h) 
\]
and write $d(S)$ for the set of duals.

It will be clear that the first non-vanishing term in sequence $N(h)$ is
$N(6)$ and that only even terms are possible. It turns out that it is
natural to consider the sequence as two alternating subsequences.
For example, for $n=0,1,2,3,...$ one subsequence starts
\[
N_n \; =  \; N(6+4n) \; = \; 1,1,2,4,5,7,11,13,...
\]
It is shown in \cite{AMM1} that this is the same sequence as the
number of ideals of 3-\el\ monoids. We now give a bijection
between the underlying sets.

\mdef
Fig.\ref{fig:hydro1} has a few examples.

\subsection{Enveloping triangles}

\mdef
Given a hollow hexagon $H$ there is a collection of six (not
necessarily distinct) affine lines that are defined by 
the edges of the dual hexagon $d(H)$.   
The hexagon $d(H)$ is the intersection of the `upper' half-planes  of
these lines. See Fig.\ref{fig:hydro2}.
The intersection of half planes 
for either alternating subsequence of three of these six lines defines
a triangle. Thus $T_u(H)$ is the up-pointing triangle and $T_d(H)$ the
down-pointing triangle of $H$. 

Given a triangle $T$ write $t=|T|$ for the linear size. 

\mdef Remark: The smaller of $T_u$, $T_d$ is the smallest triangle
containing $d(H)$. 

\mdef
Let $T(H)$ denote the smaller of $T_u$, $T_d$
(or either if $|T_u|=|T_d|$). 
The complement of $d(H)$
in $T(H)$ consists of three (possibly vanishing) triangles.
Write 
$$
\phi(H) = \; (a,b,c)
$$ 
for their linear sizes, in non-increasing order. 
This can be thought of as a three-row Young diagram.
Write $\Lambda^3$ for the set of all three-row Young diagrams.
For $t=|T(H)|$ write
$$
\phi'(H) = \; (t;a,b,c)
$$ 
Note that the map $\phi' : S \rightarrow \N_0 \times \Lambda^3$
is injective. 


\mdef
Note that for $H$ such that  $|T(H)|=t$ we have 
\beq \label{eq:abct}
a+ b+ c \; \leq \; t
\eq
Formally a more aggressive truncation of a single triangle is possible
(up to $a+b+c = 3t/2$), but it is easy to
see that the resultant hexagon has a smaller containing triangle.

Write $\Lambda^3_{\leq t}$ 
for the set of triples/Young-diagrams obeying this bound;
and $\Lambda^3_{ t}$ 
for the set of triples/Young-diagrams saturating the bound
--- the set of integer partitions of $t$ into at most 3 parts.

{\mpr{
For given triangle size $|T|=t$ the triple 
$\lambda = (a,b,c) \in \Lambda^3_{\leq t}$ 
defines a dual hexagon and hence a hollow hexagon $H = w_t(\lambda)$ up
to isomorphism. 
Indeed this inverts $\phi$ and hence 
is an injective map from $\Lambda^3_{\leq t}$ to the set $S$ of
hollow hexagons. 

The length of the boundary of the dual hexagon $d(H)$ 
for $H=w_t(\lambda)$ is 
\[
c_d(H) = 3t-(a+b+c)
\]
(recall the length of boundary of the hollow hexagon is then $h= 6+2c_d(H)$).
}}

It follows that the hollow hexagons with given $t=|T(H)|$ fall into
subsets with boundary $h= 6+4t,6+4t+2,6+4t+4, ...$.
See Fig.\ref{fig:hydro5} for example. 

\begin{figure}
\includegraphics[width=6in]{xfig/hydro5.eps}
\caption{ The partition of the set $\Delta_3$ of dual hexagons into 
subsets $\Delta^3_3$ (in which each dual hexagon has $c_{d(H)}=6$
and hence each corresponding hollow hexagon has $h=18$), 
$\Delta^2_3$, $\Delta^1_3$   and $\Delta^0_3$. 
Each dual hexagon is labelled by  its $(a,b,c)$. 
\label{fig:hydro5}}
\end{figure}


\mdef
Let $\Delta_n$ denote the set of hexagons with $|T(H)|=n$;
and $\Delta_n^{a+b+c}$ the subset with $a+b+c$ as indicated.

Note we have
\[
S(h=6+4n) = \bigcup_{i=0,1,2,...} \Delta_{n-i}^{n-3i}
\]
For example
\[
S(30) = \Delta_6^6  \; \cup \; \Delta_5^{3} \; \cup \; \Delta_4^{0}
\]
Note that for $m=n,n-1,...,0$ we have
\[
\Delta_n^m \cong \Lambda^3_m
\]
In other words we have
\[
S(h) \cong \bigcup_{i=0,1,2,...} \Lambda_{n-3i}^{3}
\]

\ignore{{
\mdef 
Define 
\[
S_n(h) = S(h) \cap \Delta_n 
\]
Thus
\[
\Delta_n = S_n(6+4n) \cup S_n(6+4n+4) \cup ...
\]
and
\[
S(h=6+4n) = S_{n}(h) \cup S_{n-1}(h) \cup ...
\]
Note we have
\[
\Delta_n^{n} \cong \Lambda^3_n
\]
}}


Thus we have a bijection
\beq \label{eq:bijop}
\phi : S(h) \rightarrow \bigcup_{i=0,1,...}  \Lambda^3_{n-3i}
\eq
For example at rank $h=18$, that is $n=3$, we have 
$(3,0,0), (2,1,0), (1,1,1)$ from $\Delta_3$ (as in the Figure);
and $(0,0,0)$ from $\Delta_2$.

\mdef There is a map  $\psi :S(h) \rightarrow \N_0^3$ defined on $H$ with
triple $\phi(H) = (a,b,c)$ by 
$$
\psi (H) = (a-b,b-c,c) .
$$

\newcommand{\gammax}{\gamma}

{\mpr{ 
The map $\psi: S(h=6+4n) \rightarrow \gammax^{3,n}$ is a bijection.
}}
\proof 
It is elementary to check well-definedness. 
We already proved the cardinalities agree, so it is enough to
prove injectivity. $\psi$ is clearly injective, so it is enough to
note the injectivity of $\phi$ as in (\ref{eq:bijop}). 
\qed

\medskip

\mdef Remark:
Note that the integer sequence of the set sizes $I_n = |\Lambda^3_n |$, 
$n=0,1,2,...$ is well known. It starts 
\[
1,1,2,3,4,5,7,8,10,...
\]
We have shown that our main sequence $N_n$ is simply the cummulant sum in
steps of 3. 

NB, One can compute the generating fuction directly from this observation.


\section{On higher $d$} 
And on lower $d$! ...


}}

\medskip

We discuss combinatorics of $\gammal$ in
\cite{AMM1}.


\subsection{Ideals generated by the $a^{\mm}$ elements}
\newcommand{\mPl}{\mP^{l}} 
\newcommand{\Ina}{I^{<\mm}} 


\mdef \label{de:Pnm}
Note that 
\[
\hashi (a^{\mm} ) = \mm
\]
Fix $l$.
Define  the subset $\mP^{\mm}_n \subset \mPl_n$ 
by 
\beq \label{eq:dePmn}
\mP^{\mm}_n =  \{ p \in \mPl_n | \hashi (p) = \mm \}
\eq


\mdef Fix $l$ and $n$  and 
define partition $W$ in $\mP^{l}_n$ by 
\beq \label{eq:W}
W \; = \; \Wl \; := \; (w w^\flip ) \otimes 1_{n-l}
\eq
NB $\Wl \; \in \; S_n a^{(n-l,0,0,...,0)} S_n $.
Example:
\[
W^3 = \raisebox{-.21in}{\includegraphics[width=1in]{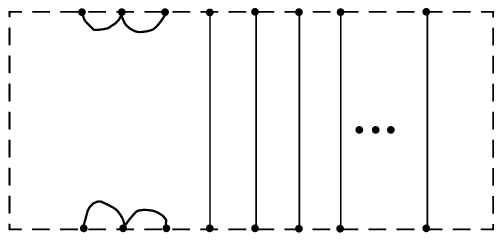}}
\]
We will also use an `idempotent version'
$W_b^l \; = \;\AAA{l\;l\!+\!1} \;W^l \; \AAA{l\;l\!+\!1}$:
\[
W_b^3 = \raisebox{-.21in}{\includegraphics[width=1in]{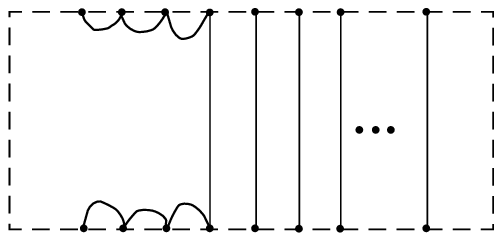}}
\]


{\mlem{ \label{lem:genX}  {\rm \cite{Tanabe97}} $\;$
For $l\in \N$, 
$P_n^l = \langle S_n, \AAA{12},W^l  \rangle $. 
\qed
}}




{\mlem{ \label{lem:1}
(I) The element $a^{\mm} \in \Pl_n$ lies in the ideal 
$\Pl_n a^{\mm'} \Pl_n$ if and only if $\mm \leq \mm'$ in $\gammal$.
(II) For $p,p' \in \mP^{l}_n$ we have $\hashi(pp') \leq \hashi(p)$.    
}}
\proof 
(I) (`if' part): 
If $\mm-\mm' = v_{ij}$ then, up to a permutation,
$a^\mm$ can be obtained from $a^{\mm'}$ by
joining a co-$i$ and a co-$j$ 
propagating part together using some $\AAA{gh}$:
$a^{\mm'} \leadsto a^{\mm'} \AAA{gh}$,
or cutting a propagating `line' using $W^l$. 
\\
(`only if' part): 
Consider $p \in \mPl_n$ and the change 
$\hash^-(p)  \leadsto \hash^-( ap)$ for 
$a \in S_n \{ 1,\AAA{12},W^l \} S_n$. The vector can only be changed 
by such an $a$ as follows. (1) by
combining propagating parts
(for example $a \in \{ \AAA{ij} \}_{ij}$; or 
$a=W^l$ and $p$  is 
as in
the following case with $l=3$:
$W^l p = $ \rb{-.21in}{\ig[width=.71in]{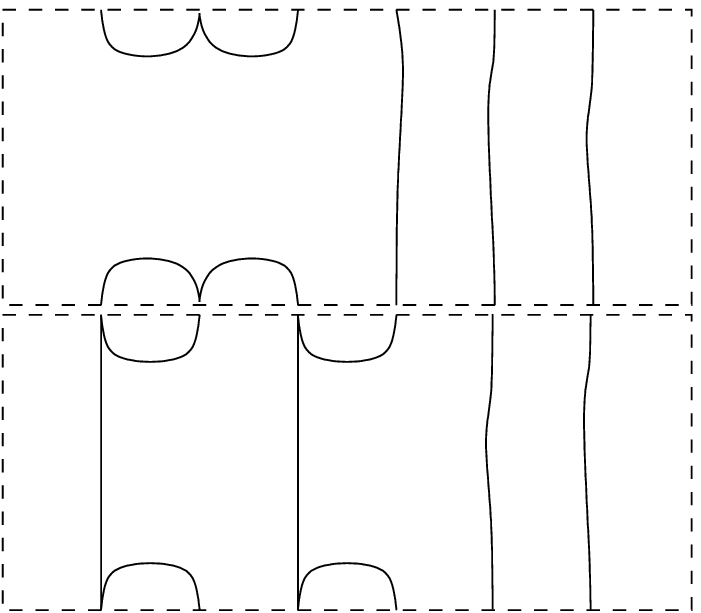}});
or (2) cutting the propagating line in a
$co-l$ part. 
But $p \leadsto ap$ (or $pa$) for {\em any} $a$ is a sequence 
of such `local' changes by Lemma~\ref{lem:genX}.
(II) follows by the same argument.
\qed



{\mlem{ \label{lem:idbas1} \label{lem:pinid-mega}
(I) The set $\bigsqcup_{\mm' \leq \mm} \mP^{\mm'}_n$ is a basis for the ideal
$\Pl_n a^\mm \Pl_n$.
(II) The partitions $p$ in the ideal $\Pl_n a^{\mm} \Pl_n$ and not in 
any  ideal $\Pl_n a^{\mm'} \Pl_n$ with $\mm' < \mm$
are precisely the subset $\mP^{\mm}_n $.  
}}
\proof This follows from Lem.\ref{lem:1}. Note that elements of 
$\mP^{\mm}_n$ are elements of $S_n a^{\mm} S_n$ and elements obtained
from these by binding a non-propagating part to a propagating one.
\qed


{\mlem{ \label{lem:12} Fix $l$. \\
(1) Every element $p$ of $\mP^{l}_n$ containing a part $p_i$ with 
$\cora(p_i) > l$ or $\rang(p_i) >l$
lies in the ideal $P^l_n W P^l_n$.

\noindent
(2) Every element of $\mP^{l}_n$ containing a non-propagating part 
lies in the ideal $\Pl_n W \Pl_n$.

\noindent
(3)

Fix $n$. For any $\mm$, if   
$r_\mm = n$ 
then  
$ \; a^\mm \; 
      \not\in \Pl_n W \Pl_n$.
}}
\proof (1) Suppose  $\cora(p_i) > l$  or $\rang(p_i) >l$.
Then there is a partition 
$p' \in \mPl_n$
differing from $p$ only in that $l$ elements of $p_i$ are in an
isolated non-propagating part. Partition $p$ evidently lies in the ideal
generated by $p'$. Now use (2).

(2) A non-propagating part $p_i$ 
in $p \in \mPl_n$ necessarily has order at least $l$.
Then there is a partition $W'$ group-conjugate  to $W$ such that one
of its non-propagating parts exactly meets a subset of $p_i$ in composition,
whereupon
$pW' = p$ or $W'p = p$. (If the order is exactly $l$ the argument is
slightly modified.)

(3) By Lem.\ref{lem:1}, noting that
$W \in S_n a^{(n-l,0,0,...,0)} S_n$.
\qed





\mdef \label{de:PPI}
For given $n,l$ and $\mm\in\gammal$
define the ideal
\[
\Ina = \sum_{\mm' < \mm} \Pl_n a^{\mm'} \Pl_n    
\]


{\mlem{ \label{le:keymm}
If $\mm' \not\geq \mm \in \gammal$ 
then $a^{\mm'} \Pl_n a^{\mm} \subset I^{< \mm}$.
}}
\proof Note that $a^{\mm'} \Pl_n a^{\mm}$ has a basis of partitions.
By \ref{lem:pinid-mega} every partition lies in a unique
{highest} ideal
of form $\Pl_n a^{\mm''} \Pl_n$.
In particular $p \in a^{\mm'} \Pl_n a^{\mm}$ lies in or below 
$\Pl_n a^{\mm} \Pl_n$, but since $\mm' \not\geq \mm$ it must be below.
\qed



\section{Representation theory generalities} \label{ss:gen}
Lemma~\ref{le:keymm} means that $\Pl_n$ has what we call the `core' property.
From this property many representation theoretic properties follow
quite generally. Here we collect the general arguments.

\subsection{Preliminaries: Green's idempotent reciprocity}

\newcommand{\defeq}{:=} 


{\mth{ \label{th:green}
    {\rm [Green localisation theorem, \cite[\S6.2]{Green80}]}
    Let $k$ be a field, let $A$ be a $k$-algebra, and $e \in A$ idempotent.
    Let $\Lambda(A), \Lambda(eAe)$ and $\Lambda(A/AeA)$ be index sets
    for classes of simple modules of the indicated algebras.
    Then there is a bijection
\[
\Lambda(A) \stackrel{\sim}{\rightarrow} \Lambda(eAe) \sqcup \Lambda(A/AeA)
\]
\qed
    }}

Fix an index set $\Lambda(A)$, and a set
$\{ M_\lambda : \lambda\in \Lambda(A) \}$
that is a complete set of simples up to isomorphism.
Let us write $\Lambda_e(A)$ for 
the subset of $\Lambda(A)$ such that
$\lambda\in \Lambda_e(A)$ implies $eM_\lambda \neq 0$
\cite{Assem06,
  Green80}.
The map $M \mapsto eM$ is a bijection from
$\{ M_\lambda : \lambda\in \Lambda_e(A) \}$
to a complete set of simples for $eAe$. 
Thus we may take $\Lambda(eAe) = \Lambda_e(A)$.
Meanwhile 
$\{ M_\lambda : \lambda\in \Lambda(A)\setminus\Lambda_e(A) \}$
is a complete set for $A/AeA$. 
This gives a natural identification
\beq \label{eq:Lambde}
\Lambda(A) = \Lambda(eAe) \sqcup \Lambda(A/AeA)
\eq


\newcommand{\eW}{e}
\newcommand{\AP}{A}
\newcommandabbreviation{\mmod}{\mbox{mod}}{mod}
\newcommand{\CoM}[1]{[ #1 ]}  

\mdef
\redx{Let $k$ be a field.}
Given a finite-dimensional $k$-algebra $A$, a simple module $L$ and a
module $M$, then $\CoM{M:L}$ denotes the composition multiplicity of $L$
in $M$. 

\medskip

\mdef \label{de:Gfunctorx}
For $A$ an algebra then $A\!-\!\mmod$ denotes its category of left-modules.
Given an idempotent $\eW$ in an algebra $\AP$  
we define, as usual
\cite{Martin96,Green80},
a functor
\[
G_\eW : \eW \AP \eW - \mmod \;\; \rightarrow \; \AP - \mmod
\]
by $ G_\eW(M) = \AP \eW \otimes_{\eW \AP \eW} M$.
Define functors
\begin{equation*}
\begin{split}
F_\eW  &: \AP - \mmod \;\; \rightarrow \; \eW \AP \eW - \mmod\\
L_\eW  &:  \mmod- \AP \;\; \rightarrow \;  \mmod-\eW \AP \eW
\end{split}
\end{equation*}
by $F_\eW(N) = \eW N$ and $L_\eW(N) =  N\eW$. We have the following standard properties
(see e.g. \cite{Green80,Assem06,Martin0915}). 

{\mth{ \label{th:Gew}
Let algebra $A$ and idempotent $e \in A$ be as above.
\\    
(I) Functor $G_\eW$ is left-adjoint to $F_\eW$; and a right-inverse to $F_\eW$.
\\
(II) Functor $G_\eW$ is right-exact; and $F_\eW$ is exact.
\\
(III) Functor $G_\eW$ preserves projectivity and indecomposability.
\\
(IV) Functor $G_\eW$ preserves simple head.
\\
(V) For $L$ a simple module $eL$ is a simple $eAe$-module or 0 and if
$eL\neq 0$ then 
\[
\CoM{ eM : eL } \; = \;  \CoM{ M:L }
\]
\qed
}}


\subsection{Algebras with the core property}
\newcommand{\lambdao}{\omega}

\newcommand{\lambdaa}{\alpha}
\newcommand{\mua}{\beta}

\ignore{{
\mdef \label{de:core} 
Let $A$ be a $k$-algebra, with $k$ a commutative ring.
Let $\gamma=(\gamma,<)$ be a finite poset and
$e_\lambdaa \in A$ an idempotent for each $\lambdaa\in \gamma$.
Let
$$
I^{<\lambdaa} = \sum_{\mua < \lambdaa} A e_\mua A
$$
(so $I^{<\lambdaa} = 0$ if $\lambdaa$ is a lowest element)
and
$$
A^\lambdaa = A/I^{<\lambdaa} .
$$
The poset $\gamma$ together with the map $e_- : \gamma \rightarrow A$
is called a {\em core} for $A$ if 
(AI): $e_\lambdaa A e_{\lambdaa'} \subseteq I^{<\lambdaa}$ for all
$\lambdaa' \not\geq \lambdaa$;
(AII): $1 \in \gamma$ (that is, with $e_1 = 1$).
}}

\mdef \label{de:core}
{
  Let  $k$ be a commutative ring.
Let $A$ be a unital $k$-algebra, with unit $1$.
Let $\gamma=(\gamma,<)$ be a finite poset and
$e_\lambdaa \in A$ an idempotent for each $\lambdaa\in \gamma$.
Let
$$
I^{<\lambdaa} = \sum_{\mua < \lambdaa} A e_\mua A
$$
(so $I^{<\lambdaa} = 0$ if $\lambdaa$ is a lowest element)
and (for later use)
$$
A^\lambdaa = A/I^{<\lambdaa} .
$$
The poset $\gamma$ together with the map $e_- : \gamma \rightarrow A$
is called a {\em core} for $A$ if 
(AI): $e_\lambdaa A e_{\lambdaa'} \subseteq I^{<\lambdaa}$ for all
$\lambdaa' \not\geq \lambdaa$;
(AII): $1 $ lies in the image of $e_-$.
}

\ignore{{
\mdef
\red{Delete this in a minute!}
Note the following. The map $e_-$ in a core $(\gamma,e_-)$ is not
necessarily injective.
The poset $\gamma$ has a unique maximal element $\top$, 
and $e_{\top} =1$.
\red{
--- No! does not follow!}
}}


{\mth{ \label{th:index1core}
Suppose  $\delta $ invertible in $\kk$.
Then the pair  $(\gammal , a^- )$ gives a core for $\Pl_n$.
}}

\proof Comparing (\ref{de:core})   
with (\ref{le:keymm}) we see
that $(\gammal , a^- )$ is a core for $\Pl_n$
up to renormalisation of the $a^-$ elements as idempotents.
\qed

\medskip

\mdef Suppose now that $k$ in (\ref{de:core}) is a field and $A$ is
finite-dimensional
with simple index set $\Lambda(A)$. 
Note the natural inclusion
$\Lambda(A^\lambdaa ) \hookrightarrow \Lambda(A)$
giving the (classes of) simple modules $M$ of $A$ such that
$e_\mua M = 0 $ for $\mua < \lambdaa$.

\label{de:LeA}

By $e_\lambdaa A^\lambdaa$ we understand $e_\lambdaa$ to act on
$A^\lambdaa$ in the natural way, i.e. as $e_\lambdaa + I^{<\lambdaa}$.
By Green's theorem (\ref{th:green}) and the construction then we take 
$\Lambda(e_\lambdaa A^\lambdaa e_\lambdaa) \hookrightarrow \Lambda(A)$
(via $\Lambda(e_\lambdaa A^\lambdaa e_\lambda) \hookrightarrow \Lambda(A^\lambdaa)
\hookrightarrow \Lambda(A)$)
to index classes of simple
modules of $A$ such that
$e_\mua M = 0 $ if $\mua < \lambdaa$ and
$e_\lambdaa M \neq 0$.



\mdef \mulem \label{le:ccomp}
Let $A$ be an algebra
over a field
with simple index set $\Lambda(A)$. 
If algebra $A$ has a core $(\gamma, e_-)$ then 
$\Lambda(A) =
   \cup_{\lambdaa\in\gamma} \Lambda(e_\lambdaa A^\lambdaa e_\lambdaa) $.

\proof
Consider the classes indexed by
$ \Lambda(e_\lambdaa A^\lambdaa e_\lambdaa) $
as in (\ref{de:LeA}).
There may be other simples with $e_\lambdaa M \neq 0$, but they have
$e_\mua M \neq 0$ for some $\mua < \lambdaa$,
and so are `counted' in some lower
$ \Lambda(e_\mua A^\mua e_\mua) $
--- noting that minimal elements $\omega \in \gamma$ finally exhaust
simples with $e_\omega M \neq 0$.
Thus the union includes all  simples with $e_\lambdaa M \neq 0$. 

But now since $1 \in \gamma$ (by axiom AII) the union includes all simples in
$\Lambda(A)$ where $1M \neq 0$, which is all.
\qed


\mdef
Remark/caveat: Note that every algebra has a `core'
for every poset with a unique minimal element
(call it $\bot$), simply
by the constant mapping $e_\lambdaa =1 $.
But all but one
$\Lambda(e_\lambdaa A^\lambdaa e_\lambdaa)$
is empty in this case
(and our
Theorem~\ref{th:genLam}
below
is trivial).
Confer 
\cite{ErdmannSaenz03,CPS,DlabRingel,GrahamLehrer96}
for certain `tighter' axiomatizations
(and indeed cf. for example
\cite{green2004standard,boltje2012twisted,enyang2017cellular}
and references therein for interesting related axiomatizations).
\redxx{Trying this caveat here now. Maybe move again? Yes, we are not
  yet over a field so this is too soon.}


In order for the core formalism to have significant utility, we will need
  something further like BH reciprocity,
as we discuss shortly.



{\mth{ \label{th:genLam}
    Let $A$ be an algebra
over a field
    with simple index set $\Lambda(A)$. 
If algebra $A$ has  a core $(\gamma,e_- )$  then
\beq \label{eq:partLA}
\Lambda(A) = 
\sqcup_{\lambdaa\in\gamma} \Lambda(e_\lambdaa A^\lambdaa e_\lambdaa )
\eq
}}
\proof
\ignore{{It follows from the theorem  and (\ref{eq:Lambde}) that
$\Lambda(eAe)$
(possibly empty) includes in
$\Lambda(A)$. We need to show disjointness and completeness of
$
\sqcup_\lambda \Lambda(e_\lambda A^\lambda e_\lambda )
$. 

1. Disjointness:
}}
Noting Lem.\ref{le:ccomp} it remains to prove disjointness. 
For $\lambdaa\in \gamma$ let
$\Lambda_\lambdaa \defeq \Lambda(e_\lambdaa A^\lambdaa e_\lambdaa)$.
\redxx{ --- that is already defined above! ---
maybe delete the one above and keep this one??? ---}
Recall the natural embedding 
$\Lambda_\lambdaa \hookrightarrow
    \Lambda(A^\lambdaa ) \hookrightarrow \Lambda(A)$.
Simple modules $M$ of $A$ coming from the subset $\Lambda_\lambdaa$ (if any) obey
$e_\lambdaa M \neq 0$ ($A e_\lambdaa M = M$)
\redx{by (\ref{th:green})/(\ref{de:LeA}).
}
But now consider $e_{\lambdaa'} M$. 
Either $\lambdaa' \not> \lambdaa$ and so $e_{\lambdaa'} M= 0$ by axiom (AI),
so $M \not\cong M_\mu$, $ \mu \in \Lambda_{\lambdaa'}$
(for which $A e_{\lambdaa'} M_\mu = M_\mu$);
or $\lambdaa' > \lambdaa$ 
so $\lambdaa \not > \lambdaa'$ 
and $e_\lambdaa N = 0$
for all $N=N_\mu$, $\mu\in\Lambda_{\lambdaa'}$,
while $A e_{\lambdaa'} N = N$ so $M \not \cong N$.
\qed 


\ignore{{
2. Completeness:
We may work by induction on the size of $\gamma$.
The base case $\gamma= \{ 1 \}$ is trivial since it asserts only that $\Lambda(A) =
\Lambda(A)$.
For the inductive step we proceed as follows.

DO WE?

Consider $\lambdao$ a lowest element in
$\gamma$.
We have $\Lambda(A) = \Lambda_\lambdao \sqcup \Lambda(A/A e_\lambdao A)$
by Green's theorem.
We claim $A/A e_\lambdao A$ satifies the axioms using the poset
$\gamma\setminus\lambdao$ and the images $\bar{e_\lambda}$ of the
idempotents from $A$.

VERIFY!

OR?...

We claim that for each $\lambda\in \gamma\setminus\omega$ we have
$
\Lambda(e_\lambda A^\lambda e_\lambda)
  = \Lambda( \bar{e_\lambda} \bar{A^\lambda} \bar{e_\lambda})
$.

OR...
  
In $\Lambda(e_\lambda A^\lambda e_\lambda)$
we have (labels for) the simple modules $M$ of $A$ such that
$e_\lambda M \neq 0$ (by the localisation)
and $e_\nu M = 0$ for $\nu < \lambda$ (by the quotient).

CAN WE NOT JUST SEE THAT THIS IS EVERYTHING???...

Thus the partial union
$\sqcup_{\lambda \leq \nu} \Lambda(e_\lambda A^\lambda e_\lambda)$
contains all simples with $e_\lambda M \neq 0$ for $\lambda \leq \nu$.

...Thus summing over the poset we have the simple modules for which
$1 M \neq 0$ --- which is all simple modules.

DONE?

In $\Lambda( \bar{e_\lambda} \bar{A^\lambda} \bar{e_\lambda})$
we have the simples modules of $A$ such that ...

VERIFY!
USE THE aXIOMS!

Let us write $\lambda=1$ for the case $e_\lambda = e_1 = 1$ (by
(AII)).
By (AI) this is necessarily the top element in the poset.
For the corresponding section we have
$\Lambda(e_\lambda A^\lambda e_\lambda ) = \Lambda (1 A^{1} 1)$.

We have an ideal filtration 

...
}}


\newcommand{\ccc}[1]{\lceil #1 \rfloor}
\newcommand{\Deltag}{\beth^{\!\gamma}}  

\subsection{$\Deltag$-modules over a core}


\mdef \label{eq:Gdeltam}
\redxx{Maybe change this so working over field comes later.
  G functor works in integral case...
Or maybe leave it like this and observe the more general case later
just when needed?}
Let $A$ be an algebra over a field with index set $\Lambda(A)$,
and  core $(\gamma, e_-)$  as above.
Let $\ccc{-}:\Lambda(A) \rightarrow \gamma$ be the map taking
$\mu\in \Lambda(A)$ to the part
$\alpha = \ccc{\mu}$
to which it belongs in \eqref{eq:partLA}.
For each $\lambdaa\in\gamma$  
suppose
$\{ S_\mu^\lambdaa \; : \; \mu \in \Lambda_\lambdaa \}$
to
be a corresponding
set of simple modules of $e_\lambdaa A^\lambdaa e_\lambdaa$.

For each of the algebras $e_\lambdaa A^\lambdaa e_\lambdaa$
we have a corresponding $G$-functor
(as in (\ref{de:Gfunctorx})).
In particular this lifts
the set of simple modules to a set of $A^\lambdaa$-modules, and hence
$A$-modules.
Thus for each $\mu \in \Lambda(A)$ we have an $A$-module
\beq \label{eq:Gdelta}
\Deltag_\mu = G_{e_{\alpha}} S_\mu^\lambdaa 
\end{equation}
for $\lambdaa = \ccc{ \mu }$ denoting the appropriate
$\lambdaa$.
Let $\Deltag = \{ \Deltag_\mu \; | \; \mu\in\Lambda(A) \}$.

\mdef \label{lem:unitri0}
\mulem {\it
Let $A$ be an algebra with index set $\Lambda(A)$,
and  core $(\gamma, e_-)$  as above.
\ignore{{
For each $\lambda\in\gamma$
let $\{ S_\mu \; : \; \mu \in \Lambda_\lambda \}$ be a corresponding
set of simple modules of $e_\lambda A^\lambda e_\lambda$.
For each of these algebras we have a corresponding $G$-functor
(as in (\ref{de:Gfunctorx})) lifting
the set of simple modules to a set of $A^\lambda$-modules, and hence
$A$-modules.
Thus for each $\mu \in \Lambda(A)$ we have an $A$-module
\red{
  $$
  \Delta_\mu = G_{e_{| \mu |}} S_\mu
  $$
for $\lambda= | \mu |$ denoting the appropriate
$\lambda$.
}}}
Let $\{ \Deltag_\mu \; : \; \mu \in \Lambda(A) \}$ be as in \eqref{eq:Gdelta}.
(In particular this supposes fixed the various sets
$\{ S^{\lambdaa}_\mu \; : \; \mu \in \Lambda_\lambdaa \}$.)
We have the following:
\\
(I) The modules $\Deltag_\mu$ have simple heads,
\redx{denoted $L_\mu$};
and their heads are exactly a complete set of simple $A$-modules.
\\
(II) 
If $\Deltag_\mu$ has a simple factor $L_\nu$ below the head then
$e_{\ccc{\mu}} L_\nu = 0$.
\ignore{{
this factor comes from higher
\red{or should this be not lower? yes probably. and this is
  essentially vacuous - given the quotient! so (II) is useless!
\\ interestingly it looks like
  --- in the tonal partition case --- these modules do have useful
unitriangular
properties, at least over $\C$, derived by a different route ...
Can they be seen from the present construction somehow???...
maybe not, given that the remark is case specific.}
in the core poset order
(i.e. if $|\mu|$ denotes the $\lambda$ from which $\mu$ comes then
$| \nu | > | \mu |$ \red{no. LOL}). 
...
}}
\\
(III) Suppose $A$ possesses an involutive antiautomorphism $a \mapsto a^\star$
and hence a contravariant duality $M \mapsto M^o$ \cite{Green80}, 
so the socle of $(\Deltag_\mu)^o$ is $L_\mu^o$.
If $e_{\ccc{\mu}} = e_{\ccc{\mu}}^\star$ then $e_{\ccc{\mu}} L_\mu^o \neq 0$,
so either $L_\mu^o$ does not appear in $\Deltag_\mu$ and there is no
module map $\Deltag_\mu \rightarrow (\Deltag_\mu)^o$; 
or $L_\mu^o = L_\mu$,
and there is exactly one module map $\Deltag_\mu \rightarrow (\Deltag_\mu)^o$
up to scalars, with image $L_\mu$.
}

\proof (I) By Theorem~\ref{th:Gew}(IV)
and (the proof of) Theorem~\ref{th:genLam}.
\\
(II) 
Follows from Theorem~\ref{th:Gew}(V).
\\
(III) If $a M \neq 0$ then $a^\star M^o \neq 0$.  
\qed


\mdef\label{uni:lem} \mulem 
\redxx{(II') something about BLOCK unitriangularity of the $\Delta$
  decomposition matrix $[\Delta_\mu : L_x ]_{\mu,x}$
  when ordered so that ... ???}
%
\label{Th:upperun}
{\it
Let $A$ be an algebra 
with core $(\gamma, e_-)$ as above. 
Let $d^\gamma_{\mu x}  := [\Deltag_\mu : L_x ]$.
Then $d^\gamma_{_{\mu \mu}}=1$ and 
if $\mu \neq x$ and 
$d^\gamma_{_{\mu x}}\neq0$ then $ \ccc{\mu}< \ccc{x}$.
}
\ignore{\red{If $e_\lambda A^\lambda e_\lambda$ is semisimple for each
    $\lambda\in\gamma$ then}}

That is to say, 
let $<'$ be any total order  on $\Lambda(A)$ in which $ \ccc{\mu}<\ccc{x}$
in $\gamma$  implies $\mu <' x$.
Then the matrix of 
$D_{A}^\gamma  :=(d^\gamma_{\mu x})$
with respect to  
$<'$
is upper-unitriangular. 
%
\proof
Let $\mu,x\in \Lambda(A)$ with $x\neq \mu$.
First note that  $d^\gamma_{_{\mu \mu}}\neq0$ by 
Lemma~\ref{lem:unitri0}(I).
Lemma~\ref{lem:unitri0}(II) then implies 
$d^\gamma_{_{\mu \mu}}=1$.

It remains to show that if $d^\gamma_{_{\mu x}}\neq0$ then
$ \ccc{\mu}< \ccc{x}$.
If 
$d^\gamma_{_{\mu x}}\neq0$ then $e_{\ccc{x}}\Deltag_\mu\not=0$.
Hence by  axiom 
(AI) of the definition of core we have $\ccc{\mu}\leq \ccc{x}.$
Note that again  by 
Lemma~\ref{lem:unitri0}(II) we have $e_{\ccc{\mu}} L_x = 0$.
Hence we cannot have 
$\ccc{x}=\ccc{\mu}$ otherwise we would get $e_{\ccc{\mu}} L_x \not=0$.
Therefore  $\ccc{\mu}< \ccc{x}$. 
\qed

\ignore{\red{if we assume $\gamma$ is a chain then either $|\mu|=|x|$(which is not possible by the above paragraph) or $|\mu|<|x|$. Hence $\mu<'x.$ }}


\redxx{No QED yet. Just needs tidying up a little more.}

\redxx{Maybe DELETE this?:
  Caveat. 
  $\Deltag$ modules are not generally what one might expect.
  They do not generally have a pivotal property as described
  in (\ref{de:pivot})
(and they are not useful without such a property).}

\mdef \redx{
In general $\Deltag$ modules are  
not particularly useful, for
  computing the Cartan decomposition matrix for example.}
 Below we discuss conditions under which they become useful.

\redxx{Maybe DELETE this?:  There is another thing to check.
Originally we worked very hard to show semisimplicity in
`layers' of the tonal algebra, and this was crucial to the proof. Where
is this work being done now?
This is the same question as \ref{necBH}.}
\redxx{
(An overall problem is that of lack of examples in this section,
as already noted. Could fix by moving sec.2. as discussed.)}


\redxx{Maybe DELETE this?:
 \em Necessary conditions for BH reciprocity.
Here we illustrate caveats on the core property by discussing
conditions under which the $\Deltag$ modules can NOT
be pivotal in the sense of (\ref{de:pivot}).
}

\redxx{
I list some points which I am not sure about, the list may change, possibly the points are all related  
\\1. As we discussed the BH may fail when $e_\lambda=1$
\\2. If we do not have something like eq \ref{eq:sumsquares}}


\subsection{Modular systems and pivotal sets} 

\newcommand{\DeltaD}{{\mathcal D}}  

\mdef \label{de:pivot}
Consider an algebra $A$ over a field
with simple modules $L_\mu$, $\mu\in\Lambda(A)$ and 
indecomposable projective covers  $P_\mu$, $\mu\in\Lambda(A)$.
Suppose we have another set
$\DeltaD$ 
of $A$-modules $\DeltaD_\mu$, $\mu\in\Lambda\supseteq\Lambda(A)$,
that span the Grothendieck group.
Suppose in particular that there is an expression for the characters
of the projective modules $P$ in terms of $\DeltaD$,
with coefficients denoted $(P:\DeltaD_\mu)$.
We say $A$ has the
{\em BH reciprocity} property (as in Brauer--Humphreys) if
\beq \label{eq:BH0}
[\DeltaD_\mu : L_\nu ] =  ( P_\nu : \DeltaD_\mu )
\eq
In particular we say  
$\DeltaD$ is a {\em pivotal} set
for the BH property.

If $\DeltaD$ gives a basis for the Grothendieck group we call this {\em
  strong} BH property.

\redxx{There is a `weak' (i.e. not-necessarily-strong) version
  in which the index set for $\DeltaD$ modules can
  be bigger than for simples. Then this set of modules contains a
  basis for the Grothendieck group (and the well-definedness of
  $  ( P_\nu : \DeltaD_\mu ) $ is
  subtler --- see later).
  And the $\DeltaD$-decomposition
  matrix is not necessarily square.
  We restrict to the square case because we have the highest weight
  category case in mind. See later.
  }

\mdef \label{de:Bms}
A Brauer modular system \cite{Benson95} for a $k$-algebra $A$ is a triple of
commutative rings $(K,K_0,k)$
where $K$ is an integral domain, $K_0$ the field of fractions,
and $k$ the quotient by some maximal ideal,
with the following properties.
Firstly there is an `integral' version $A^K$ of $A$ over $K$,
and $A$ itself  
is obtained by base change, $A=k\otimes_K A^K$, and is split.
Secondly the base change instead to $K_0$,
the `ordinary' version,
is split semisimple.

\newcommand{\most}{BH}

\mdef \label{de:complete}
A {\em \most\ module} for a Brauer modular system
is the  image of a simple module of $K_0 \otimes_K A^K$ 
under an integral lift to $A^K$ followed by the base change to $k$.

There may be many such lifts in general. But a set of \most\ modules
that is the image of a complete set of ordinary simples is called
complete. 

\mdef \label{modBH}
A sufficient condition for the BH property,
\redx{ with pivotal set $\Delta$, }
is that the set $\Delta$ is
a complete set of \most\ modules.
(See e.g. \cite[Prop.1.9.6]{Benson95}.)




\subsection{Modular core property and highest weight categories}
\newcommand{\Sord}{\hat{S}}  
\newcommand{\podular}{podular}
\newcommand{\Sint}{{\mathsf S}} 

\redx{In this section we define the (strong)
  modular core property for algebras and show that
  it implies a highest weight category.
}

\mdef A {\em modular core} is  
\redxx{ either 
  (A) a Brauer modular system for an algebra with core $(\gamma, e_-)$ in
  which there is a complete  set (hence a pivotal set)
  in the sense of \ref{de:complete} that coincides with $\Deltag$ ???;
  or (B) }
 an algebra $A$ and
a triple of rings $(K,K_0,k)$  as in (\ref{de:Bms}) giving a Brauer modular system for
$A$
(thus $A= k \otimes_K A^K$)
together with core data $(\gamma, e_-)$   giving a core for
$A^K$.
\\
\redx{A  {\em strong} modular core is a modular core that is strong modular,
  i.e. the simple modules over $k$ and $K^0$ have the same index set.}

  
\mdef Remark.
In the classical modular theory for finite groups one finds the
integral representations
(that will be reduced mod. the prime defining $k$)
as lattices inside the
simples of the rational case. In diagram algebras suitable integral
representations can generally be constructed directly (by the good
basis properties of such algebras).
Nonetheless their base changes to the rational case are the simples of
that case. In other words it is the index set for simples in the
rational case that labels the `Brauer/Specht modules'.

This means there are some `carts and horses'
(hypotheses and conclusions) 
that must be placed carefully in the
right order!

\mdef
Note that if $(K,K_0,k)$ is such a triple of rings then the core
property over $K$ will base change to the field cases.
\redxx{Correct???}

In particular then
consider  the  
algebra $A^0 = K_0 \otimes_K A^K$.
This is now split semisimple by hypothesis. 
We have
$$
\Lambda(A^0) =
\sqcup_{\lambdaa\in\gamma} \Lambda(e_\lambdaa (A^0)^\lambdaa e_\lambdaa )
$$
by Th.\ref{th:genLam}.
An idempotent subalgebra of a quotient of a semisimple algebra is
semisimple \redxx{ right? }
so
in this scenario $e_\lambdaa (A^0)^\lambdaa e_\lambdaa$ is semisimple.
For each $\lambdaa\in\gamma$ 
let us write
$ \{ \Sord^\lambdaa_\mu \; | \; \mu\in\Lambda^0_{\lambdaa} \}$
for a complete set of simples of
$e_\lambdaa (A^0)^\lambdaa e_\lambdaa$.

\redx{Caveat: Given a modular core there is no reason to suppose that
  we are able to construct such a set directly.}

In this  case the corresponding collection of sets of $\Deltag$-modules
(here specifically denoted $\hat{\Deltag}$ to distinguish from the
other base rings)
as in \eqref{eq:Gdelta}  
give, by 
Lemma~\ref{lem:unitri0},
%
a complete set of simples for $A^0$. 


\medskip

Working now over $K$,
we have the various algebras $e_\lambdaa (A^K)^\lambdaa e_\lambdaa$.
We call a collection of 
sets of $K$-algebra 
modules `\podular' if the appropriate embedding
$G_{e_\lambdaa}\!-\;$ followed by
base change to $K_0$ takes them to
a complete set of simples as above.
\\
Note in particular then that this collection may be indexed in the
same way:
$\Sint = 
\{ \{ \Sint^\lambdaa_\mu \; | \; \mu\in\Lambda^0_{\lambdaa} \} \; | \;
\lambdaa\in \gamma\}$.
\\
(This is the same as to say that the base change of the \podular\ set
\redx{  gives } 
a collection of sets
$ \{ \{ \Sord^\lambdaa_\mu \; | \; \mu\in\Lambda^0_{\lambdaa} \} \; | \;
\lambdaa\in \gamma\}$.)

\newcommand{\Deltagg}{\check{\Deltag}}

\mdef \label{pa:ispivot}
We now have a collection of $A^K$-modules
\redxx{ these we should not call $\Deltag$ without comment!
  they do not match the def. see below. --- so I have changed it. }
$$
\Deltagg = 
\bigcup_{\lambdaa\in \gamma}
\{ \Deltagg_\mu \; | \; \mu\in\Lambda^0_{\lambdaa} \} \; 
$$
obtained by applying the appropriate $G_{e_\lambdaa}\!-\;$ functors to
the collection $\Sint$.
\redxx{Note that calling these $\Deltag$ is an abuse of notation.
it creates a problem later in 2.28.}
\\
As already noted, this collection passes by base change to $K_0$ to the
$\hat{\Deltag}$-modules.
\\
This collection $\Deltagg$ passes by the other base change to $k$ to a set of $A$-modules
$\bar{\Deltag}$.
Once again Lemma~\ref{lem:unitri0} applies.
Note also that these are a complete set of BH modules,
thus they are a pivotal set.
\redxx{Correct!?!?}

\mdef Remark. Note that the implication of semisimplicity flows in one
direction only, so the $\bar{\Deltag}$ modules are not simple in
general.

\medskip\medskip




\newcommand{\ltx}{\triangleleft}  
\newcommand{\ltxeq}{\trianglelefteq} 
\newcommand{\Deltaa}{\varDelta^{\!\!\!\ltx}}  

\mdef\label{HW-def}
We recall the definition of highest weight category  from
\cite{KS1999} (the original reference is
\cite[Definition 3.1]{CPS} but we use an equivalent definition stated
in
\cite{KS1999} as it fits our purpose). 

Let $A$ be a finite dimensional algebra over a field $k$. Fix a
complete set of pairwise non-isomorphic simple $A$-modules
$\{S_\lambda\,\mid \lambda \in \Lambda(A)\}$.
Let  $\ltxeq$ be a partial order on the index set $\Lambda(A)$.
The category $A-\!\!\!\!\mod$ is a \textit{highest weight category}
\textit{with respect to $\ltxeq$}
if for each $\lambda \in \Lambda(A)$
there exists a left $A$-module $\Deltaa_\lambda$,
called \textit{standard module}, such that  
\begin{enumerate}
\item
  There exists a surjective morphism
  $\theta_\lambda: \Deltaa_\lambda\to S_\lambda$, such that if $S_\mu$
  is a composition factor of the ${ker}(\theta_\lambda)$ then
  $\mu \ltx \lambda$.
\item
  Let $P_\lambda$  be the  projective cover of $S_\lambda$.
  There exists a surjective morphism
  $\vartheta_\lambda:P_\lambda\to \Deltaa_\lambda$ such that the
  kernel of
  $\vartheta_\lambda$ is filtered by modules $\Deltaa_\mu$ with
  $\lambda \ltx \mu$.   
\end{enumerate}

\newcommand{\llambda}{\underline{\lambda}}

\ignore{{
\mdef \mulem \label{lem:modHWC}
If a set of  
\most\ 
modules of a modular system for $A$

[obey axiom~1]
\red{(perhaps we want to say this in a core way!...)},

are the $\Deltag$-modules of a core $(\gamma,\leq)$

and filter projectives, then 
$A$-mod is a HWC with respect to ?????.
\red{CORRECT? --- if so then this is a more general proof of the
  Theorem~\ref{th:qh5}.}
\redxx{ The idea of the proof is to use BH
  reciprocity. ...exactly as below!}
\proof
By the upper-unitriangular property of the $\Deltaa$-modules,
i.e. axiom 1 in (\ref{HW-def}),
and reciprocity
\red{equation~(\ref{eq:BH0}) / equation~(\ref{eq:BH1})},
which follows from the modular property as noted in \ref{modBH},
we have $\lmu\geq \llambda$
whenever $(P_{{\llambda}}:\Delta_{\lmu} )\neq 0$
\redx{with $\Delta_{\llambda}$ only at the `head'}
(since the head is determined by $\lambda$ and the multiplicity is 1).
\qed
}}

\mdef \mulem\label{lemm:modHWC}
{\it  Let $A$ be a $k$-algebra with a core $(\gamma, e_-)$.
Let $\Deltag = \{ \Deltag_\mu \; : \; \mu \in \Lambda(A) \}$
be the set of $\beth$-modules of $A$ for the core $(\gamma, e_-)$
(in the sense of Equation~\ref{eq:Gdelta}).
Suppose  the set $\Deltag$ filters projectives and is a 
pivotal set for the strong BH property.
Now let $\prec$ be any
order on $\Lambda(A)$ in which $ \ccc{\mu}<\ccc{x}$
in $\gamma$  implies  
$x \prec \mu$.
Then the category $A$-mod is a highest weight category 
(HWC) with respect to   
$\prec$, with 
 $ \Deltag$
as the set of  standard modules.
}

\proof
\ignore{First note that the definition of HWC requires that  for each
$\lambda\in \Lambda(A)$ there must be a $\Deltaa$-module. By  \ref{eq:Gdeltam} the set $\{ \Deltag_\mu \; : \; \mu \in \Lambda(A) \}$ satisfies this requirement.}
\redxx{ That  this requirement is satisfied ??  --- I do not understand
  --- do you mean just that there is a set of the right size??
  and do you mean to use the $\Deltag$s for this?
}
\ignore{{
We will prove the result by 
considering the modules
$ \Deltag_\mu $  
as candidates for the modules $\Deltaa_\mu$.
}}
For Axiom~1 note from
Lemma~\ref{lem:unitri0}(I) that there is a map $\theta_\lambda$;
and from  Lemma~\ref{Th:upperun} that
$[\ker(\theta_\lambda ) : L_\mu]=0$ unless
$ \ccc{\lambda}<\ccc{\mu}$,
hence
$\mu \prec \lambda$.
For Axiom~2,
let $P_\lambda$ be an indecomposable projective module of $A$,
then by the assumptions there is a filtration 
\[
0\subseteq  M_1   \subseteq M_2\subseteq\cdots\subseteq M_k=P_\lambda
\]
of $P_\lambda$ by $\Deltag_\mu$-modules. We have
$P_\lambda/M_{k-1}\simeq \Deltag_\lambda$, since the $\Deltag$-modules
are pairwise non-isomorphic and their heads form a complete set of
simple modules of $A$. 
\ignore{For each $\lambda, \mu\in \Lambda(A)$ the multiplicity
  $(P_\lambda, \Deltag_\mu)$ is welldefined, since    $\{ \Deltag_\mu
  \; : \; \mu \in \Lambda(A) \}$ is pivotal for the BH property.}
\redx{ By the strong BH {\em basis} property the multiplicities in the
  filtrations are uniquely defined,
  and coincide with the character definition of  $(P_\lambda: \Deltag_\mu)$. }
Thus $(P_\lambda: \Deltag_\mu)=[\Deltag_\mu :S_\lambda]$
by the strong BH property. \ignore{ $(P_\lambda, \Deltag_\mu)$ is defined and}
\redxx{ say how this is defined?? } Hence, if
$(M_{k-1}: \Deltag_\mu)\neq0$
\redxx{ do you mean $M_{k-1}$ ?? } 
then Lemma~\ref{Th:upperun} implies that $\lambda\prec\mu.$
\redxx{ which way around??! }
\qed

\ignore{
\blu{An alternative version of the of the previous Lemma}

\mdef \mulem\label{lemm:modHWC}
{\it Let $A$ be a $k$-algebra with a core
$(\gamma, e_-)$.   Let $\{ \Deltag_\mu \; : \; \mu \in \Lambda(A) \}$
be the set of $\Deltag$-modules of $A$ for the core $(\gamma, e_-)$
(in the sense of Equation~\ref{eq:Gdelta}).
If the set of $\Deltag$-modules filters projectives and
 is a 
pivotal set for the strong BH property, then the
category $A$-mod is a
highest weight category 
(HWC) with respect the reverse of any total order 
$\prec$ on $\Lambda(A)$ such that $ \ccc{\mu}<\ccc{x}$
in $\gamma$  implies $\mu \prec x$. In particular, $\{ \Deltag_\mu \; : \; \mu \in \Lambda(A) \}$  forms a complete set of  standard modules of $A.$
}

\proof
\ignore{First note that the definition of HWC requires that  for each
$\lambda\in \Lambda(A)$ there must be a $\Deltaa$-module. By  \ref{eq:Gdeltam} the set $\{ \Deltag_\mu \; : \; \mu \in \Lambda(A) \}$ satisfies this requirement.}
\redxx{ That  this requirement is satisfied ??  --- I do not understand
  --- do you mean just that there is a set of the right size??
  and do you mean to use the $\Deltag$s for this?
}
\ignore{{
We will prove the result by 
considering the modules
$ \Deltag_\mu $  
as candidates for the modules $\Deltaa_\mu$.
}}
We will
show that the modules $ \Deltag_\mu $
satisfy the axioms of ~\ref{HW-def} given an order as specified.
 By Lemma~\ref{Th:upperun} the modules   $ \Deltag_\mu $ satisfy Axiom~1 of ~\ref{HW-def}.
For Axiom~2,
let $P_\lambda$ be an indecomposable projective module of $A$,
then by the assumptions there is a filtration 
\[
0\subseteq  M_1   \subseteq M_2\subseteq\cdots\subseteq M_k=P_\lambda
\]
of $P_\lambda$ by $\Deltag_\mu$-modules. We have
$P_\lambda/M_{k-1}\simeq \Deltag_\lambda$, since the $\Deltag$-modules
are pairwise non-isomorphic and their heads form a complete set of
simple modules of $A$. 
\ignore{For each $\lambda, \mu\in \Lambda(A)$ the multiplicity
  $(P_\lambda, \Deltag_\mu)$ is welldefined, since    $\{ \Deltag_\mu
  \; : \; \mu \in \Lambda(A) \}$ is pivotal for the BH property.}
 By the strong BH basis property the multiplicities in the
  filtrations are uniquely defined,
  and coincide with the character definition of  $(P_\lambda: \Deltag_\mu)$. Note that $M_{k-1}$ is the kernel of the natural surjection $P_\lambda\to P_\lambda/M_{k-1}\simeq \Deltag_\lambda$, and 
  \[
0\subseteq  M_1   \subseteq M_2\subseteq\cdots\subseteq M_{k-1}
\]
is a filtration of $M_{k-1}$ by $\Deltag_\mu$-modules.
 Further, by the strong BH basis property $(M_{k-1}: \Deltag_\mu)=(P_\lambda: \Deltag_\mu).$
Hence if $(M_{k-1}: \Deltag_\mu)\neq0$
then $[\Deltag_\mu: L_\lambda]=(P_\lambda: \Deltag_\mu)\neq0$. Finally,  Lemma~\ref{Th:upperun} implies that $  \mu\prec\lambda.$
\qed
\blu{end of the alternative version}}

\ignore{ \red{Finally we want a strong modular core version????
 It looks like the above plus (\ref{pa:ispivot}) is exactly what we
  need. Yes???
  Do you want to state the Theorem?}
}

\mth \label{th:strongmodular}
Let  $A$ be  a $k$-algebra. Suppose:
\\ (I) $A$ is a 
strong  modular core algebra with  
modular system $(K,K_0,k)$ and
core $(\gamma,e_-)$.
Let $\bar{\Deltag}$
be the set of $\Deltag$-modules of the $k$-algebra $A$ for the core
$(\gamma, e_-)$, obtained by base change as in \ref{pa:ispivot}.
Let $\prec$ be any  
order on $\Lambda(A)$ in which
$ \ccc{\mu}<\ccc{x}$   in $\gamma$  implies $x \prec \mu.$
\\ (II)  
the set $\bar{\Deltag}$ filters  projective $A$-modules.
\\
Then $A$-mod is a HWC with respect to $\prec$,
with  $\bar{\Deltag}$ as the set of  standard modules.

\proof
By  \ref{pa:ispivot} the set $\bar{\Deltag}$ is pivotal set for the BH property. Now apply Lemma~\ref{lemm:modHWC}. \qed



\ignore{{
(Probably best to ignore this draft `intro' for now!)
 These are {\bf rough} {\em notes}. 

The layout   
is as follows. 
First we define the even partition category and algebra, and introduce
some notation. Then we prove theorems determining the fundamental
invariants of representation theory in the `generic' case. 
Next we
 set up theorems needed to address the `exceptional' cases.

A primary motivation for this study is its intrinsic interest,
which will be clear below. Other motivations   
include  potential applications
analogous to those of the partition algebra itself (for example in
statistical mechanics and invariant theory \cite{,}).
But we also have in mind that several diagram algebras have a
geometrical aspect to their representation theory, analogous to
Kazhdan--Lusztig theory \cite{,,}, 
originally observed post-hoc \cite{,}, 
which, while given a beautiful interpretation in certain
individual cases for example in \cite{Stroppel}, remains intriguing as
a collective phenomenon. 
The aim is
to provide further examples, of comparable complexity to the Brauer and
partition cases, and hence investigate if this geometry can be given a
natural uniform explanation. 

In order to prove our main theorems  
we first note some basic representation theoretic
properties. Then we investigate the structure of a key quotient
algebra. Here we use some gram matrix calculations, assuming
familiarity with symmetric group representation theory as in James
\cite{}. 
Next we determine the Bratteli diagram \cite{,,} of a natural tower of
inclusions of even partition algebras.
... combinatorics ... translation functors ... 
Then we construct functors relating part of the representation theory
to that of known objects.
Finally we investigate gram matrices in the general case. ...

}}


\section{Simple index theorem for $\Pl_n$} \label{ss:ideal2}
\subsection{The quotient algebras $P^{\mm}_n$ of $\Pl_n$}


\mdef Define quotient algebra
\[
P_n^{\mm} \; = \; \Pl_n /  \Ina  
\]
Schematically:
\[
\includegraphics[width=1in]{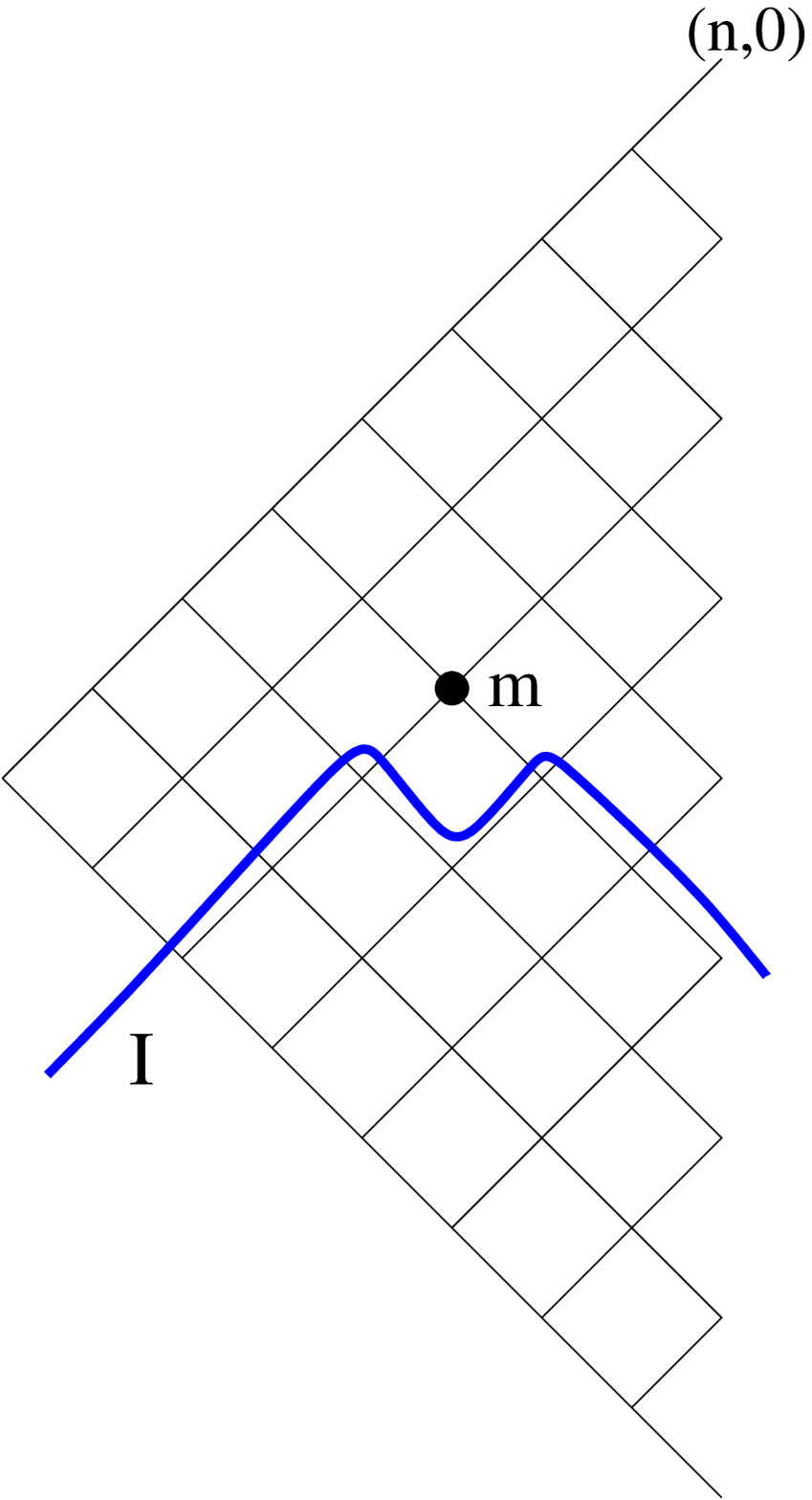}
\]


\mdef By $a^\mm \in P_n^\mm$ we understand the element of which
$a^\mm$ is a representative
(and similarly for $b^\mm$). 

{{\mlem\label{lem:basis:mP^mm}
The ideal $P_n^{\mm}  a^{\mm} P_n^{\mm} $ in $P_n^{\mm}$ has basis 
 $\mP^{\mm}_n$.  
}}
\proof This follows from Lem.\ref{lem:idbas1} and the construction. \qed


\newcommand{\std}{Specht}
\newcommand{\mmu}{\underline{\mu}}
\newcommand{\ttSp}{{\mathtt {Sp}}}  

\newcommand{\Smmm}{\times_{i=1}^{l} S_{m_i}} 
\newcommand{\Smmmrev}{\times_{i=0}^{l-1} S_{m_i}}
\newcommand{\Smm}{S_{\mm}}
\newcommand{\sigmaz}{\rho}


\mdef
For $\mm\in\gammal$ define
\beq \label{eq:Smm1}
\Smm \; = \; \Smmm
\eq
For $\sigmaz = (\sigmaz_1 , \sigmaz_2, ..., \sigmaz_l ) \in\Smm$,
define
$w^a_\sigmaz = 
w_\sigmaz \in P^l_n$ 
as the image of $\sigmaz$  
realised on the propagating lines in $a^{\mm}$
as follows. 
Here is an example for $l=2$:
\beq \label{eq:wsig}
w_\sigmaz \; =\;  \;
\raisebox{-.21in}{
\includegraphics[width=5cm]{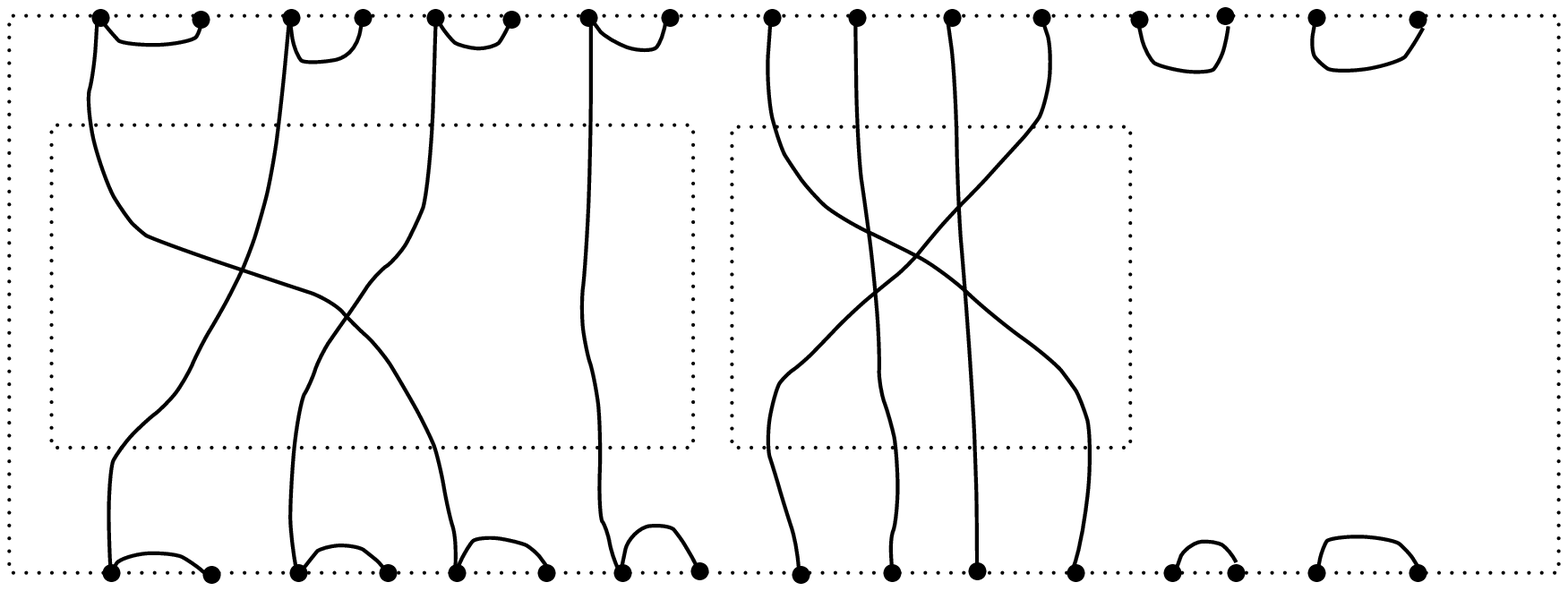}}
\eq
--- we put $\sigmaz_1$ on the co-1 propagating lines (the second
interior box of lines in the example) in the natural way;
$\sigmaz_2$ on the co-2 lines
(the first interior box in the example); and so on.

\mdef \label{pa:aba0}
For $\mm\neq 0$ we define $w^b_\sigmaz \in P^l_n$ analogously to
$w_\sigmaz$, but realised on the propagating lines of $b^\mm$.
We may similarly define $w^{ab}_\sigmaz$ on the propagating lines of
$a^\mm b^\mm$; and $w^{ba}_\sigmaz$ analogously.
Note that
\beq \label{eq:wsigx}
w^b_\sigmaz w^b_{\sigmaz'} = w^b_{\sigmaz \sigmaz'}  ,
\hspace{.6in}
w^{a}_\sigmaz w^{ba}_{\sigmaz'} = w^{a}_{\sigmaz \sigmaz'}  
\eq
and so on.



\begin{figure}
\[
\includegraphics[width=4.4cm]{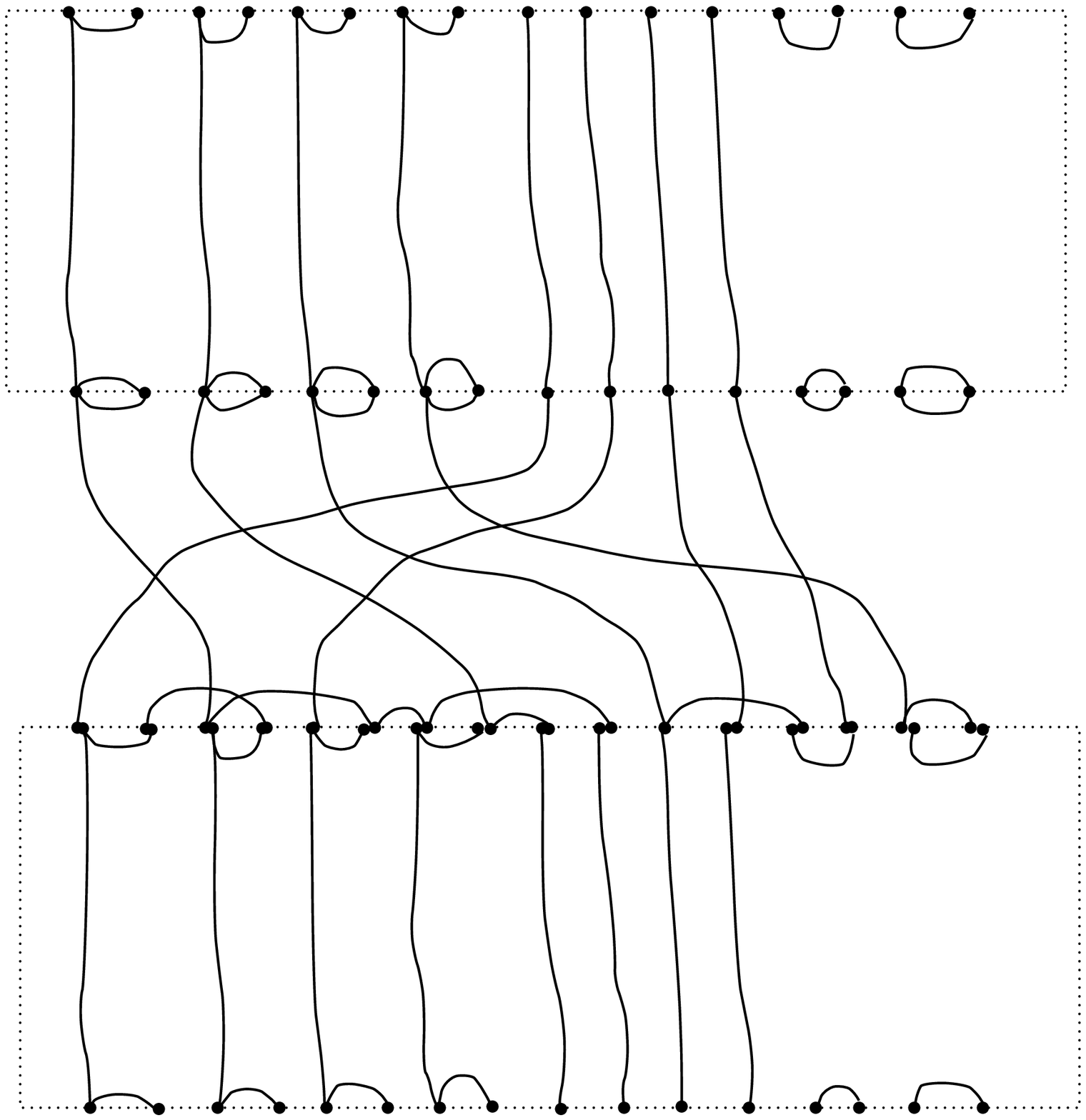}
\qquad \qquad
\includegraphics[width=4.4cm]{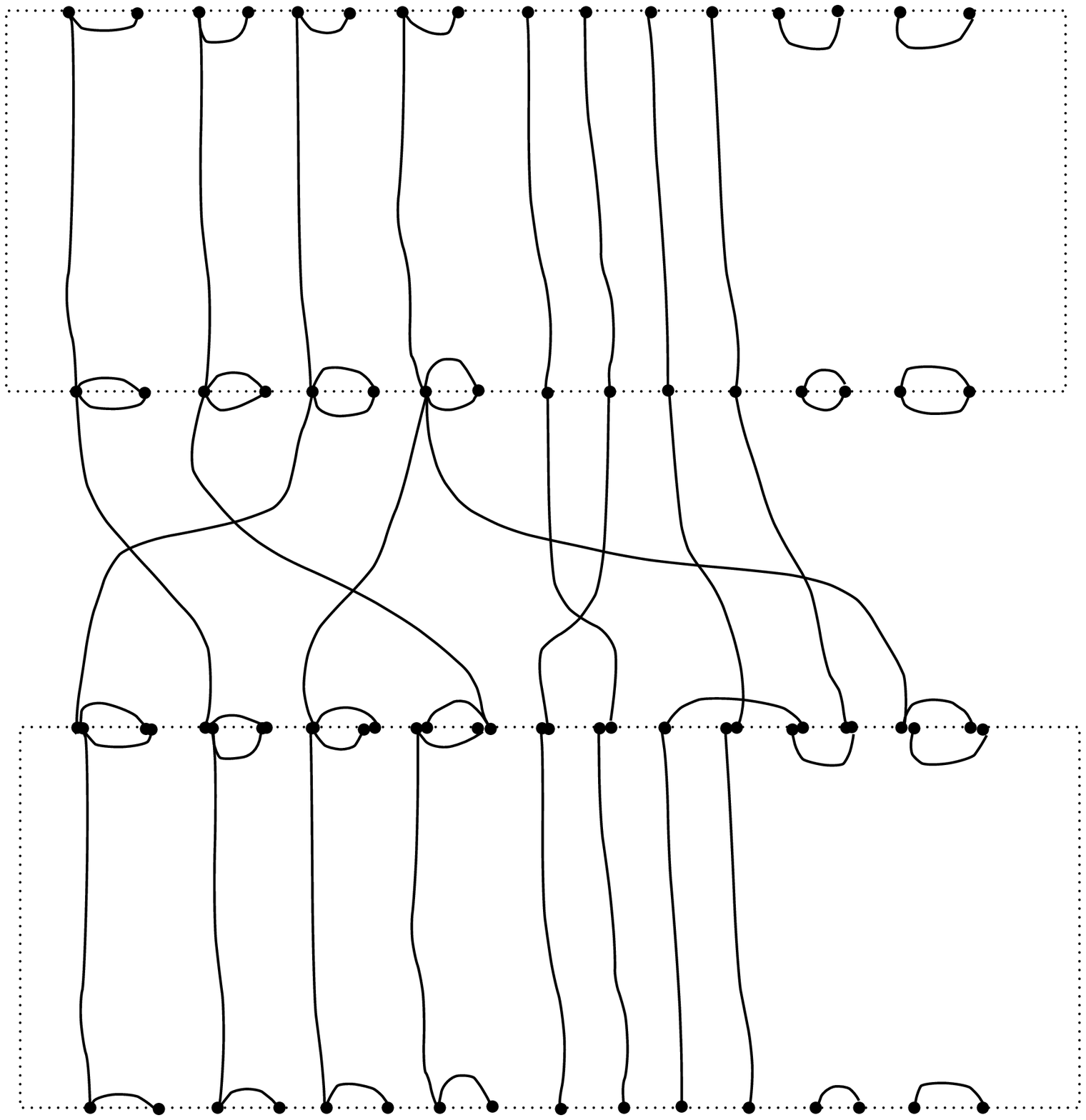}
\]
\caption{ Partitions illustrating $a^{\mm} P^{\mm}_n a^{\mm}$
in case $l=2$.
\label{fig:groupc1}}
\end{figure}

\mdef \label{pa:arithmetic1}
Consider 
partitions of form $q=a^\mm p a^\mm$, as illustrated in
figure~\ref{fig:groupc1}.
(Partitions of form $q = b^\mm p b^\mm$ are directly similar.)
Note
that the number of co-$i$ parts of $q$ cannot be
greater than that of $a^\mm$ for any $i$,
unless this is the result of two or more propagating parts coming
together,
such that $i' + i'' \equiv i$ mod. $l$ (or $\sum_j i_j \equiv i$).
Now cf. (\ref{de:gammal}) and the definition of $P^\mm_n$. 
Thus for example 
the first element in the figure, $q_1 = a^\mm p_1 a^\mm$ say, is
in $I^{<\mm}$ and hence zero in $ P^\mm_n$  by the quotient.
To see this explicitly consider $\hashi(q_1)$, and in particular $\hashi(q_1)_1$.
We may follow down the leftmost co-1 part from the top $a^\mm$
factor.
If this meets a $co-i>1$ part from the lower $a^\mm$ (as here)
then, by the pigeonhole principle, the number of
propagating parts decreases.
 
The second example is essentially a `permutation' of form $w_\sigma$
(in this case up to factors of $\delta$).
Note then that only permutations 
within each of the $l$ group factors are possible, 
if we work in the $P^{\mm}_n$ quotient.


{\mlem{  \label{lem:groupc} \label{lem:aDa1}
Let $\kk$ be a commutative ring, and $\delta\in\kk$.
Fix $l$ and $n$, and $\mm \in \gammal \setminus \{ 0 \}$.
Then we have the following.
\\ (bI) The map from $S_\mm$ to $P^l_n$ given by $\sigma \mapsto w^b_\sigma$ 
has image in $ b^\mm P^l_n b^\mm$.
\\ (bII) 
Idempotent subalgebra
$ b^{\mm} P^{\mm}_n b^{\mm} $ in $P^\mm_n$ has basis
$\{ [w^b_\sigma ] = w^b_\sigma +I^{<\mm} \; | \;  \sigma \in S_\mm \}$.
\\ (bII') The map
$\sigma \mapsto  [w^b_\sigma ]$ 
defines a map from $ S_\mm$ to
the basis in $ b^{\mm} P^{\mm}_n b^{\mm} $
that is a group isomorphism.
\\ (bIII)
This map gives an algebra isomorphism
\[
b^{\mm} P^{\mm}_n b^{\mm} \; \cong \; \kk S_\mm   .
\]
}}
\proof
(bI)
Note that $w^b_\sigma \in P_n^l$ and $b^\mm w^b_\sigma b^\mm = w^b_\sigma$. 
(bII) Noting (\ref{eq:abaa}),
this follows essentially from Lemma~\ref{lem:1} and the pigeonhole
principle,
as in (\ref{pa:arithmetic1}) (replacing $a^\mm$ with $b^\mm$).
(bII') This follows from (bII) and (\ref{eq:wsigx}).
(bIII) Follows from (bII) and (bII'). 
\qed

\mdef \label{pa:aba}
Note that by the same argument the idempotent subalgebras
$b^\mm P^\mm_n a^\mm$, 
$a^\mm P^\mm_n b^\mm$, and
$a^\mm P^\mm_n a^\mm$
are all isomorphic to $\kk S_\mm$ (except that the `$aa$' case requires
$\delta$ to be a unit). 



{\mth{ \label{th:index1}
Suppose $\kk$ is a field with $\delta \neq 0$.
Then we have simple index set 
$ \Lambda(\Pl_n ) = \sqcup_{\mm \in \gammal} \Lambda(\kk S_\mm ) . 
    $
}}

\proof Noting Theorem~\ref{th:index1core},
now apply Th.(\ref{th:genLam}) and 
substitute using (\ref{lem:groupc})
\redx{and (\ref{pa:aba}).}
\qed

\mdef Remark. There is a straightforward strengthening of
Th.\ref{th:index1} to the case $\delta=0$ in most cases.
We omit this for brevity.



\subsection{Aside on symmetric groups and Specht modules} \label{ss:ss}

\label{pa:reviewSS}
In light of Lem.\ref{lem:groupc} {\em et seq.}, 
it is useful to recall some properties of the group algebra
$\kk S_\mm$ (as defined in (\ref{eq:Smm1});
for example 
$\kk S_\mm = \kk (S_{n-2i} \times S_{i})$
in case $\mm = (n-2i,i)$). 
We focus for brevity in this exposition on the case $l=2$.
The generalisation is straightforward.

\mdef
Define $\Lambda_i = \{ \lambda \vdash i \}$,
the set of integer partitions; and 
$$
\Lambda_{\mm} = \; \times_{i=1}^l \Lambda_{m_i}  
  \; = \; \Lambda_{m_1} \times \Lambda_{m_2} \times ... \times \Lambda_{m_l}
$$


\mdef \label{de:pre} 
Here an element $e$ 
in a $k$-algebra $S$
($k$ some commutative ring)
 is `preidempotent' 
if $ee=ce$ for some $c \in k$.
If $c$ a unit then $e$ may be renormalised as an idempotent.
`Primitivity' of $e \in S$ 
means that $ewe=c_w e$ for some $c_w \in k$ for all 
$w\in S$.

\mdef
Recall (cf. \cite{James} and \cite[\S43]{CurtisReiner62}, say) that for each 
$\underline{\lambda} = 
      (\lambda^1, \lambda^2 , ..., \lambda^l ) \in \Lambda_{\mm}$
there exists a 
primitive preidempotent 
\beq \label{eq:preidSS}
e_{\underline{\lambda}}  =  e_{\lambda} = e'_{\lambda^1} e''_{\lambda^2} ...
\eq 
in 
$\Z S_{\mm} $
(primes indicate belonging to different factors)
such that the left ideal
\[
{\ttSp}_{\lambda}
   = \; \kk S_{\mm}   \; e_{\lambda} 
\]
is a Specht module for $l=1$ and hence a generalised Specht module otherwise.

In our case we may choose $e=e_\lambda$ so that $e=e^{op}$, as in (\ref{de:flip})
(this follows from one of the 
well-known constructions for preidempotents in $\Z S_n$ \cite{JamesKerber}). 


\mdef \label{de:blambdabasis}
For definiteness we have in mind a tableau-labelled basis 
$$
b_{\lambda} = b_{\lambda^1} \times b_{\lambda^2} \times ...
$$ 
for
$\ttSp_{\lambda}$. 
This is a basis encoded as $l$-tuples of standard
sequences such as $b_{(3,1)} = \{ 1112,1121,1211 \}$. 
The details of the corresponding explicit basis contruction 
can be found for example in \cite{James},
but the full details will not be needed here. 

\medskip

\redxx{ CAN we say that (left-ideal) Specht modules have a free basis that may be
  extended to a free basis of $\kk S_\mm$??
  I guess yes. That might simplify things later, e.g. in 6.7.}


\medskip

\mdef \label{de:gcvf}
Returning to (\ref{de:pre}),
suppose algebra $S$ has
a $k$-linear map $op:S \rightarrow S$, written  $s \mapsto s^{op}$,
that is an involutive antiautomorphism.
Suppose primitive preidempotent  $e=e^{op}$.
This then allows us to define a $k$-bilinear form $(-,-)_e$
on $Se$ as follows: 
for $se, s'e \in Se$ we have 
\beq \label{eq:gcvf1}
(se)^{op} s'e = es^{op} s'e = \; c_{s^{op} s'} \; e .
\eq
Now set $(se,s'e)_e = c_{s^{op} s'}$.
Note that for $a \in S$ we have 
$(se,as'e)_e = c_{s^{op} as'} = c_{(a^{op} s)^{op} s'} = (a^{op} se , s'e)_e$.
That is,  the form is {\em contravariant} with respect to $op$ \cite{Green80}.

\mdef
This form is useful in  studying the corresponding $S$-module morphism from $Se$ to its
contravariant dual \cite{Green80}.
Note that the choice of $e$ is not unique in our ideal construction $M=Se$, 
and although $Se$ does not
depend on the choice (up to isomorphism) the form does depend on it by
an overall factor. Thus the form is not canonical on $M$.
However for symmetric groups there is a
good choice of form, due to James, that encodes representation theory
within a very useful organisational scheme \cite{James}.
Over fields of char.0 the constant $c$ is always a unit and these subtleties can be
ignored,
as we will see in \S\ref{ss:cvf}.



\section{Polar decomposition of partitions in $\mPl_n$} \label{ss:pol}

We describe a version for $P^l_n$ of the $P_n$ polar decomposition 
\cite[p79-80]{Martin94}.


\mdef \label{pa:freerightS}
Let $\mm \in \gammal$. 
Let $B^{\mm}$ denote  the natural `diagram' basis of $P_n^{\mm} a^{\mm}$
(i.e. the basis of certain partitions $p \in P^l_n a^\mm$ 
where $p$ is understood to 
mean $p+I^{<\mm}$).
Next we describe $B^\mm$. 
Note that the left ideal 
$P_n^{\mm} a^{\mm}$ has a natural right action of $\Smm$
upon it (see e.g. Lem.\ref{lem:aDa1} and (\ref{pa:aba})). 
Indeed it is a free right $\kk S_\mm$-module. 
It follows that we may partition 
$B^{\mm}$ 
into orbits of the right action of $S_\mm$.

\mdef
The basis $B^{\mm}$ consists of elements $p$ representable in the following
form.
\\
(1) The restriction of $p$ to the 
`top' (unprimed) subset of vertices 
consists of:
\\
for each $i=0,l-1,l-2,...,1$,
$m_i$ co-$i$ parts that belong to propagating parts; 
and some further number of co-0 non-propagating parts. \\
The restriction of $p$ to the `bottom' (primed) subset of vertices consists of:
\\
for each $i$,
$m_i$ order-$i$ parts that belong to propagating parts; 
and $(n-r_\mm )/l$ further  order-$l$ non-propagating parts. \\
(2) The propagating connection from
top to bottom  
for a propagating part
may be drawn as a line
from the first
(lowest numbered) vertex
of the part on the top
to the first on the bottom. 


\mdef \label{de:relnonc}
We say that $p \in B^{\mm}$ is {\em relatively non-crossing} if, for
each $i$,  
the co-$i$
propagating lines in the representation above are pairwise non-crossing.  


An example of a relatively non-crossing partition
    {(in case $l=2$)}
is given by:
\[
\reflectbox{ \rotatebox[origin=c]{180}{
  \includegraphics[width=5.4cm]{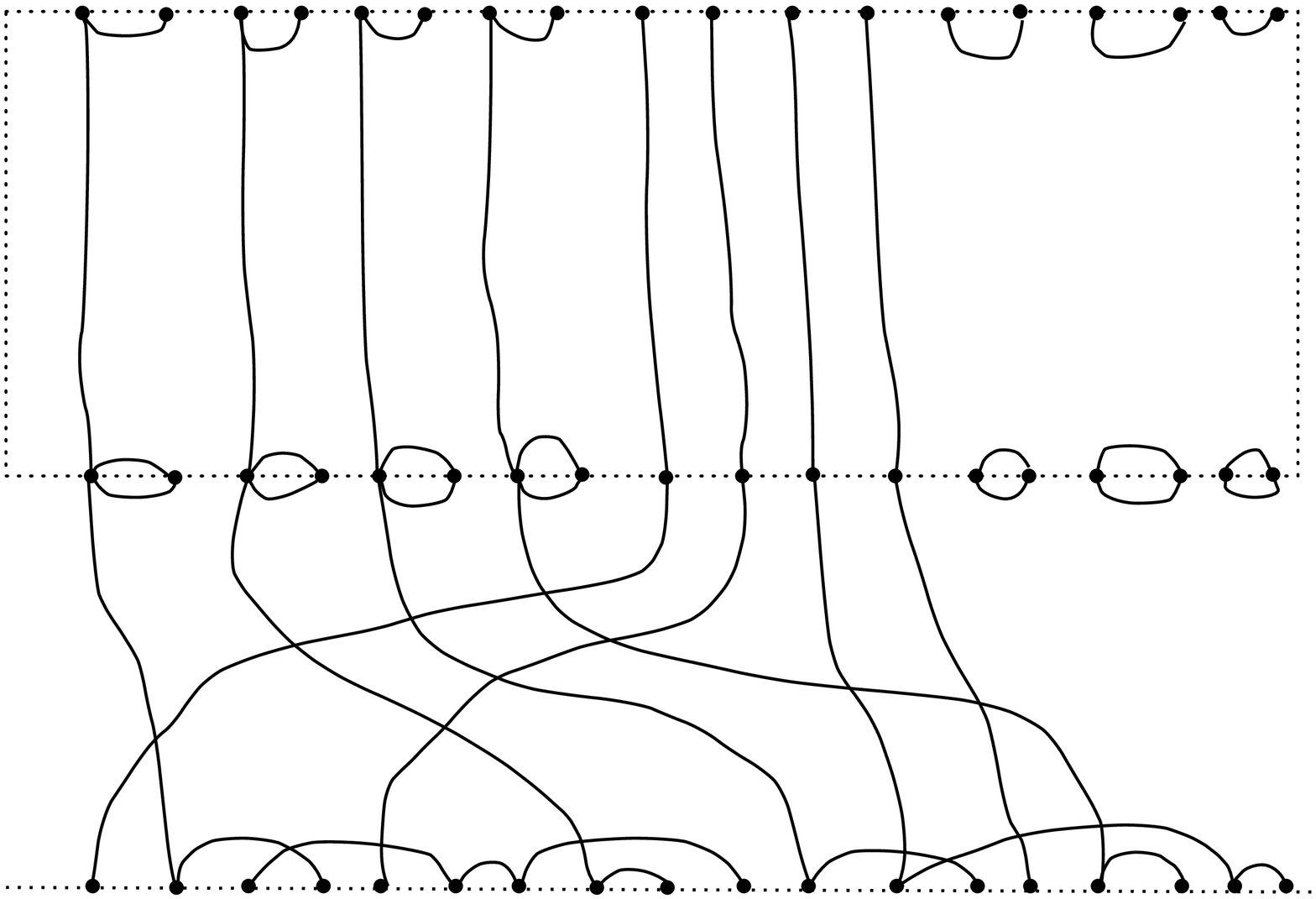}
  }}
\]
The isolated loops in this picture can be ignored
(we assume $\delta \neq 0$ here for simplicity)
or replaced with a suitable `meander'. 
They are drawn to demonstrate
that the element lies in $P_n^{\mm} a^{\mm}$ 
(here in case $\mm = (m_0,m_1)=(4,4)$). 

{\mlem{ \label{lem:tranv1}
Let $\mm \in \gammal$. 
Each orbit of the
right
$S_\mm$ action on basis $B^{\mm}$  of 
$P_n^{\mm} a^{\mm}$
as in (\ref{pa:freerightS})
contains a unique representative element with the 
relative non-crossing property.
Let $T^{\mm}$ denote the relative non-crossing transversal. Then 
\[
B^{\mm} = \{ pw \; | \; p \in T^{\mm}; \; w \in S_\mm \}
\]
where $w$ acts in the natural way.
\blux{In particular $T^\mm$ is a $\kiss$- basis of the free 
right $\kiss$-module $P_n^{\mm}a^{\mm}$.}
\qed
}}


\mdef \label{lem:nu4}
It will be apparent that any partition $p$ in $\mP^{\mm}_n$
(as defined in (\ref{eq:dePmn}))
can be written
in a generalisation of the usual partition algebra {\em polar decomposition}.
That is, we have the following. 

{\mlem{ \label{lem:polarPm}
Each $p \in \mP^{\mm}_n$ can be written 
in a factored form as 
$$
p= a w_\sigma b
$$ 
where $a$ is relatively non-crossing (as in (\ref{de:relnonc})),
i.e. $a \in T^{\mm}$;
$w_\sigma=(w_1, w_2,...)$ 
is as in (\ref{eq:wsig});
and $b$ is a flipped relatively non-crossing partition
($b^\star \in T^\mm$). 
The factorisation in this form is unique.
\qed
}}

\mdef
It is convenient to denote the decomposition by
\beq \label{eq:ketcbra}
p \mapsto | p \rangle \; w(p)  \; \langle p |
\eq


\section{Basic integral representation theory of $P^l_n$} \label{ss:int}
\newcommand{\lnu}{\underline{\nu}}

Here
we
\redx{
aim to 
  construct
a modular system for $P^l_n$ based on 
analogues of Specht modules.
}


\begin{figure}
\[
\ig[width=2.176in]{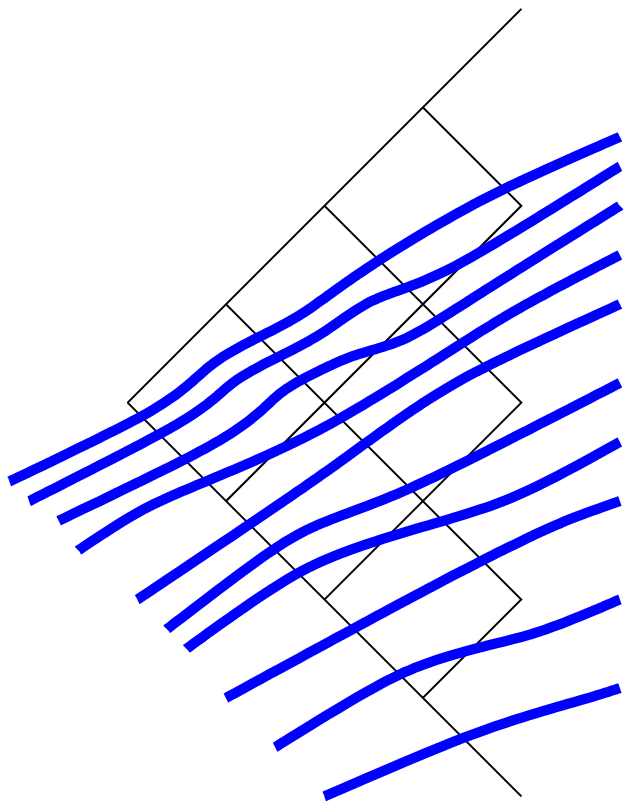}
\hspace{2.1cm}
\ig[width=1.56in]{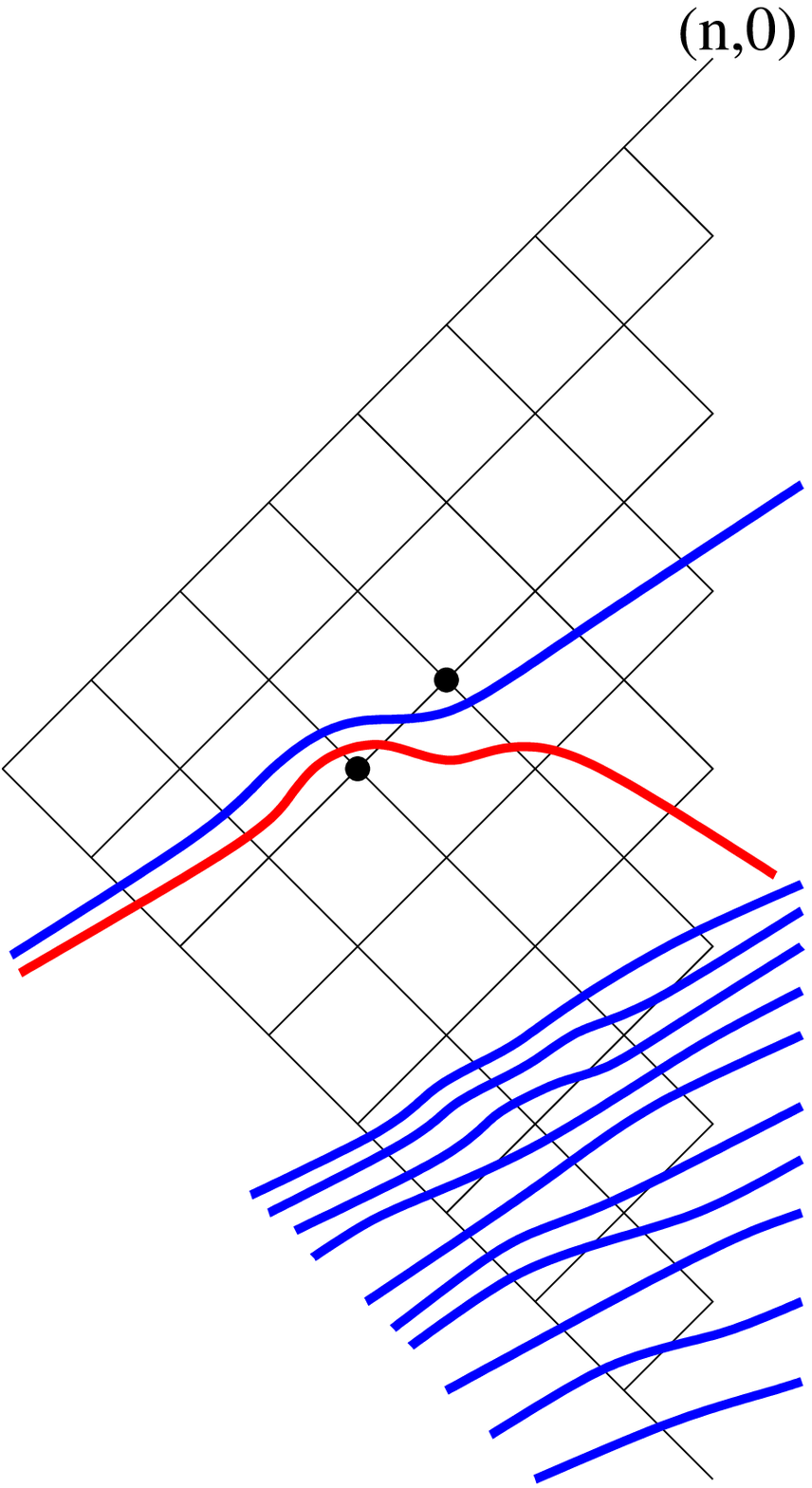}
\hspace{.1in}
\ig[width=1.56in]{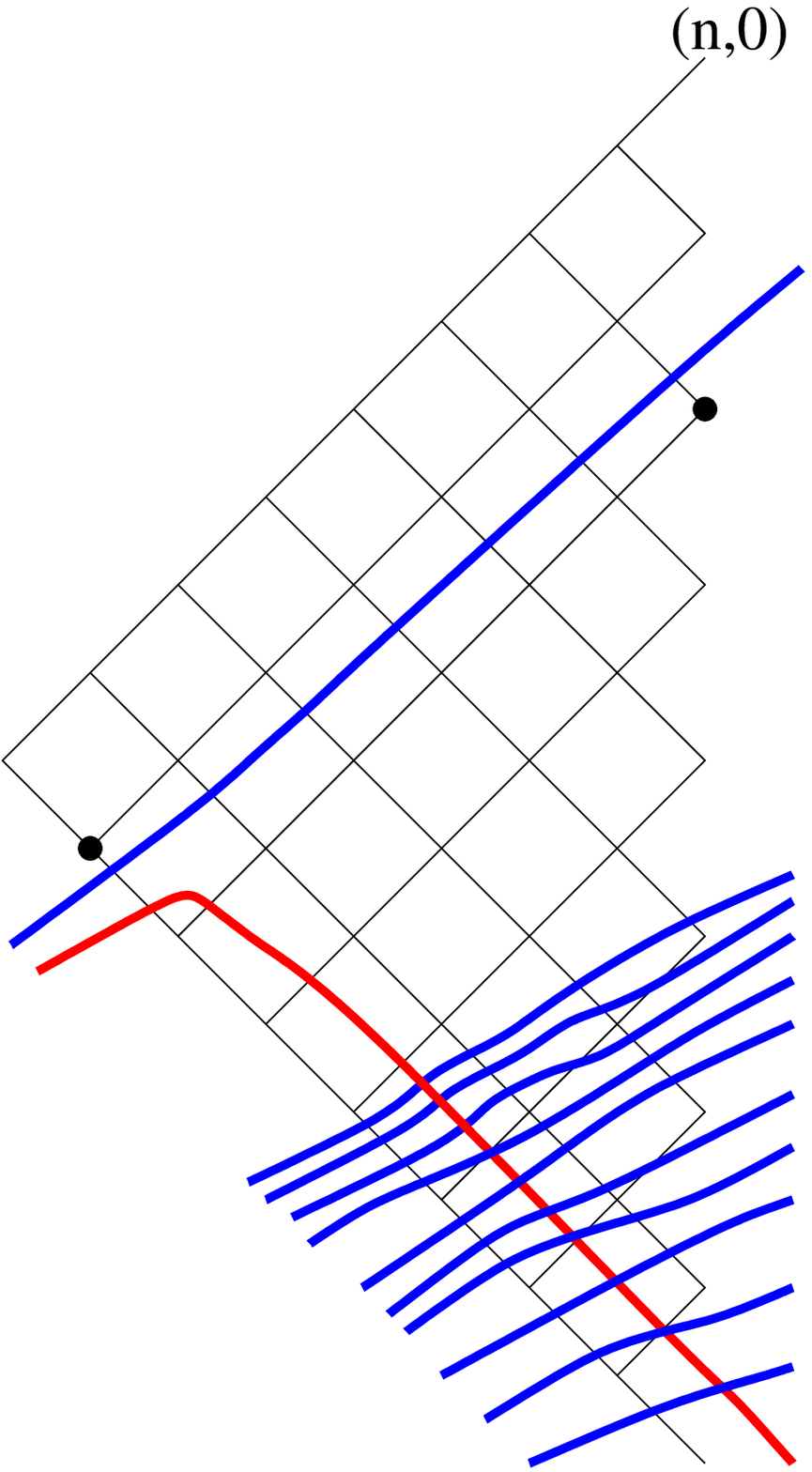}
\]
\caption{\label{fig:ePnpx}
  (a) Schematic for $\rhd$-sectioning the $\gammal$ poset ($l=2$).
  (b) Schematics for comparing $J_{\mm(i)}$ (below blue line)
  with $I^{<\mm(i+1)}$ (below red line).}
\end{figure}


\mdef\label{id-chain}
Recall that any poset can be refined to a total order. 
Let $(\gammal, \rhd)$ be a total order refining $(\gammal , >)$
(e.g. as indicated by Fig.\ref{fig:ePnpx}). 
Let us define
\[
J_{\mm} = \sum_{\mm' \unlhd \mm} P^l_n a^{\mm'} P^l_n
, \qquad \qquad
J_{\lhd\mm} = \sum_{\mm' \lhd \mm} P^l_n a^{\mm'} P^l_n
\]
For convenience define $\mm(1), \mm(2) , ...$ as the elements of
$\gammal$ in the total order
so
$$
J_{\mm(1)} \subset J_{\mm(2)} \subset J_{\mm(3)} \subset ...
\redx{ \subset J_{(n,0,0,...,0)} = P^l_n  }
$$


\mdef \label{lem:step2} \mulem 
{\it
Write $P=P^l_n$ for a moment.
We have an isomorphism of bimodules
\beq \label{eq:bison}
J_{\mm(i+1)} / J_{\mm(i)}
\cong
P a^{\mm(i+1)} P / I^{<\mm(i+1)} .
\eq
}

\noindent
Proof:
We have
\[
J_{\mm(i+1)} / J_{\mm(i)} = 
  \frac{P a^{\mm(i+1)} P + P a^{\mm(i)} P +... }{P a^{\mm(i)} P +P a^{\mm(i-1)} P +...}
 \hspace{1in} \]
\[ \hspace{1in}
\cong
 \frac{P a^{\mm(i+1)} P  }{ P a^{\mm(i+1)} P \cap (P a^{\mm(i)} P +P a^{\mm(i-1)} P +...)} 
=
 \frac{P a^{\mm(i+1)} P  }{ P a^{\mm(i+1)} P \cap J_{\mm(i)}} 
\]
by the second isomorphism theorem. 
Thus we may consider the `numerators' in (\ref{eq:bison}) to be the same, up to isomorphism;
and compare the
`denominators' (the submodules that are quotiented by).
The argument proceeds in two steps.
\\
(I)
For any order $\rhd$ refining $>$ we have that $\mm(i+1) > \mm$
implies $\mm(i+1) \rhd \mm$. But for any total order this implies
$\mm(i) \unrhd \mm$. 
Thus $J_{\mm(i)} \supset I^{<\mm(i+1)}$.
(Cf. Fig.\ref{fig:ePnpx}(b).) 

\noindent
(II) Consider the denominator in the third expression.
In particular consider $\;$ 
$
P a^{\mm(i+1)} P \cap P a^{\mm(j)} P
$
for $j=i, i-1, ...$.
By  Lemma~\ref{lem:idbas1} every partition $p$ lies in a
unique highest ideal of the form $Pa^\mm P$,
and there is a basis of partitions.
 Since $\mm(j) \not\geq \mm(i+1)$   we have (cf. (\ref{le:keymm}))
that
\[
P a^{\mm(i+1)} P \cap P a^{\mm(j)} P
\subseteq I^{<\mm(i+1)}
\qquad (j\leq i)
\]
\newcommand{\mmp}{\mm(i+1)} 
Combining with the inequality in the other direction from (I),
noting that $ P a^{\mm(i+1)} P \supset  I^{<\mm(i+1)} $, we see that
\[
J_{\mm(i+1)}/ J_{\mm(i)} \cong \frac{P a^{\mmp} P }{I^{<\mmp}} = P^{\mmp} a^{\mmp} P^{\mmp}
\]
\qed



\subsection{Globalisation functors and  quotient algebras $A^l_n$}

\newcommand{\barr}[1]{\overline{#1}}
\newcommand{\iow}{\iota_{W}}


Note that a partition of a set $S$ determines a partition of a subset
$S'$ by restriction.
In particular an element $p$ of $\mP_n$ determines an
element $p|_{[1,n-l]}$ of $\mP_{n-l}$ by restricting to the first
$n-l$ pairs of elements `top and bottom'.
Similarly $p|_{[l+1,n]}$ restricts to the last $n-l$ pairs.
(Again this determines an element of  $\mP_{n-l}$ in the obvious way.)
Note in particular for $n>l$ that the restriction
$p|_{[l+1,n]}$
of a partition $p$ in
$  W_b^l \Pl_n W_b^l$
removes $l$ top elements from the same part and $l$ bottom elements
from the same part, and hence $p|_{[1,n-l]}$ 
lies in $ \Pl_{n-l}$.
A similar property holds on restricting $  W_{}^l \Pl_n W_{}^l$,
or indeed $  \barr{W_b^l} \Pl_n \barr{W_b^l}$
(using (\ref{de:latflip})), and so on.
Given one of these cases, 
let us write $\iota_{W}$ for this restriction map.
Schematically ($W^3$ case):
\[
\rb{-2.75cm}{
\includegraphics[width=1.1704in]{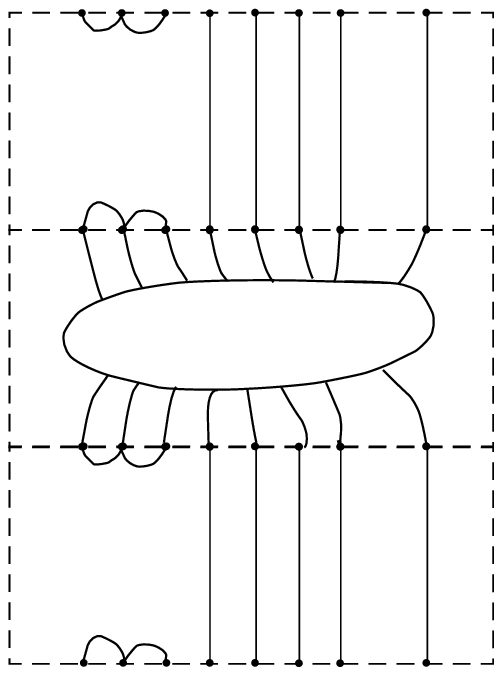}}
\;\; \longrightarrow \;\;
\rb{-.81cm}{
\includegraphics[width=1.1704in]{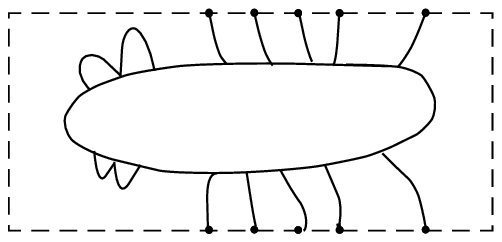}
}
\;\; = \;\;
\rb{-.81cm}{
\includegraphics[width=.91404in]{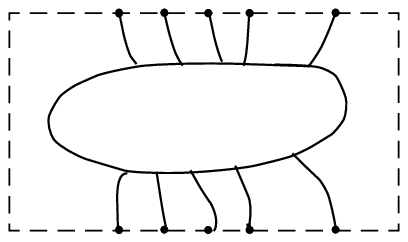}
}
\]


{\mth{ \label{th:UEU}
Fix $l$. 
For   $n>l$    
the maps $\iow$ give 
isomorphisms of algebras:
$$
\;  W_b^l \Pl_n W_b^l \cong       \Pl_{n-l}
\hspace{.1in} \mbox{ and }\hspace{.1in} 
 \;  \barr{W_b^l} \Pl_n  \barr{W_b^l} \cong       \Pl_{n-l} .
$$
Similarly for $n \geq l$ 
we have the (not necessarily unital) algebra isomorphisms
$
 \;  W^l \Pl_n W^l \cong    \delta   \Pl_{n-l}  .
$
and
 $
 \;  \barr{W^l} \Pl_n \barr{W^l} \cong    \delta   \Pl_{n-l}  .
$
}}
\proof Exactly analogous to the $P_n$ case as in \cite{Martin94}.
Schematically ($W^3$ case), consider the figure above.
\qed


\mdef \label{de:Gfunctor}
As in (\ref{de:Gfunctorx}) 
we define  `short' functors
$$
G_{\barr{W}} :
\barr{\Wlb} \Pl_n \barr{\Wlb} - \!\mmod \;\; \rightarrow \;\; \Pl_n - \!\mmod
$$
and $F_{\barr{W}} $.
And similarly (for $\delta$ invertible)
$ \;
G_W : \Wlb \Pl_n \Wlb - \!\mmod \;\; \rightarrow \;\; \Pl_n - \!\mmod
$.
By Theorem~\ref{th:UEU} we will consider these as functors between
$\Pl_{n-l} - \!\mmod$  and $\Pl_{n} - \!\mmod$.




\mdef Fix $l$. Define the quotient algebra 
$$
A^l_n \;  
  \; = \; \Pl_n / \Pl_n a^{(n-l,0,0,..,0)} \Pl_n
  \; = \; \Pl_n / \Pl_n \Wl \Pl_n 
  \; = \; \Pl_n / \Pl_n \Wl_b \Pl_n .
$$

\mdef Example: In case $l=1$, $A^1_n = kS_n$.

{\mlem{ \label{lem:Abasis}
There is 
a basis for $A^l_n$ consisting of 
partitions in which every part is propagating 
with $\cora(p_i) = \rang(p_i) $
and  no part
has $\cora(p_i) > l$ or  $\rang(p_i) > l$.
}}
\proof By Lem.\ref{lem:12}. \qed

\newcommand{\oddsandendsone}{
\footnote{Is this true and do we care?!:

{\mth{ \label{th:EaE}
Fix $l=3$. 
For $n \geq l$, 
$$ \;  
\Pl_n / \Pl_n a^{(n-l,0,0,..,0,1)} \Pl_n \; \cong  \;  A^{l-1}_n  .
$$
}}
\proof 
By Lem.\ref{lem:1} the ideal $\Pl_n a^{(n-l,0,0,..,0,1)} \Pl_n$ 
contains every  partition in
$\Pl_n$ except
for those in which every part is propagating but no part
has $\cora(p_i) > l-1$ or  $\rang(p_i) > l-1$.
Now compare the definition of $ A^{l}_n  $ with Lem.\ref{lem:Abasis}.
\qed

}
}


\mdef Let 
$
H^l_n
$ 
denote the subset of idempotents $a^\mm$ with $\mm\in \hl_n$.

\newcommand{\mPlp}{\mP^{l-}}

Let  $\mPlp_n \subset \mPl_n$ denote the subset of $l$-\tone\ 
partitions of form $ShT$, where $S,T$ are permutations and $h \in H^l_n$.
That is, $\mPlp_n$ is the set of 
$l$-\tone\ partitions having parts with at most $l$ elements per row,
and all parts propagating. 

Note from Lem.\ref{lem:1} the following.

{\mlem{ \label{lem:11}
Fix $l$. Algebra 
$A^l_n$ has basis the subset $\mPlp_n$. 
\qed
}}


\mdef \label{de:compatible}
An element of $\mPlp_n$ 
is partially characterised by the restricted partition of
the upper (resp. lower) row into parts of size $l,l-1,...,1$.
We say that an ordered pair $(a,b) \in \mPl_n$ are {\em compatible} if,
considering the rows that meet in composition $ab$, every 
part 
from $a$
(resp. $b$) is 
a union of parts 
from $b$ (resp. $a$). 

Note that $ab \equiv 0$ in $A^l_n$ unless compatible. 


\medskip

In case $l>n$ there can be no non-propagating part in a partition in
$\mPl_n$. In this case $\Pl_n$ coincides with  
Kosuda's party algebra \cite{Kosuda06}. 
We will be interested, though, in general $n$ for each fixed $l$.


\subsection{Long functors from $\kk S_{\mm}-\!\mmod $} 

\medskip

\newcommand{\GF}{G}

\mdef \label{de:Ga}
Note from Lem.\ref{lem:groupc}
and (\ref{pa:aba})
that $P^\mm_n a^\mm $ is a right $\kk S_\mm$ module.
\redxx{-We should explain the action. In the algebra isomorphism case it
  is clear. But otherwise?
There are various ways to do this. In the end the detail is probably
not so important to us here. But we are trying to learn something, by
playing. }
 (The case $\mm=0$, $\delta=0$
 can be included by identifying
 $\kk S_0$ with the ground ring.)

Now for given $\mm\in\gammal$ define the {\em long functor} 
$$
\GF_a : \kk S_{\mm}-\!\mmod \; \rightarrow \; P^{\mm}_n-\!\mmod
$$
by 
$
 \GF_a M =  P^{\mm}_n a^{\mm} \otimes_{\kk S_\mm} M .
$



\medskip

Here we want to recall the isomorphism between $G_a M$ and
$P^\mm_n a^\mm M$ that holds when $M$ is a left ideal, so that standard modules
inherit the properties from both constructions.
Noting the invariant basis number (IBN) property we write $rank_\kk(M)$ for
the basis number of a free $\kk$-module with finite basis.

\ignore{{
 

{\mlem\label{lem:FM}
Let $A$ be a finite dimensional $k$-algebra,
$M$ be a left ideal ${A}$-module  with $k$-basis $B$,
and $F$ be a right free $A$-module of finite rank with $A$-basis $T$.
Then
(I) $TB =  \{tb \;|\; t\in T ; b \in B \}$  is basis of $FM.$
(II)  $T\otimes B =  \{t\otimes b \;|\; t\in T ; b \in B \}$  is basis
    of $F\otimes_A M.$
(III) The map $\mu: T\otimes B \rightarrow TB$ given by $t\otimes
b \mapsto tb$ lifts to a well-defined $k$-space isomorphism
$\mu : F\otimes_A M \rightarrow FM$.
If $F$ is an $A'$-$A$-bimodule then $\mu$ is an isomorphism of left $A'$-modules.
\ignore{{
\\
(I) $FM=\bigoplus_{p\in T} p  M;$
\\
(II) for each $p\in T$ we have $dim( pM) =dim(M);$
\\
(III)   the set $\{pw\;|\; p\in T ;w\in B_M\}$  is basis of $FM.$
\\
(IV) $M' = F \otimes_{A} M$ has a basis with a factored form 
\[
B_M'  \; = \{ p \otimes w \; | \; p \in T ; \; w \in B_M  \}.
\]
}}
}

}}

 
{\mlem\label{lem:FM}
Let $A$ be a finite rank $\kk$-algebra,
$M$ be an ${A}$-module  with $\kk$-basis $B$,
and $F$ be a right free $A$-module of finite rank with $A$-basis $T$.
Then
\\
(I)  $T\otimes B =  \{t\otimes b \;|\; t\in T ; b \in B \}$  is a basis
    of $F\otimes_A M.$
\\
If in addition $M$ is a left ideal then
\\
(II) $TB =  \{tb \;|\; t\in T ; b \in B \}$  is a basis of $FM.$
\\
(III) The map $\mu: T\otimes B \rightarrow TB$ given by $t\otimes
b \mapsto tb$ lifts to a well-defined  
$\kk$-module isomorphism
\redxx{right?
(two finite rank free $\kk$-modules are isomorphic - as $\kk$-modules/as free
$\kk$-modules - iff same rank??)}
$\mu : F\otimes_A M \rightarrow FM$.
If $F$ is an $A'$-$A$-bimodule then $\mu$ is an isomorphism of left $A'$-modules.
\ignore{{
\\
(I) $FM=\bigoplus_{p\in T} p  M;$
\\
(II) for each $p\in T$ we have $dim( pM) =dim(M);$
\\
(III)   the set $\{pw\;|\; p\in T ;w\in B_M\}$  is basis of $FM.$
\\
(IV) $M' = F \otimes_{A} M$ has a basis with a factored form 
\[
B_M'  \; = \{ p \otimes w \; | \; p \in T ; \; w \in B_M  \}.
\]
}}
}


\proof
(I) The $A$-basis property 
says $F = \oplus_{t\in T} tA$
where each $tA$ is a copy of the right regular $A$-module. 
Thus $F \otimes_A M = ( \oplus_{t\in T} tA ) \otimes_A M
 = \oplus_{t \in T} ( tA \otimes_A M ) $
by tensor-distributivity (see e.g. \cite[\S 12]{CurtisReiner62})
and hence 
$ F \otimes_A M = \oplus_{t \in T} t \otimes M$ since $A \otimes_A {}_A M \cong M$.
Thus
$rank_{\kk}(F \otimes_A M) = | T |.|B |$.
Furthermore
if $G$ is a $\kk$-basis of $A$ then 
$F$ has $\kk$-basis $TG$, so 
$\{ tg \otimes b \; : \;
t \in T$, $g \in G , \; b \in B   \}$
spans $F \otimes_A M$.
But $ta \otimes b = t \otimes ab$ so, fixing $t$, the set
$\{ t \otimes b \; : \;  b \in B \}$
spans the same  
$\kk$-module  as
$\{ ta \otimes b \; : \; a\in A , b \in B \}$.
Thus $T\otimes B$  spans $F \otimes_A M$.
Since $| T\otimes B|  = | T |.|B |$ we are done.
\\
(II)
\redxx{OK?... }
We may compare for example with 
\cite[4.12]{lam2012lectures}.
The set  
$FM$ is a $\kk$-module, so $FM= \oplus_{t} tAM$ is a direct sum of
$\kk$-modules - specifically $AM=M$ is a $\kk$-module and $t$ is a
formal symbol.
Thus $FM= \oplus_t tM $ is a formal direct sum of copies of $M$.
Since $M$ is itself $\kk$-free, we have a formal direct sum of formal
direct sums.
\ignore{{
By the Exchange Theorem there exists a $k$-basis $G$ of $A$ with
$G \supseteq B$, so $TB \subseteq TG$. $TG$ is a $k$-basis of $F$ by
freeness,
so $TB$ is linearly independent. But also $FM = (TA)M = T(AM)$
and $AM = M$ and $B$ spans $M$, so $TB$ spans $FM$.
}}
\\
(III) Well-definedness follows from the balanced map property of
tensor products.
The map is surjective by construction, and hence an
isomorphism by
IBN.
Commutativity of the map with the
algebra action is also by construction.
\qed



\ignore{{

\proof  
(I) First note that by the freeness of $F$ we have that $T A=FA$. Hence
 $$F M=TA M\subseteq +_{p\in T } p  M\subseteq F M,$$
 which implies that $FM =+_{p\in T} p M.$ To show the sum $+_{p\in T} p M$ is direct, 
let  $p, q\in T$ with $p\neq q$. Then $p  M\cap q  M\subseteq pA\cap qA$ since $M$ is a left ideal of $A$.  We have $pA\cap qA=  \{0\}$ since $p$ and $q$ are $A$-linearly independent.    Hence $p  M\cap q  M=\{0\}$  

(II) 
\blu{another proof: To prove II it is enough to show that  $pB_M=\{pw\; |\; w\in B_M\}$ is a basis of $pM$, since $|B_M|=|pB_M|$. To this end, extend the basis  $B_M$ of $M$ to a basis of $A$, call it $B'_M$. Then $pB'_M=\{pw\; |\; w\in B'_M\}$ is a basis of $pA$, since it spans $pA$ and $dim(A)=dim(pA)$. Now $pB_M\subseteq pB_M'$ hence $pB_M$ is linearly independent. Hence  $pB_M$ is a basis of $pM$.}
  The set map
\begin{equation*}
\begin{split}
\mu_p:  B_M &\to p M\\
w&\mapsto pw
\end{split}
\end{equation*}
  can be extended  uniquely in to a linear map from $M$ to $p M$.  Further, the map $pw\mapsto w$ \blu{(do we
 need to say more about this map or it is okay like that?)}\blu{ the argument goes like this: First we show that that $pB_M$ is a basis of $pM$. To this end, extend the basis  $B_M$ of $M$ to a basis of $A$, call it $B'_M$. Then $pB'_M=\{pw\; |\; w\in B'_M\}$ is a basis of $pA$, since it spans $pA$ and $dim(A)=dim(pA)$. Now $pB_M\subseteq pB_M'$ hence $pB_M$ is linearly independent. Hence  $pB_M$ is a basis of $pM$. Then \begin{equation*}
\begin{split}
\mu'_p:  pB_M &\to  M\\
pw&\mapsto w
\end{split}
\end{equation*} is a set map  which 
  }
  is the inverse of (the linear extension of ) $\mu_p$.  Hence $dim( p M) =dim(M)$.

III) By part I the given set spans $ F  M.$ Again by using parts I and II of this Lemma   we have $dim(FM)=|T |\cdot|B_M |.$ Hence the set $\{pw\;|\; p\in T ;w\in B_M\}$ is linearly independent. 

IV) By the freeness of $F$ we have a spanning set of form $pw \otimes b$ with
$p \in T $, $w \in A$ and $b \in B_M$.
But $pw \otimes b = p \otimes wb$ so, fixing $p$, the set
$\{ p \otimes b \; : \;  b \in B_M \}$
spans the same space as
$\{ pw \otimes b \; : \; w\in A, b \in B_M \}$.
Thus the claimed set is spanning.
To see that it is independent note that $F$ is a free
right
$A$-module
with free basis $T$
and so
(since $A \otimes_A {}_A M \cong M$) we have
$dim(M' ) = | T |\cdot|B_M |$.

\qed

{\mlem{ \label{de:basMM} \label{lem:basisc} \label{lem:ny1}
Let  $A,M$ and $F$ be as given in Lemma~\ref{lem:FM}.   
Let $B$ be finite dimensional algebra.  If $F$ is  an $A$-$B$-bimodule and $M$ is a left ideal $A$-module then
  the surjective multiplication map
$\mu: p \otimes w \mapsto pw \in FM$
   is injective, so $$ F \otimes_{A} M\cong  F  M$$  as left $B$-modules. 
}}
\proof
If $M$ is a left ideal $A$-module then   the set map 
\begin{equation*}
\begin{split}
F \times M&\to FM\\
(p,w)&\mapsto pw
\end{split}
\end{equation*}
is $A$-balanced. Hence by the universal property of tensor products there exists a group homomorphism
\begin{equation*}
\begin{split}
\mu:F \otimes M&\to F M\\
p\otimes w&\mapsto pw.
\end{split}
\end{equation*}
%
%
  Further, for any $q\in B$ we have 
$$\mu(q(p\otimes w))=\mu(qp\otimes w))=(qp)w=q(pw)=q\mu(p\otimes w)$$
  hence  $\mu$ is a left $B$-module morphism.
By Lemma~\ref{lem:FM}
$FM$ and $ F \otimes_{A} M  $ have the same dimension, hence surjectivity of $\mu$ implies that $\mu$ is injective.
}}

%

\ignore{{
(I). By Lemma~\ref{lem:tranv1}  we have a spanning set 
$\{ pw \otimes b \; : \;
p \in T^\mm$, $w \in S_\mm , \; b \in B_M   \}$.
But $pw \otimes b = p \otimes wb$ so, fixing $p$, the set
$\{ p \otimes b \; : \;  b \in B_M \}$
spans the same space as
$\{ pw \otimes b \; : \; w\in S_\mm, b \in B_M \}$.
Thus the claimed set is spanning.
To see that it is independent note that $P^\mm_n a^\mm$ is a free
right
$k S_\mm$-module
\red{(see (\ref{pa:freerightS}))}
with free basis $T^\mm$
\red{can we not perhaps use this property in some more general way here?}
and so
(since $A \otimes_A {}_A M \cong M$) we have
$dim(M' ) = | T^\mm |.|B_M |$.
}}

\ignore{{
\qed

{\mlem{ \label{de:basM1} \label{lem:basisc} \label{lem:ny1}\blu{Corollary?}
(I) Any $P^\mm_n$-module of form 
\[
M' = P^{\mm}_n a^{\mm} \otimes_{\kiss} M
\]
where $M$ is a  
${\kiss}$-module,
has a basis with a factored form 
\[
B_M^{\mm} \; = \{ p \otimes w \; | \; p \in T^{\mm}; \; w \in B_M \}
\]
where $B_M$ is a basis for $M$.

(II) If $M$ is a left ideal ${\kiss}$-module then
$M' \cong  P^{\mm}_n a^{\mm}  M$ as left $P^{\mm}_n$-modules.
}}

\proof Part I is a special case of part IV of Lemma~\ref{lem:FM}.  II) By  Lemma~\ref{lem:tranv1} the right $\kiss$-module $P^\mm a^\mm$ is free, now use Part II of Lemma~\ref{de:basMM}. 
\qed


}}



\subsection{Standard/Specht modules} 
\newcommand{\sspecht}{\ttSp} 
\newcommand{\pond}{\operatornamewithlimits{\otimes}}
\newcommand{\oxn}{\underset{n-l}{\otimes}} 
\newcommand{\ox}{\underset{\raisebox{1pt}{\tiny$n\!\!-\!\!l$}}{\otimes}} 
\newcommand{\oxm}{\underset{{\raisebox{5pt}{\tiny$\mm$}}}{\otimes}}
\newcommand{\WW}{W^l_{}}


\mdef \label{de:std1}
Fix $n$ and $l$. For 
$\mm \in \gammal$ 
and 
$
\mmu = (\mu^{1},\mu^{2},...,\mu^{l})  \in \Lambda_\mm
$
we may define the  
$\mmu$-`\std' module 
of $\Pl_n$ as the module obtained by
applying the functor $\GF_a$ in (\ref{de:Ga})
to the 
$\mmu$-Specht module $\sspecht_{\mmu}$  of 
 ${\kiss}$:
\beq \label{eq:spechtinflate}
\Delta_{\mmu} \; = \; \GF_a \sspecht_{\mmu}
\eq
(see \S\ref{ss:ss} for details of $\sspecht_{\mmu}$).
Note that this is a $\Pl_n$-module since it is a $P_n^{\mm}$-module. 
Define  $\Lambda_0(P^l_n) = \cup_{\mm\in\gammal} \Lambda_\mm $  and 
let $\Delta^{l,n} = \{\Delta_{\mmu} \; | \; \mmu\in
 \Lambda_0(P^l_n) \}$.



\mdef \mulem \label{de:Llmu}
\redxx{(We are possibly going to remove and replace this a bit later.)}
{\it
For $k=\C$ and $\delta \neq 0$ the set
$\{ L_{\mmu} = \head \Delta_{\mmu} \}_{\mmu\in\Lambda_0(P^l_n)}$
is a complete set of simple modules. 
}
\proof
  Compare Th.\ref{th:index1}(I) with
  Lem.\ref{lem:unitri0}(I) and  (\ref{de:std1}).
  \qed


\mdef \label{de:specht}
Let $\mm\in\gammal$ and $\underline{\mu} \in \Lambda_\mm$,
and 
$e_{\lmu} \; = \; \prod_i e_{\mu^i}^i$ 
as in (\ref{eq:preidSS}).
Define 
\[
\specht^n_{\underline{\mu}} 
  \; = \; P^{\mm}_n a^{\mm}_n
             w^{ba}_{e_{\underline{\mu}}} 
\]
where  $e_{\underline{\mu}}$  
thus acts 
as in (\ref{eq:wsig} - \ref{eq:wsigx}).

\redxx{Now we should point to EXAMPLES! There are some in next section.}

\newcommand{\www}{w^{ba}_}

{{\mlem \label{lem:basisBSp}
Let  $b_{\lmu}$ be a basis of $\mmu$-Specht module $\sspecht_{\mmu}$  of 
$\kk S_\mm$ (cf. \ref{de:blambdabasis}, \cite{James});      
and $T^{\mm}$ be
a non-crossing transversal in $P^\mm_n$
as defined in Lemma~\ref{lem:tranv1}.
Then 
$$
B^{\lmu}_{\sspecht{}} \;
   := \{ \;  p\, \www\omega \; | \; p\in T^\mm, \; \omega \in b_{\lmu} \}  
$$
is a basis of $\specht^n_{\underline{\mu}} .$
}
\proof By ~\ref{lem:tranv1} the set $T^{\mm}$ is a  
$\kk S_\mm$-basis
of $P^{\mm}a^{\mm}$.
Now in  Lemma~\ref{lem:FM} part  II let
$A=\kk S_\mm$, $M=\sspecht_{\mmu}$ and $F=P_n^{\mm}a^{\mm}$ to obtain the result.
\redxx{$k$ to $\kk$????}
\redxx{The anomaly at $\delta=0$: the issue is $F$ being a free
  module. Formally this is false, but with our added comments above it seems
  fine. OK??}
\qed
}

\begin{figure}
\newcommand{\tmpsize}{.68in}
\[
\includegraphics[width=\tmpsize]{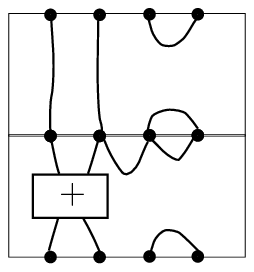}
\hspace{.3cm}\includegraphics[width=\tmpsize]{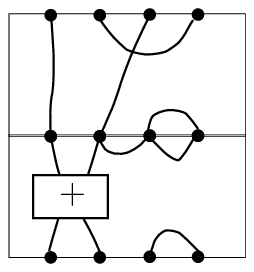}
\hspace{.3cm}\includegraphics[width=\tmpsize]{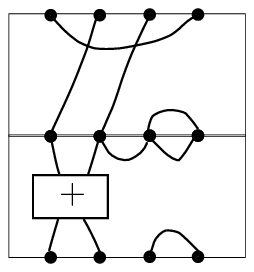}
\hspace{.3cm}\includegraphics[width=\tmpsize]{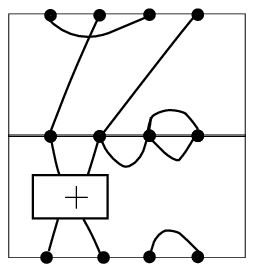}
\hspace{.3cm}\includegraphics[width=\tmpsize]{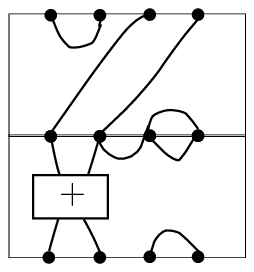}
\hspace{.3cm}\includegraphics[width=\tmpsize]{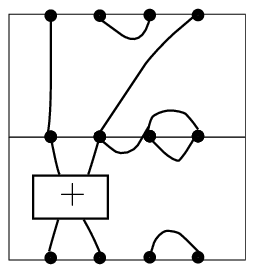}
\] \[
\hspace{.0cm}\includegraphics[width=\tmpsize]{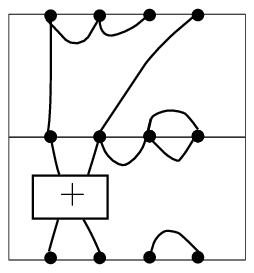}
\hspace{.3cm}\includegraphics[width=\tmpsize]{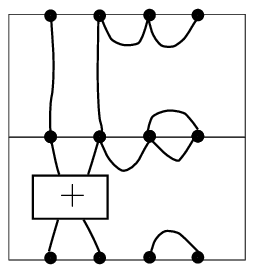}
\hspace{.3cm}\includegraphics[width=\tmpsize]{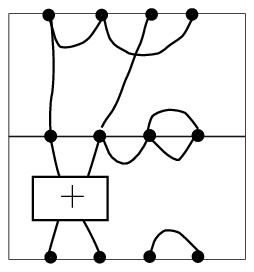}
\hspace{.3cm}\includegraphics[width=\tmpsize]{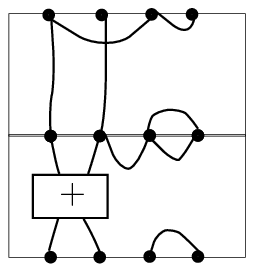}
\]
\caption{Diagrams for the basis for $\specht^4_{((2),\emptyset)} $.
  See main text. \label{fig:BofS}}
\end{figure}

\mdef
\redxx{Either we do the ba case (in which the lower box has an extra arc)
  or we need to explain more! }
Example.
Fig.\ref{fig:BofS} shows
a diagrammatic realisation of the basis for
$ \specht^4_{{((2),\emptyset)}}$.
The box labeled $+$ denotes the $\Z$-linear combination
corresponding to the preidempotent
$e_{(2)} = 1-\sigma_1 \in \Z S_2$.
Note cf. \cite{Martin94} that this
diagram calculus is well defined.
So far, then, the pictures give combinations of partitions
--- but then
finally these partitions are understood to represent the classes in
the module of which they are representative. 

\redxx{COMMENT: Formally the parenthesis won't work for $\delta=0$
  (gives 0).
  This is where
  (in my thinking) we note that $S_0$ is trivial, so we can basically
  omit $e_\mu$.
But the real question is `what we use all this for?'. Probably the
interesting point is the tower of recollement? But we do not get that
far in this paper!}


\begin{figure}
\newcommand{\tmpsize}{.8in}
\hspace{.653cm}\vspace{-.71cm}\includegraphics[width=\tmpsize]{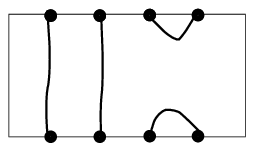}
\hspace{.2cm}\includegraphics[width=\tmpsize]{latexdraw/barrW.eps}
\hspace{.12cm}\includegraphics[width=\tmpsize]{latexdraw/barrW.eps}
\hspace{.2cm}\includegraphics[width=\tmpsize]{latexdraw/barrW.eps}
\hspace{.12cm}\includegraphics[width=\tmpsize]{latexdraw/barrW.eps}
\hspace{.2cm}\includegraphics[width=\tmpsize]{latexdraw/barrW.eps}
\[
\includegraphics[width=\tmpsize]{latexdraw/BofS01.eps}
\hspace{.3cm}\includegraphics[width=\tmpsize]{latexdraw/BofS02.eps}
\hspace{.3cm}\includegraphics[width=\tmpsize]{latexdraw/BofS03.eps}
\hspace{.3cm}\includegraphics[width=\tmpsize]{latexdraw/BofS04.eps}
\hspace{.3cm}\includegraphics[width=\tmpsize]{latexdraw/BofS05.eps}
\hspace{.3cm}\includegraphics[width=\tmpsize]{latexdraw/BofS06.eps}
\]
\vspace{1in}
\hspace{2.85cm}\vspace{-.71cm}\includegraphics[width=\tmpsize]{latexdraw/barrW.eps}
\hspace{.2cm}\includegraphics[width=\tmpsize]{latexdraw/barrW.eps}
\hspace{.12cm}\includegraphics[width=\tmpsize]{latexdraw/barrW.eps}
\hspace{.2cm}\includegraphics[width=\tmpsize]{latexdraw/barrW.eps}
\vspace{-.95in}
\[
\hspace{.0cm}\includegraphics[width=\tmpsize]{latexdraw/BofS07.eps}
\hspace{.3cm}\includegraphics[width=\tmpsize]{latexdraw/BofS08.eps}
\hspace{.3cm}\includegraphics[width=\tmpsize]{latexdraw/BofS09.eps}
\hspace{.3cm}\includegraphics[width=\tmpsize]{latexdraw/BofS10.eps}
\]
\caption{Action of $\barr{W}$ on the basis
  for $\specht^4_{((2),\emptyset)} $. \label{fig:WonS}}
\end{figure}


{\mlem{
\label{lem:FWS}
      {\rm (I)}   For
$\mm \in h^l_{n}$ and $\lmu\in\Lambda_\mm$,
      \redxx{$\lmu \in \sqcup_{\mm\in h^l_{n}} \Lambda_\mm$, }
    $F_{\barr{W}} \specht_{\lmu}^{n}    =0$.
{\rm (II)}   For
$\mm \in \gamma^l_{n-l}$ and $\lmu\in\Lambda_\mm$,
    \redxx{$\lmu \in \sqcup_{\mm\in \gamma^l_{n-l}} \Lambda_\mm$, }
    $F_{\barr{W}} \specht_{\lmu}^{n}
    = {\barr{W}} \specht_{\lmu}^{n}
    \cong \specht_{\lmu}^{n-l}$.
    (Here we assume for simplicity that $\delta\neq 0$.)
}}
\proof
(I)
Here $\barr{W} P^\mm_n = 0$.
(II)
Note that for $x \in I^{<\mm}_{n-l}$ then $x\otimes ww^\star \in I^{<\mm}_{n}$.
Considering the bases 
one finds that
$b \mapsto b\otimes ww^\star$
(with $b$ of form $d + I^{<\mm}_{n-l}$)
gives rise to an injection
$\specht_{\lmu}^{n-l} \rightarrow  {\barr{W}} \specht_{\lmu}^{n} $.
The image is fixed by $\barr{W}$ so $  F_{\barr{W}} \specht_{\lmu}^{n}$
contains $ \specht_{\lmu}^{n-l}$ as a submodule.
We \redx{
observe complementarily that the map
$\barr{W} P^l_n \barr{W} \rightarrow P^l_{n-l}$
(here NB the standing assumption)
takes $p \in \mP^{\mm}_n$ to $p' \in \mP^{\mm}_{n-l}$.
Indeed $\barr{W} P^l_n \barr{W} \cap \mP^{\mm}_n
= \mP^{\mm}_{n-l} \otimes ww^\star$. 
Furthermore by construction
each element of $\barr{W} \specht^n_{\lmu}$
is a subset of 
$\barr{W} P^l_n \barr{W}$
(to see this consider e.g. Fig.\ref{fig:WonS}),
so we can apply  $ x \leftrightarrow x \otimes ww^\star$.
On the other hand we know by Lemma~\ref{lem:1}(II)
that whenever the propagating index
is changed it is reduced.
Thus
the image of the injection is also spanning.
\\
(To see this less formally, 
from the example --- Fig.\ref{fig:WonS} ---
one can see the following types of cases for $\barr{W} d$
for $d \in B^{\lmu}_{\sspecht}$:
(1) last $l$ vertices already connected in $d$;
(2) propagating lines are at most permuted;
(3) the propagating index is reduced.)
} \qed


{\mlem{ \label{le:Diso}
    We have an isomorphism of $P^l_n$-modules:
$\Delta_{\lmu} \cong \specht^n_{\lmu}$.
}}
\proof
\redxx{PM: Although our story is probably logically correct, it seems
  maybe it was easier to follow when we made explicit for this case
  the basis of the modules (not 
  for the proof here, but for the combinatorics) --- even though they
  are implied directly by 6.7.
In the previous draft of course the bases were explicit. $\;$ 
Also} \redxx{
There are some
  things I don't yet fully find clear here. In the paper we define $G_e$
  when $e$ is idempotent (should point to it now! we did not do this even
  in the previous draft, but it is a long way earlier in
  the paper) and $G_a$ more
  generally (or at least under different conditions).
  In what follows I deduce we are using the former?
And that it agrees with $G_a$, where defined.
  But I
  think this should be stated... particularly since it introduces a
  condition that we sometimes like to argue about!}
By \ref{pa:freerightS} the $P^l_n-\kiss$ bimodule $P_n^\mm a^\mm$ is
free as right $\kk S_\mm$-module.
Now in Lemma~\ref{lem:FM} part III let   $F$  be $P^\mm_na^\mm$
and $ M$ be the left ideal $\mmu$-Specht module $\sspecht_{\lmu}$
to obtain the desired isomorphism.
In particular the
basis of $ \Delta_{\lmu}$ 
of form
$$
B^{\lmu}_{\Delta} = \{ p \oxm w \; | \; p \in T^\mm , w \in {b_{\lmu}} \}
$$
(Lemma~\ref{lem:FM}(I);
we write $\oxm$ here to distinguish from the other types of tensor product in
this Section)
is taken element-wise to $B^{\lmu}_\sspecht$.
\qed



{\mlem\label{lem:step3}{
Let $\delta\neq0$ and $k=\mathbb{C}$.
Let $\mm\in\gammal$.
As a left $P^l_n$-module
the quotient $P^\mm_na^\mm P^\mm_n$
{  is a direct sum of $\Delta$-modules}. 
}}
\begin{proof}
By Lemma~\ref{lem:tranv1} we have the following isomorphism of left $P^l_n$-modules
\beq \label{eq:PaP1}
P^\mm_na^\mm P^\mm_n =
\bigoplus_{p \in T^\mm}P^\mm_na^\mm p^\star
\simeq \bigoplus_{T^\mm}P^\mm_na^\mm
\eq
where $p^\star$ is as defined in \ref{de:flip}. 
Furthermore   Lemma~\ref{le:Diso} implies that 
\begin{equation} \label{eq:PaP2}
P^\mm_na^\mm\simeq  P^{\mm}_n a^{\mm} \otimes_{\kiss}  {\kiss}  \simeq
\bigoplus_{\underline\mu\in\Lambda_\mm}\left(P^\mm_na^\mm e_{\underline\mu}\right)^{dim(\sspecht_{\lmu})}
\end{equation}
    as a $P^l_n$-module.
\end{proof}


\mdef 
Comparing \ref{lem:idbas1}, \ref{lem:polarPm} and
\eqref{eq:PaP2}
we see that
\beq \label{eq:sumsquares}
\dim(P^l_n ) = \sum_{\mmu \in \Lambda(P^l_n)} \dim( \Delta_{\mmu} )^2 
\eq
(cf. a cellular basis of $P^l_n$ in the sense of \cite{GrahamLehrer96}).



{{\mlem \label{th:weaktri}
Suppose $k=\C$
and $\delta \neq 0$.
Recall the simple $P^l_n$-modules $L_{\lmu} = \head\specht^n_{\lmu}$
from (\ref{de:Llmu}). 
The modules
$\{ \specht^n_{\lmu} \; : \; \lmu \in \Lambda(P^l_n) \}$
have a 
lower-unitriangular decomposition matrix
$([\specht^n_{\lmu} : L_{\lnu}])_{\lmu,\lnu \in \Lambda(P^l_n)}$
with respect to 
any order $\prec$ on  $\Lambda(P^l_n)$ in which
$\mm<\mm'$  implies
$\lnu \prec \lmu$
for
$\lmu \in \Lambda_\mm$ and $\lnu\in\Lambda_{\mm'}$.
}}
\proof
By  Theorem~\ref{th:index1} the pair $(\gammal , a^- )$
is a core of  $P^l_n$.
When $k=\C$ for each $\lmu\in \Lambda_\mm$ the $\mmu$-Specht module
$\sspecht_{\mmu}$  is simple as $\C S_\mm$-module.
By Equation~\ref{eq:spechtinflate} and Lemma~\ref{le:Diso}
the set of modules
$\{\specht^n_{\lmu}\mid \lmu \in \Lambda(P^l_n)\}$
are the corresponding
\redxx{ what does this mean: }
long $G_a$-functor
$\Delta$-modules as in Equation~\ref{eq:Gdelta}.
Now the result follows from  Lemma~\ref{Th:upperun}.\qed


\section{Globalisation of standard modules}
\newcommand{\W}{{\barr{W}}}
\newcommand{\pT}{{\mathsf T}}

In this section we study the effect of the $G_\W$ functor on standard
modules. That is we study the $P^l_n$-module
$G_\W \; \specht^{n}_{\lmu} $.
This is particularly interesting because the non-trivial core property
leads to some new departures from the partition algebra argument.
The main result is Proposition~\ref{lem:Gspecht} below.



\mdef
Recalling  $T^\mm_{}$, 
let $\pT^\mm_n$ denote  the 
set of representative relative-non-crossing
partitions $p$ rather than the classes $p+I^{<\mm}$.

\newcommand{\Deltasm}{\Delta_{\lmu}}



\mdef NB our convention is that if $R$ is a ring and $S$ a set then
$RS$ generally denotes the free $R$-module with basis $S$. However if
$S$ is given as a subset of an $R$-module $M$ then $RS$ means the
$R$-span of $S$ in $M$.


\mdef
{\mupr { \label{lem:Gspecht}
    Let
    $\mm\in\gamma^l_{n-l}$ and $\lmu\in\Lambda_\mm$.
Applying $G_\W$  from (\ref{de:Gfunctor}),
consider the subset of $G_\W \; \specht^{n-l}_{\lmu} $ given by
$B^{\lmu}_G = 
\{ t \ox a^\mm_{n-l} \www\omega \; | \; t \in \pT^\mm_n,  \omega \in b_{\lmu} \}$ 
where
we use $\ox$ for the $G_\W$ tensor product.
Then $B^{\lmu}_G$
is a basis of $G_\W \; \specht^{n-l}_{\lmu} $
and
\redxx{
\\ (II) It induces an identical representation to that induced by 
$B^{\lmu}_{\sspecht}$ in (\ref{lem:basisBSp}).
That is,
for $b \in B=B^{\lmu}_\specht$,
$\;\psi(b)$ its natural image in $B^{\lmu}_G$,
and $ab = \sum_{b'\in B} \alpha^a_{b,b'} b'$,
we have $a \psi(b) = \sum_{b'\in B} \alpha^a_{b,b'} \psi(b')$.
\\ (III) Thus
}
\[
\specht^n_{\lmu} \;\; \cong \;\; G_\W \; \specht^{n-l}_{\lmu} \; .
\]
}}
\proof

\noindent
Note from the construction that
$ P^l_n\W\specht^n_{\lmu}= \specht^n_{\lmu}$.
By Lemma \ref{lem:FWS} 
$G_\W \; \specht^{n-l}_{\lmu} \cong P^l_n\W \ox \W \specht^n_{\lmu} $.
For the latter form we have 
a multiplication map $p\otimes s \mapsto ps$ to
$P^l_n\W\specht^n_{\lmu}= \specht^n_{\lmu}$.
This gives a surjective $P^l_{n}$-module homomorphism. 
Since $| B^{\lmu}_{G} |=|B^{\lmu}_{\sspecht}  |$ it is enough to show that
$B^{\lmu}_{G} $ is spanning.
The basis $ B^{\lmu}_{\sspecht}$ gives
\beq \label{eq:run2u}
\specht^n_{\lmu} = P^\mm_n a^\mm_n w^{ba}_{e_\mu}
    = \kk \{  t \www\omega \;|\; t \in T^\mm_n; \; \omega \in b_{\lmu}  \}
\eq
This holds for any $n$ but, as indicated here
the non-crossing transversal $T^\mm_n$ in $P^\mm_n a^\mm_n$ of course depends on $n$.
Applying $G_\W$ to the case with $n$ replaced by $n-l$
we have
\beq \label{eq:run3u}
G_\W \; \specht^{n-l}_{\lmu} =
P^l_n \W \ox \specht^{n-l}_{\lmu}    
    =
    \kk \{ P^l_n \W \ox t \www\omega \;|\; t \in T^\mm_{n-l}; \; \omega \in b_{\lmu}  \}
\eq \[ \hspace{2.1in}
=
    \kk \{ d \W \ox t \www\omega \;
            |\; d \in \mP^l_n; \; t \in T^\mm_{n-l}; \; \omega \in b_{\lmu}  \}
\]
\redxx{more bars to do!
(NB. For now we pretend that $W$ has its $ww^*$ tensor factor at the end
rather than the begining!)}
Let us `move' $t$ through the tensor product in \eqref{eq:run3u}:
\[ \hspace{1in}
    = \kk \{ d \W ((t+I^{<\mm}_{n-l})  \otimes 1_l)
             \ox a^\mm_{n-l}  \www\omega 
               \; | \; d \in \mP^l_n; \; t \in \pT^\mm_{n-l}; \; \omega \in b_{\lmu}  \}
\]
\beq \label{eq:run4u}  \hspace{1.6in}
    = \kk \{ d  
         (t \otimes ww^*)
            \ox a^\mm_{n-l} \www\omega 
              \; |\; d \in \mP^l_n; \; t \in \pT^\mm_{n-l};  \omega \in b_{\lmu}  \}
\eq
Note the recasting of $t$.
We can omit the $+ I^{<\mm}_{n-l} $ since it does not affect the element. 
We aim to show that this is spanned by terms of form
$
  t \ox a^\mm_{n-l} \www\omega 
$
where $t \in \pT^\mm_n$.
\redxx{Right!? Hmm. I am still not sure
this is exactly right! I think the $ww^\star$ should not be there.} 
Note that we may assume the module is generated by elements of the
given form. 
We proceed as follows.
Consider the action of the generators from 4.9 on
an element of the claimed spanning
set.
Let $d$ be such a generator and consider first the `factor'
$dt$ in $dt   \ox a^\mm_{n-l} \www\omega$.
Noting Lemma~\ref{lem:1}(II)
there are two ways in which $dt$ might pass out of the
relative noncrossing transversal: 
either (A) it has a component with lower propagating index,
i.e. a component in $I^{<\mm}$;
or (B) a relative crossing is introduced.


\redx{
In Case (A): by definition such a component of $dt$ is spanned by elements of
the form $s a^{\mm'}_n q$ with $\mm' < \mm$.
Indeed since $\mm\in\gamma^l_{n-l}$ we have $t \propto t \barr{W}$ and
so $dt$ is spanned by elements of the form
$s a^{\mm'}_n q (1_{n-l}\otimes ww^\star)$.
Since $ a^{\mm'}_n = a^{\mm'}_{n-l}\otimes ww^\star$
$=  (a^{\mm'}_{n-l}\otimes 1_l) (1_{n-l}\otimes ww^\star)$
this becomes
$ s (a^{\mm'}_{n-l}\otimes 1_l)
 (1_{n-l}\otimes ww^\star) q (1_{n-l}\otimes ww^\star)$.
Noting that
$$
(1_{n-l}\otimes ww^\star) q (1_{n-l}\otimes ww^\star)
\propto \;   q|_{[1,n-l]} \otimes ww^\star ,
$$
we see that
$
dt \propto \;\;
s(( a^{\mm'}_{n-l}  q|_{[1,n-l]}) \otimes ww^\star ) 
$.
But
$$
(( a^{\mm'}_{n-l}  q|_{[1,n-l]}) \otimes ww^\star ) \ox  a^\mm_{n-l} \www\omega
\propto\;\;
( a^{\mm'}_{n-l}   \otimes ww^\star )
     \ox   (a^{\mm'}_{n-l}  q|_{[1,n-l]}) a^\mm_{n-l} \www\omega
\;\; =0 
$$
so $dt\ox  a^\mm_{n-l} \www\omega =0$.
\\
In Case (B): by the $k S_{\mm}$-freeness property
(Lemma \ref{lem:tranv1}) such a
crossing may be factored out and passed through the tensor product.
}

Thus neither case takes us out of the span, and we are done.
\qed


\medskip

We conclude this section with 
a remark
on our working assumptions. 

\mdef
Consider  bases as in the proof of Lemma~\ref{le:Diso}, 
\redxx{(the bases are a little hidden there! what did this used to point to? )}
noting that $G_W G_{a^\mm_{n-l}} \cong G_{a^\mm_n}$.
That is,
$P^l_n W^l \otimes_{P^l_{n-l}} P^\mm_{n-l} a^\mm_{n-l} \cong   P^\mm_n a^\mm_n$,
since $a^\mm_{n-l} \otimes ww^\star = a^\mm_n  $
by (\ref{eq:amm}).
\redx{
  At this point the case $\mm=0$, $\delta=0$ has
  extra interest.
}


Note that  we cannot apply
$P^l_nW^l_b\otimes_{P^l_{n-l}} -$
to $\specht^0_{\lmu}$ since $W^l_b$ requires $n>l$.
In this case we could attempt to use $W^l$ instead.
This works straightforwardly if $\delta$ is a unit.
But if $\delta=0$ 
then the setup is slightly but interestingly different.
Of course our mechanism for making $ P^\mm_l a^\mm_l $ a right $P^l_0$
module does not work.
And $W^l$ is not normalisable as an idempotent
(although the functor
given by allowing ${P^l_{n-l}} $ to act on the right of
$P^l_nW^l_b$ by restriction
is still well defined).
It is an interesting exercise to see what
happens if we simply allow $P^l_0 \cong \kk$ to act as $\kk$.
In that case we are comparing
$G_W 
G_{a^\mm_{n-l}} $ with the direct long functor $G_{a^\mm_{n}}$.
We will leave this for a separate work.



\section{Properties of standard modules}
\subsection{Standard module contravariant form} \label{ss:cvf}

\begin{figure}
\includegraphics[width=5.76543in]{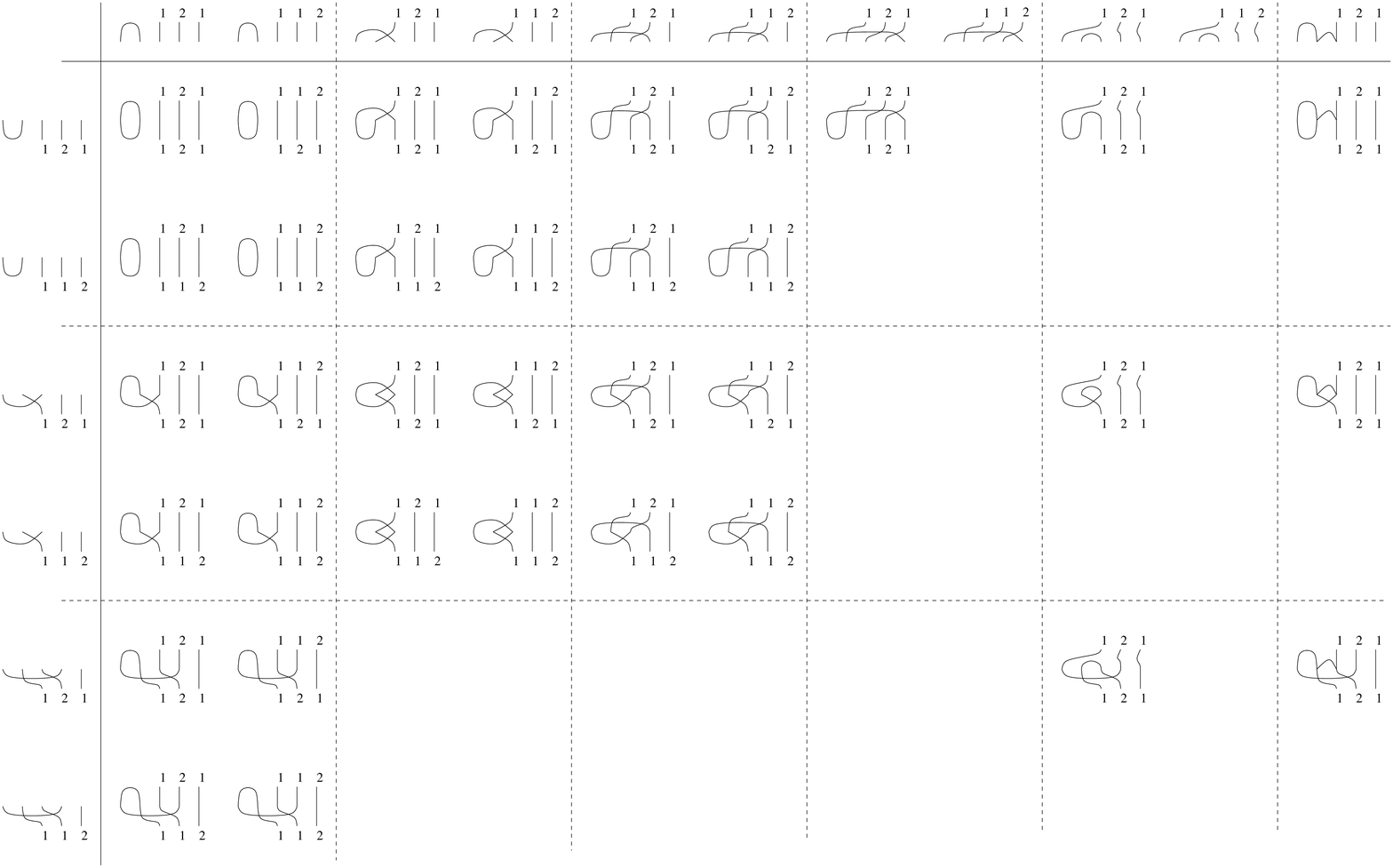}
\caption{Example partial gram matrix calculation. 
The table shows part of the 
case $n=5$, $(\lambda^1,\lambda^2) = ( (2,1) , \emptyset )$.
We use the tableaux basis for $(2,1)$: 
$\{ 112, 121 \}$ \cite{James}.
\label{fig:gramex1i}}
\end{figure}


\begin{figure}
\[
\includegraphics[width=6cm]{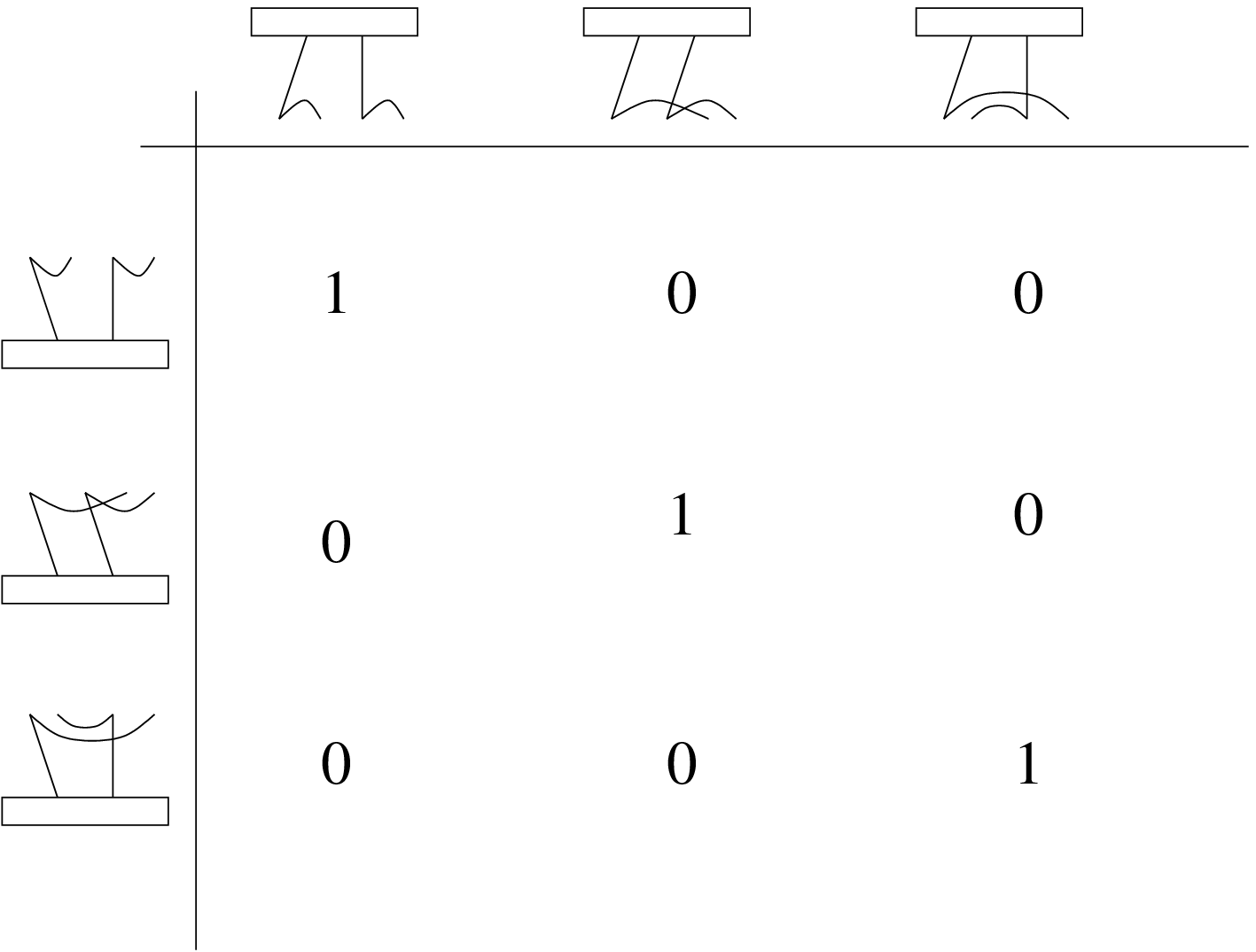}
\qquad
\includegraphics[width=7.496cm]{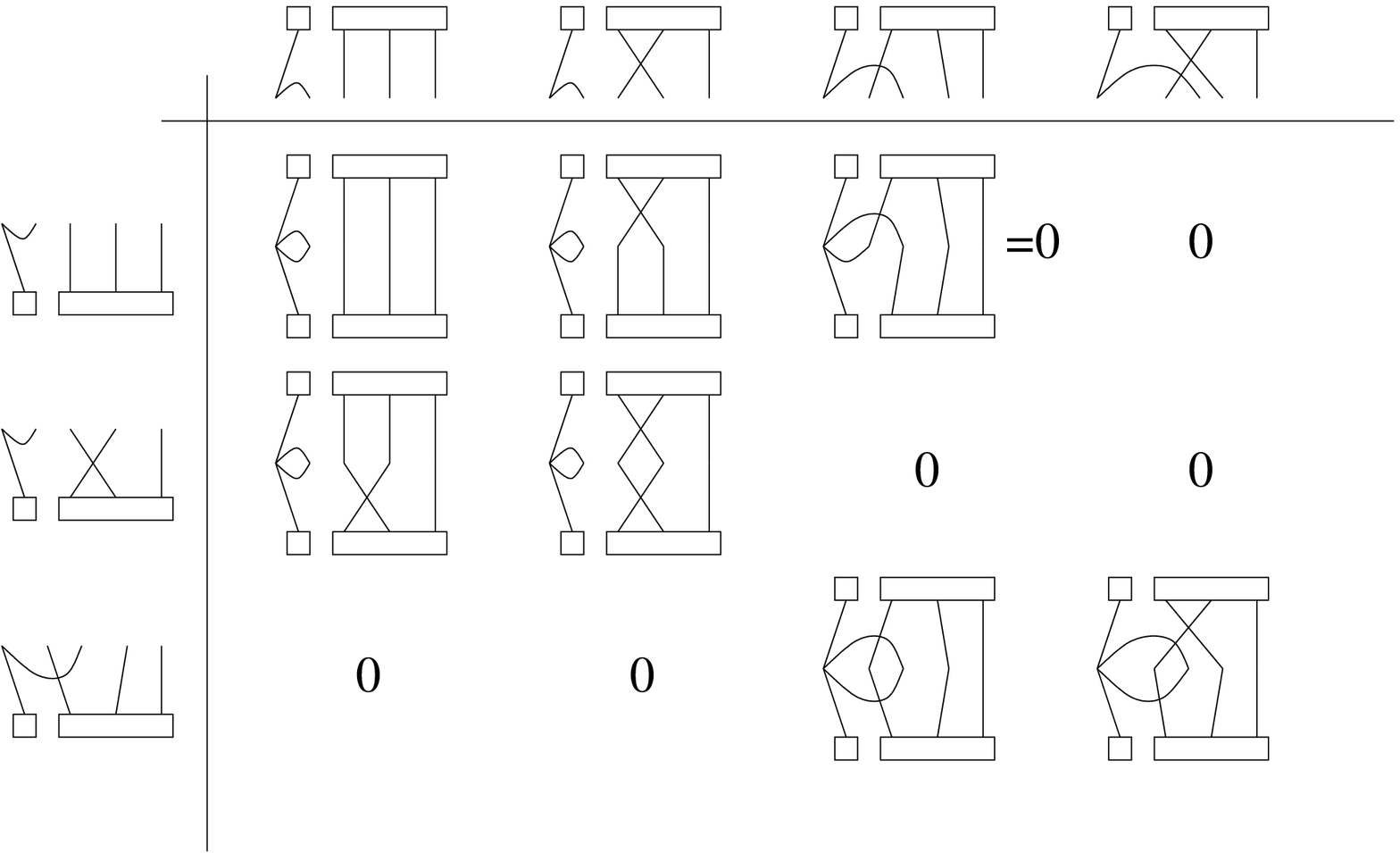}
\]
\caption{ \label{fig:exa:grams}
(a) Here we have $\lambda^1 = \emptyset$ and $\lambda^2 =(2)$ or $(1^2)$ in the box (the picture is 
effectively  the same in each case).
(b)
$\specht_{((2,1),(1))}^5$ 
is 20-dimensional. The top left-hand corner of
the gram matrix is shown.
Here we have $\lambda^1 =(2,1)$ in the long box, with a choice of idempotent
from $kS_3$ such that $e$ and $\sigma_1 e$ span $kS_3 e$.
}
\end{figure}


{\mlem{ \label{lem:ndgfz0}
    There is a contravariant form $(-,-)_e$ on each
    $\specht^n_{\lmu}$
defined by 
\beq \label{eq:f12}
\left( x a^{\mm} e_{\lmu} \right)^{op} \; y a^{\mm} e_{\lmu} 
\;\; = \;\;  \left( e_{\lmu} \right)^{op} a^{\mm} x^{op} y 
   (a^{\mm} e_{\lmu})
\eq \[ 
\;\;\;\; = \;\;  \left( e_{\lmu} \right)^{op} (a^{\mm} x^{op} y 
   a^{\mm})( e_{\lmu})
\;\; = \;\; (x,y)_e \; a^{\mm} e_{\lmu}
\]
}}
\proof 
(For (\ref{lem:ndgfz0}) and (\ref{lem:ndgfz})
we have a direct generalisation of the usual partition algebra
argument as in, for example,  \cite{Martin96}.) 
Note that we are working in $P^\mm_n$. 
The well-definedness of the form follows from the construction
as in (\ref{de:gcvf}), 
using Lemma~\ref{lem:groupc} and primitivity at the last step.
\qed

Examples: See Fig.\ref{fig:gramex1i} and \ref{fig:exa:grams}.


{\mlem{ \label{lem:ndgfz}
For $k = \Z[\delta]$ the
determinant of the gram matrix of the form $(,)_e$
on $\specht^n_{\lmu}$ ($\lmu \in \Lambda(P^l_n )$)
is nonzero.
}}
\proof 
First organise the basis into blocks according to
Lemma~\ref{lem:basisBSp}  
--- i.e. each block has a fixed non-crossing partition $p$, with only the
permutation group module basis part $w$ varying. See Fig.\ref{fig:gramex1i}.
Note: 
(QI) If we work for the moment over $\Q [ \delta]$ 
(as it will be clear that we can in investigating the nonzero
property) 
then $k S_\mm$ is (split) semisimple \cite{James} and 
we may use a basis for the permutation group part in which the gram
matrix of this part is diagonal.
We could, for example, use tableau bases \cite{James},
as in (\ref{de:blambdabasis}),
as illustrated in the figure.
The numerical {\em details} 
of this part of the construction will not be needed here.

For any given choice of ordered basis we arrive at the gram matrix, 
denoted $G^n_{\lmu}$. 
In our $ ((2,1),\emptyset) $ example we have 
\[
G^n_{((2,1),\emptyset)}  \; = \; \mat{c|c|c|c|c|ccc}
\delta G & G       & G'  & G'' & G & G & ... \\ \hline
G        &\delta G & G   & ... & G & G & ... \\ \hline
G'       & ...     & ... & ... & 0 & 0 & ... \\ \hline
...      & ...     & ... & 
\tam 
\]
where the entries shown are the block submatrices. 

In every row of the gram matrix proper, every entry is a
polynomial in $\delta$, indeed an integer multiple of a power of $\delta$. 
The diagonal entry of the basis element $pw$
(say --- cf.
Lem.\ref{lem:basisBSp}
and its proof),
determined by $(w^{op} p^{op}) p w = w^{op} (p^{op} p) w$, 
is a nonzero polynomial
whose degree is not exceeded by any other entry. 
The degree is the same through a given block; and there is at least
one row where the block of the diagonal has strictly the highest of
all degrees in the row.
NB, In case $\lmu \in \mm \in h^l_n$ all entries are constant, but all
the off-diagonal blocks are zero by the $I^{< \mm}$ quotient. 
Finally, by (QI) the blocks on the block-diagonal are diagonal. 
Combining these we see that the gram matrix has full rank for
indeterminate $\delta$. Thus the determinant is a nonzero polynomial
in $\delta$.
\qed


\mdef \label{lem:Stanss}
 \mupr 
(I) The standard modules 
$\specht_{\lmu}$ with
$\lmu \in \sqcup_{\mm \in h^l_n} \Lambda_\mm$
are $\Al_n$-modules.
(II) Over $k = \C$ every such   $\specht_{\lmu}$  is simple;
and $A^l_n$ is semisimple.
\proof
(I) $W \specht_{\lmu} = 0 $ iff $| \lmu | \in h^l_n$. 
(II) The contravariant form is non-degenerate here.
Indeed, by the compatibility condition (\ref{de:compatible}),
$(x,y)$ is only non-zero if `row-parts' 
\redxx{(was pairs - fix)}
match, whereupon the gram matrix consists of blocks
corresponding to matches. Within these blocks the entries are the same
as for the gram matrices for the (product of) symmetric groups
--- which gram matrices are of full rank over $\C$ \cite{James}.
Finally, 
by (\ref{de:Llmu}) 
we have completeness and pairwise nonisomorphism.
 \qed



\subsection{Standard module decomposition matrix properties} 
\label{ss:uut}

\redxx{Delete one of the following 2 Theorems!: } 

\redxx{Paul version:}

{\mth{ \label{de:remcx}
Suppose $k=\C$
and $\delta \neq 0$.
Recall the simple $P^l_n$-modules $L_{\lmu} = \head\specht^n_{\lmu}$
from (\ref{de:Llmu}). 
The modules
$\{ \specht^n_{\lmu} \; : \; \lmu \in \Lambda(P^l_n) \}$
have an upper-unitriangular decomposition
matrix
$([\specht^n_{\lmu} : L_{\lnu}])_{\lmu,\lnu \in \Lambda(P^l_n)}$
with respect
any order
$(\Lambda(P^l_n ) , \preceq)$
in which
$\lmu \prec \lnu$
if
$r_{ |\lmu |} < r_{| \lnu |}$. 
}}
\proof
This follows from Proposition~\ref{lem:Gspecht},
Proposition~\ref{lem:Stanss}
and the construction
using Lem.\ref{lem:unitri0} and Th.\ref{th:index1}. 
\qed


\mdef Example: If $r_{| \mu |} = n$ then $\specht^n_{\mu} = L_\mu$
(since there is no $\nu \in \Lambda(P^l_n)$ with $r_{ | \nu |} > n$).
Meanwhile $(\specht_0 : L_0) =1 $ and no other composition factor is
precluded for $\specht_0$ by this Theorem
(and indeed none can be without specifying $\delta$).



{\mth{ \label{th:ssss}
Consider $k=\C$ and $\delta \in \C$. 
Each module $\specht_{\lmu}$ is simple for all but finitely
many values of $\delta$. 
}}
\proof 
By 
(\ref{lem:unitri0})
it is enough to show that there is a nondegenerate
contravariant form on each module
(if one of these modules is isomorphic to its contravariant dual then it
contains the dual of the head $L$ in the socle; but if these are not the
same module then they are not isomorphic,
by the unitriangular property
\redx{Th.\ref{de:remcx}}; 
and by completeness there is another such module
with head $L^o$ and socle $L$, contradicting upper-triangularity). 
Now note Lem.\ref{lem:ndgfz}.
\qed



{\mth{\label{th:ssssc}
Consider $k=\C$ and $\delta \in \C^*$. 
For all but finitely many values of $\delta$: 
(I) The set  
$\{ \specht_{\lmu} \; | \; \lmu \in \Lambda(P^l_n) \}$ 
is  a complete set of  simple modules of $P^l_n$;
(II) $P^l_n$ is semisimple. 
}}
\proof (I) These modules are (sufficiently often) simple by Th.\ref{th:ssss}.
By \ref{de:Llmu} it is a complete set.
By the embedding property using the functors
$F,G$
there are no duplicates (pairwise isomorphisms) in the set. 
\\
(Alternatively we may argue using the dimension count (\ref{eq:sumsquares})
and either the completeness or the pairwise nonisomorphism.)
\\
(II) follows from (I) and (\ref{eq:sumsquares}).
\qed


\mdef
By Theorem~\ref{th:ssssc}
and the construction 
we have, over $\C$, a modular system
for each $\delta$ 
(see e.g. \cite[\S1.9]{Benson95}).
Specifically we may take $K= \C[x]$ for the integral ground ring;
the field of fractions $K_0$ as the ordinary case; and $\C$ with $x$
evaluated at $\delta$ as the modular case.
Thus we have
Brauer--Humphreys reciprocity:
\beq \label{eq:BH1}
(P_{\lambda}:\specht_{\mu} ) = [ \specht_{\mu} :L_{\lambda} ]
\eq
where $(P_{\lambda}:\specht_{\mu})$ denotes the `composition multiplicity' of 
$\specht_{\mu}$ in $P_{\lambda}$ 
(as usual this makes strict sense over the rational field
via an idempotent lift, and as a
multiplicity in the Grothendieck group in general). 

In particular if $L_{\lambda}$ 
is a composition factor of $\specht_{\mu}$ then $\specht_{\mu}$ 
is a filtration factor of $P_{\lambda}$.


\section{On quasi-heredity}

Here we prove, in Theorem~\ref{th:qh5}, that the
$P^l_n$ module categories are highest weight
categories
(in the sense of Cline, Parshall and Scott \cite{CPS})
when $\delta \neq 0$ and $k=\C$.
Given Theorem~\ref{th:strongmodular}, Theorem~\ref{de:remcx} and so on,
it is enough to show that projective modules are
filtered by $\Delta$-modules. We do this next. 
Recall the following.

\subsection*{General Lemmas}

\mdef \mulem \cite{Martin0915}
\label{lem:filtfilt}
(I) Let $A$ be an algebra, $M$ an $A$-module and $S$, $T$ sets of
$A$-modules. If $M$ has an $S$-filtration and every $N \in S$ has a
$T$-filtration then $M$ has a $T$-filtration.




{\mlem{ \label{lem:Mfissummand}
    Let $A$ be an algebra, $f$ an idempotent, and $M$ a bimodule.
    Then there are left-module maps $Mf \rightarrow M$ given by
    inclusion;
    and $M \rightarrow Mf$ given by $m \mapsto mf$. Indeed the
    sequence
    \[
0 \rightarrow Mf \rightarrow M \rightarrow M(1-f) \rightarrow 0
    \]
    is short-exact and split.
}}

{\mlem{ \label{lem:gentry}
    Let $A$ be a finite-dimenional algebra and suppose
    $0\subset J_1 \subset J_2 \subset ...\subset J_k =A $
    is a filtration by ideals.
    Let $X$ be the set of indecomposable summands of all the left-modules
    $J_i / J_{i-1}$ up to isomorphism. Then every projective left-module of
    $A$ is filtered by $X$.
}}
\proof
It is enough to show for indecomposable projectives, and hence for
modules of form $Af$ where $f$ is idempotent.
For such a module we have a (possibly degenerate) filtration
   $0\subseteq J_1 f \subseteq J_2 f \subseteq ...\subseteq J_k f =Af $
by assumption.
Suppose $M \supset N$ are $A$-bimodules.
We claim that there is a left-$A$-module map $(M/N)f \rightarrow Mf/Nf$ given by
$\{ mf+n : n \in N \} \mapsto \{ mf+nf : n \in N \}$
(i.e. act on the set elementwise by $f$ on the right),
and that this is a left-$A$-module isomorphism.
To see this note that
(i) the image lies in $Mf/Nf$;
(ii) this gives a vector space isomorphism;
(iii) it is a left-$A$-module morphism
(the map uses the action on the right which, by the bimodule property
commutes with the action on the left).

Thus in particular
$(J_i / J_{i-1})f \cong   
J_i f / J_{i-1} f$,
so there is a sectioning of $Af$ with sections isomorphic to modules
$(J_i / J_{i-1})f $.
By
Lemma~\ref{lem:Mfissummand}
$(J_i / J_{i-1})f $ is a sum of (some) direct summands of an
indecomposable direct summand decomposition of $J_i / J_{i-1}$.
And by  the Krull--Schmidt Theorem for modules,
and our working assumptions,
every such decomposition is a sum from $X$.
The Lemma now follows routinely using (\ref{lem:filtfilt}),
since a direct sum has a filtration by its summands.
\qed




\subsection{Quasiheredity/HWC for the tonal algebras}

{\mlem{ \label{lem:Pfiltered}
    For $\delta\neq0$, the indecomposable projective
    $P^l_n $-modules are  filtered by the set $\Delta^{l,n}$ of
    $\Delta$-modules (as in (\ref{de:std1})). 
}}

\proof \redxx{Using Lemma~\ref{lem:gentry}.}
Consider the ideal chain from (\ref{id-chain}).
The sections are as in Lemma~\ref{lem:step2}.
By
Equations~\eqref{eq:PaP1} and ~\eqref{eq:PaP2}
these are sums of certain modules, and by (\ref{le:Diso}) these
modules are $\Delta$-modules.
Now use Lemma~\ref{lem:gentry}.
\qed


\redxx{DELETE?: Noting (\ref{HW-def})  and \cite[\S x]{CPS} we may now prove:}

{\mth{ \label{th:qh5}
    If  
    $k = \C$ and $\delta\neq 0$ 
    then $P^l_n-\!\!\!\!\mod$ is a highest weight category
    \redx{ with respect to}
any order $\prec$ on  $\Lambda(P^l_n)$ in which
$\mm<\mm'$  implies
$\lnu \prec \lmu$
for $\lmu \in \Lambda_\mm$ and $\lnu\in\Lambda_{\mm'}$,
    with the set of standard modules
    $\{ \Delta_{\lmu} \; | \; \lmu \in \Lambda(P^l_n) \}.$
 }}



\redxx{{
\proof
We will use Lemma~\ref{lemm:modHWC}.
\red{We thus require to show: (I), (II), (III),.......}
By  Theorem~\ref{th:index1} the pair $(\gammal , a^- )$
is a core of  $P^l_n$.
By Equation~\ref{eq:spechtinflate}
the set of modules
$\{ \Delta_{\lmu} \; | \; \lmu \in \Lambda(P^l_n) \}$
are the corresponding $G$-functor $\Delta$-modules.
\red{(Does this mean $\Deltag$?)}
Lemma~\ref{lem:Pfiltered} implies that the indecomposable projectives of
$P^l_n$ are filtered by the $\Delta$-modules.
By Equation~\ref{eq:BH1} the $\Delta$-modules are pivotal set
for the BH property. Now the result follows from
Lemma~\ref{lemm:modHWC}
\red{using Lemma~\ref{th:weaktri}????}.
\qed
}}

\redxx{We also have in mind a PROOF version using the last core theorem
  \ref{th:strongmodular}?
  - where we need the ground ring to be localised at $\delta=0$
  to satisfy the modular core property (so we cannot base change
  to have $\delta=0$ --- and indeed this case is not (quite) a HWC).}

\proof
\redxx{Take 2: } 
We will use Th.~\ref{th:strongmodular}.
We thus require to show (I) strong modular core; (II) projective
filtration.
\\ (I)
Consider $K$ to be the localisation of $\C[\delta]$ at $\delta$, i.e. the ring of
Laurent polynomials.
Note by  Theorem~\ref{th:index1core} that the pair $(\gammal , a^- )$
is a core of  $P^l_n$ over $K$.
We take $K_0$ to be the extension $\C(\delta)$.
Then the system is strong modular by
Theorem~\ref{th:ssssc}.
\redxx{ Is this ok? It is not usual modular system! } 
Note that our localisation does not prohibit base change to $k$ as
required (but does prohibit $\delta=0$).
\redxx{We thus require firstly a
  modular system  $(K,K_0,k)$ for $P^l_n$.} 
\redxx{We thus require to show: (I), (II), (III),.......}
\\ (II)
By Equation~\ref{eq:spechtinflate}
the set of modules
$\{ \Delta_{\lmu} \; | \; \lmu \in \Lambda_0(P^l_n) \}$
are the corresponding 
$\Deltag$-modules.
By Lemma~\ref{lem:Pfiltered}  the indecomposable projectives of
$P^l_n$ are filtered by the $\Delta$-modules.
\qed

\[ \]


\newcommand{\mule}{\underline{\mu}}
\newcommand{\mula}{\mu,\lambda}
\newcommand{\mulax}{\lambda,\mu}
\newcommand{\mux}[2]{(#2,#1)}
\newcommand{\spechtx}[2]{\specht_{#2,#1}}
\newcommand{\rest}{{\mathtt {res}}} 
\newcommand{\ind}{{\mathtt {ind}}}
\newcommand{\EE}{P^2}
\newcommand{\E}{E}

\section{Standard branching rules for 
           $P^l_n \hookrightarrow P^l_{n+1}$} \label{ss:branch}


Consider the algebra inclusion $P^l_n \hookrightarrow P^l_{n+1}$ given by 
$d \mapsto 1_1 \otimes d$ 
(or equivalently $d \mapsto d \otimes 1_1$). 
We give here the 
`standard-module branching rule'  
associated to the corresponding tower of algebras. 

In the generic/semisimple case this is the simple branching rule
--- the edge rule for the Bratteli diagram, and is given by
Kosuda for example in \cite{Kosuda08}. In the semisimple case the
restriction of a simple modules is a direct sum of simple modules. In
our case we cannot expect this. 
We will need to make some preparations. 

\medskip

Here we write
$M = M_1 + M_2$,
or say $M = \bigplus_i m_i M_i $,
if module $M$ has a filtration by 
a set $\{ M_i \}$,
with the indicated multiplicities. This notation does {\em not} 
give the filtration series order.
(In general there may be a
filtration with different multiplicities, but not, say, if the set is a
basis for the Grothendieck group.)

We may write $S_{\mm -\epsilon_j}$ for the subgroup of $S_\mm$ in which factor
$S_{m_j}$ is replaced by $S_{m_j -1}$. 


We first deal with the case $l=2$,
i.e. $\E_n = \EE_n$,
then, more briefly, with the general case. 


\subsection{The case  $\EE_n = \E_n \hookrightarrow \E_{n+1}$}

\ignore{{

Consider the algebra inclusion $E_n \hookrightarrow E_{n+1}$ given by 
$d \mapsto 1_1 \otimes d$ 
(or equivalently $d \mapsto d \otimes 1_1$). 
We give here the 
`standard-module branching rule'  
associated to the corresponding tower of algebras. 

In the generic/semisimple case this is the simple branching rule
--- the edge rule for the Bratteli diagram, and is given by
Kosuda for example in \cite{Kosuda08}. In the semisimple case the
restriction of a simple modules is a direct sum of simple modules. In
our case we cannot expect this. 
We will need to make some preparations. 

Here we write $M = M_1 + M_2$ if module $M$ has a filtration by 
a set $\{ M_i \}$,
with the indicated multiplicities. This notation does {\em not} 
give the filtration series order.
(In general there may be a
filtration with different multiplicities, but not, say, if the set is a
basis for the Grothendieck group.)

We may write $S_{\mm -\epsilon_j}$ for the subgroup in which factor
$S_{m_j}$ is replaced by $S_{m_j -1}$. 

}}


\begin{figure}
\[ \includegraphics[width=4.3cm]{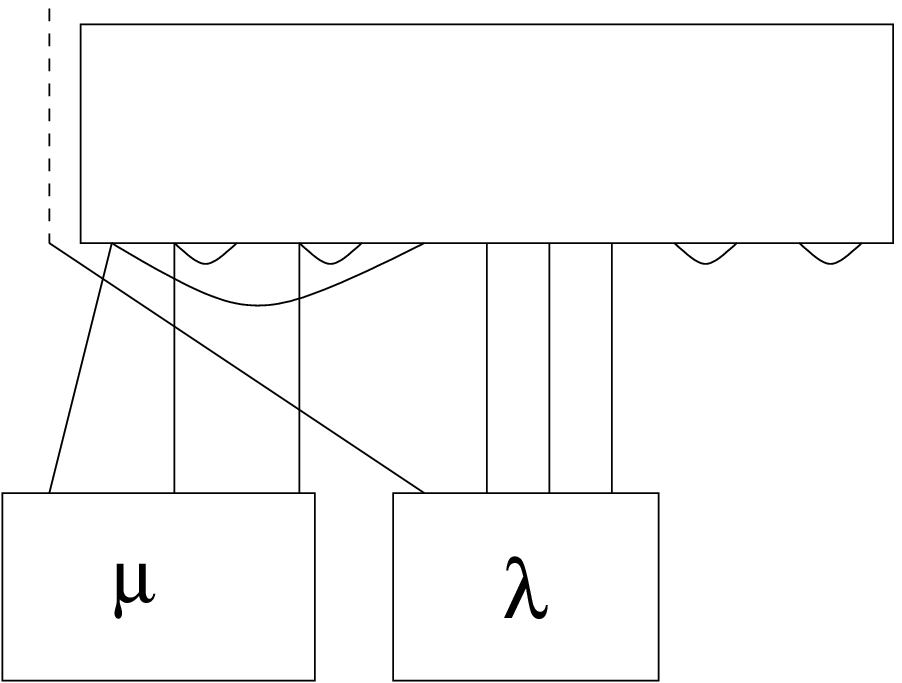} 
\qquad
\includegraphics[width=4.3cm]{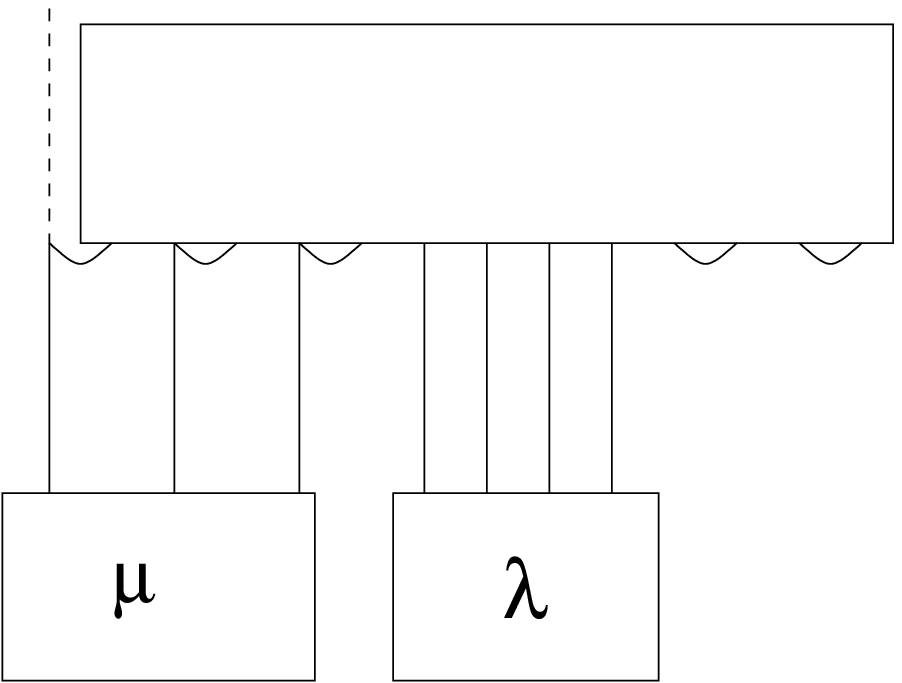} 
\qquad
\includegraphics[width=4.3cm]{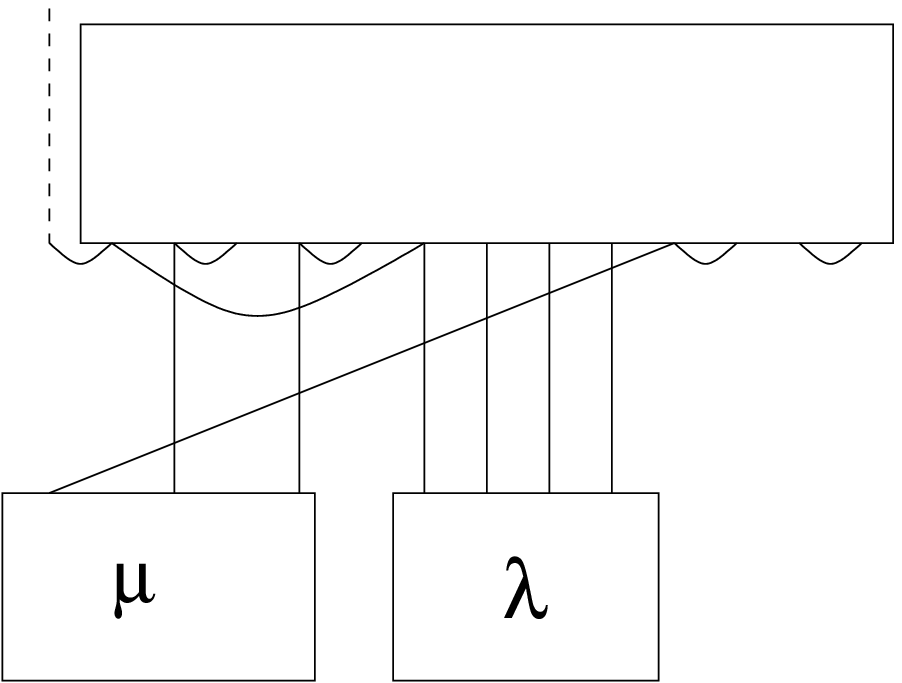} 
\]
\caption{Restricted action on basis elements: cases (1-3).
The undecorated box represents the action of the given inclusion 
$P^2_{n-1} \hookrightarrow P^2_n$ on the element. 
\label{fig:Xres1}}
\end{figure}


\mdef
Set $\mule = \mux{\mu}{\lambda}$.  
Consider the set $B^{\mule}_{\ttSp}$ of 
basis diagrams of $\specht_{\mux{\mu}{\lambda}}$.
(Because of the $e_\lambda$ and $e_\mu $ in the construction these
`diagrams' are linear combinations of partitions in general, but the
well-definedness of the following manipulations will be clear.)
We may organise the set $B^{\mule}_{\ttSp}$ into four 
subsets $ B^{\mule}_{\ttSp}(i)$, $i=1,2,3,4$, 
containing diagrams in which the first vertex is: 
\\
(1) a `propagating singleton'
(Fig.\ref{fig:Xres1} (1)); 
\\
(2) part of a 
co-2 part (Fig.\ref{fig:Xres1} (2)); 
\\
(3) part of a propagating part of higher odd half-order
(Fig.\ref{fig:Xres1} (3)); 
\\
(4) part of a non-propagating part (necessarily of even order). 
\medskip


For example, for $\specht^3_{((1),0)}$ we have $B^{}_{\ttSp}(2) = \emptyset$,
\\
$B^{}_{\ttSp}(3)= \{$
\rb{-.121in}{\includegraphics[width=.561in]{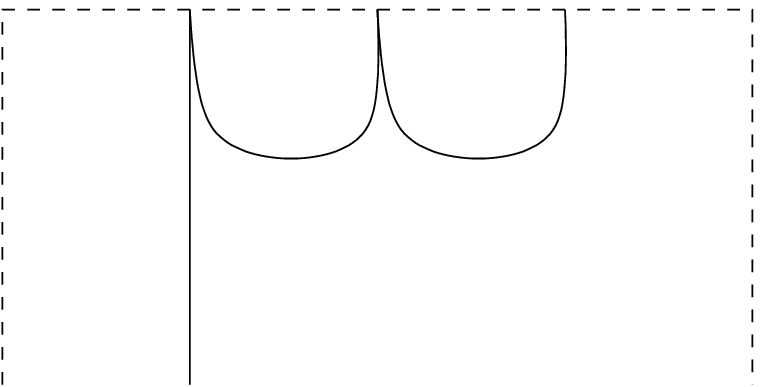}}
$\}$;
$B^{}_{\ttSp}(1)= \{$
\rb{-.121in}{\includegraphics[width=.561in]{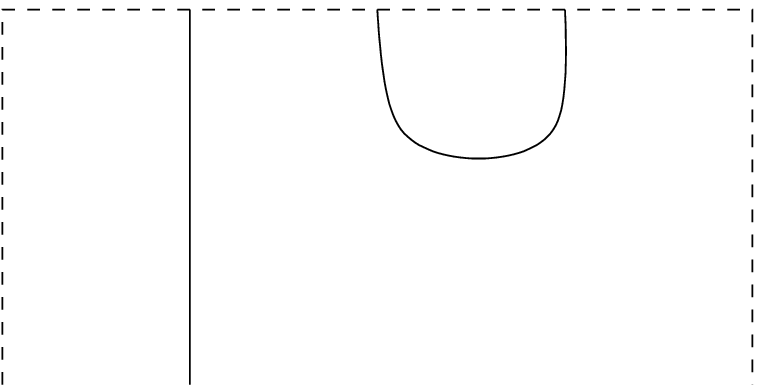}}
$\}$; and
$B^{}_{\ttSp}(4)= \{$
\rb{-.121in}{{\includegraphics[width=.561in]{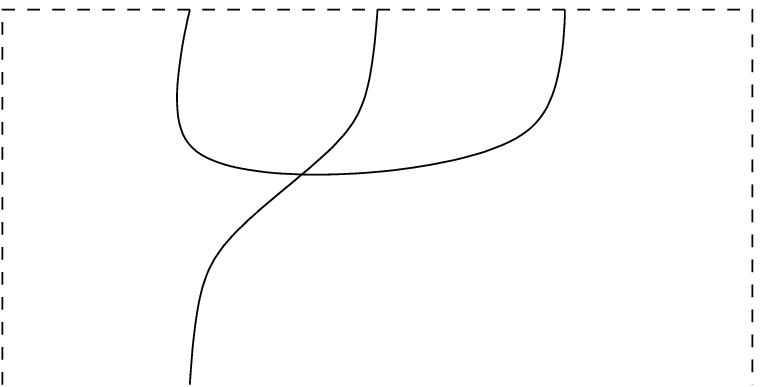} }} , 
\rb{-.121in}{\includegraphics[width=.561in]{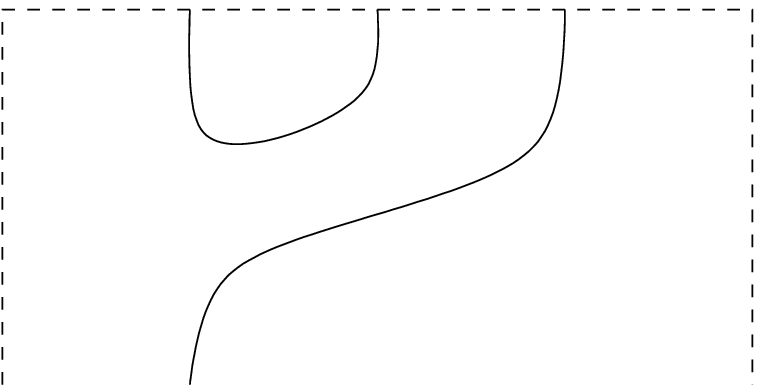}}
$\}$.

\mdef \label{pr:hidden1}
Consider the action  of $1_1 \otimes \E_n$ on
$\specht^{n+1}_{\lmu}$.
\\
 (I) Since the action of $1_1 \otimes \E_n$ is `trivial' on the first
vertex, and so cannot change the propagating singleton property, it 
closes on the $B(1)$ part of
the basis for $\specht^{n+1}_{\lmu}$.
\\
(II) The $B(2)$ part is also closed --- the co-2 property at the first
vertex can only be changed in
principle by combining with a co-1 part, but this gives 0 by the
$I^{<\mm}$ quotient.
\\
(III) On the other hand $B(3)$ is not closed in general, as illustrated by the example:
\rb{-.21in}{\includegraphics[width=.61in]{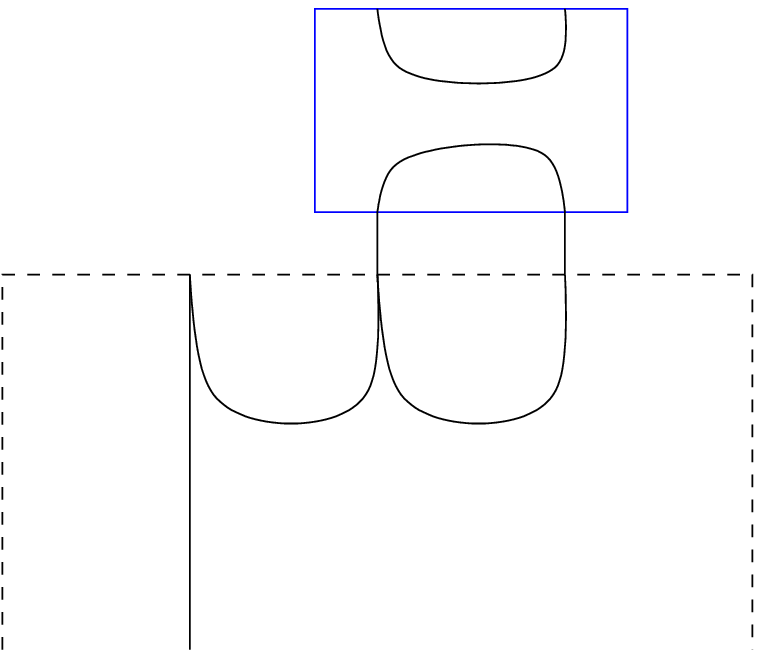}}
which takes $q \in B(3)$ into $B(1)$.
 \medskip \\
(IV) Notice  
however that the subspace spanned by
\beq \label{eq:ASp1}
A_{\ttSp} \; := \; \sqcup_{i=1,2,3} B_{\ttSp}(i)
\eq
is a $1_1 \otimes \E_n$-submodule.



{\mth{ \label{th:brat1} 
Let $\rest_n : E_{n+1}-mod \rightarrow E_n-mod$ denote the
natural restriction corresponding to the inclusion
$E_n    \hookrightarrow E_{n+1}$ given by $d \mapsto 1_1 \otimes d$.
Then $\rest_n \specht_{\mux{\mu}{\lambda}}$ has a filtration by
\redx{standard/Specht}
modules
(as defined in \ref{de:specht}).
The multiplicities are given as follows.
Firstly we have a short exact sequence of $E_n$-modules
\beq \label{eq:resses2}
0 \rightarrow k A_{\ttSp} \rightarrow
  \rest_n \specht^{n+1}_{\mux{\mu}{\lambda}}
    \rightarrow k B_{\ttSp}(4) \rightarrow 0 
\eq
where $A_{\ttSp}$ is the subset of
$ \specht^{n+1}_{\mux{\mu}{\lambda}} $ defined in (\ref{eq:ASp1}).
Then
\[
k A_{\ttSp}  \; = \;
  \bigplus_i \; \specht_{\mux{\mu}{\lambda-e_i}} \; + \; 
   \bigplus_{i,j} \specht_{\mux{\mu+e_i}{\lambda-e_j}} \; + \;  
   \bigplus_{i,j} \specht_{\mux{\mu-e_i}{\lambda+e_j}}  
\]
and
\[
k B_{\ttSp}(4) =   \bigplus_i \; \specht_{\mux{\mu}{\lambda+e_i}}   
\]
\ignore{{
\[
\rest_n \specht_{\mux{\mu}{\lambda}} \; = \; 
   \bigplus_i \; \specht_{\mux{\mu}{\lambda-e_i}} \; + \; 
   \bigplus_{i,j} \specht_{\mux{\mu+e_i}{\lambda-e_j}} \; + \;  
   \bigplus_{i,j} \specht_{\mux{\mu-e_i}{\lambda+e_j}}  \; + \; 
   \bigplus_i \; \specht_{\mux{\mu}{\lambda+e_i}}   
\]
}}
where sums over $\lambda-e_i$, say, denote sums over all ways of
removing a box from the Young diagram of $\lambda$.
}}
\proof 

\noindent
Firstly \redx{(\ref{eq:resses2})}
follows from (\ref{pr:hidden1}). 
For the filtration factors, 
we consider the action of the subalgebra $E_n$ on each of the $ B^{\mule}_{\ttSp}(i)$ in turn. 


\medskip

Case (1): Here we see from Fig.\ref{fig:Xres1}(1)
that the propagating singleton can essentially be ignored.
Each basis element $pw$ (as in (\ref{lem:basisBSp}))
is then like a basis element of 
$+_i \specht_{\mux{\mu}{\lambda-e_i}}$. In particular the `inflation' factor
associated to the $\lambda$ label is the restriction of the Specht
module of $S_{| \lambda|}$ to $S_{| \lambda| -1}$. One then uses the
usual symmetric group restriction rule 
(valid even integrally \cite{Peel75}):
$$
\rest^{S_{|\lambda |}}_{S_{|\lambda| -1}} \ttSp_{\lambda} \;  
     = \;\; \bigplus_i \; \ttSp_{\lambda - e_i}
$$
At the level of bases recall from (\ref{de:blambdabasis})
that $b_\lambda$ is the set of sequence
that are perms of 1112233... (say) satisfying the tableau condition.
Such a perm still satisfies the condition on removing the last
`letter' $i$ say, leaving a basis element in $b_{\lambda-e_i}$.

We thus  have
\[
\bigplus_i \; \specht^n_{\mux{\mu}{\lambda-e_i}} \; 
     \hookrightarrow \; \rest_n \specht^{n+1}_{(\mulax)}
\]
as an injection of $P^2_n$-modules; and a bijection
\beq \label{eq:basiso1}
\bigsqcup_i B^{\mux{\mu}{\lambda -e_i}}_{\ttSp} 
     \stackrel{\sim}{\longrightarrow} B^{(\mulax)}_{\ttSp} (1)
\eq 
The map, on an element $pw$ (as in Lem.\ref{lem:tranv1}), 
is to add a string starting on the left,
then passing over to the box labelled $\lambda-e_i$ and hence to the
part of $w$ that is a tableaux sequence for $\lambda-e_i$, then add 
$i$ to this sequence. 

\medskip


Case (2): Here we can see from the figure that 
(as far as the restricted action is concerned)
the number of propagating
even-half-order components   
effectively goes down by 1 and the number of propagating
odd-half-order components goes up by 1.
Broadly analogously to the previous case we then have a symmetric
group factor with an induction (rather than restriction) on $\lambda$,
and a restriction on $\mu$. 

It will be convenient to have a small complete example to refer
to:  
\[
\underbrace{
\includegraphics[width=.51in]{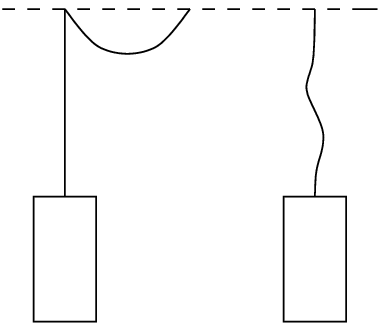}
\qquad
\includegraphics[width=.51in]{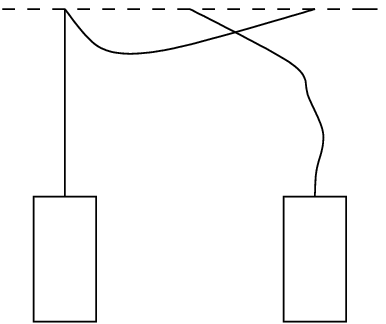}
}_{B_{\ttSp}^{(1,1)}(2)}
\qquad
\underbrace{
\includegraphics[width=.51in]{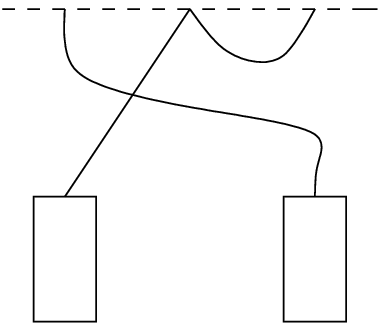}
}_{B_{\ttSp}^{(1,1)}(1)}
\qquad
B_{\ttSp}(3) = \emptyset,
\qquad
B_{\ttSp}(4) = \emptyset
\]
We consider this $\specht^{n+1}_{(1,1)}$ in case $n=2$ 
as a $P^l_{n}$ module
\redx{(N.B. $l=2$)}.
In this case if we quotient
by the $B_{\ttSp}(1)$ submodule, which is isomorphic to
$\specht_{(0,1)}$,  
then the quotient module is already the sum of Specht
modules  
$\specht_{((2),0)} + \specht_{((1^2),0)}$. 
(But we can consider how this is realised. 
We see from (\ref{pr:hidden1})(II) 
that we  
cannot
leave the $B_{\ttSp}(2)$ part.) 

\medskip

Case (3): Here the number of odd lines goes down and the number of even
lines goes up. In the language of the previous cases we have an
induction on one factor and a restriction on another. 
\\
(However, this time the 
construction, (\ref{pr:hidden1})(II), does not preclude a non-split extension. 
And indeed the extension is non-split in general.
To see this consider the example in (\ref{exa:3nonsplit}) below.)

\medskip


Case (4): Here the number of even propagating lines stays the same and the
number of odd lines goes up by 1.
Thus we may put the basis elements in correspondence with an induction
in the odd position. This gives
\beq \label{eq:basiso4}
\bigsqcup_i B^{\mux{\mu}{\lambda +e_i}}_{\ttSp} 
     \stackrel{\sim}{\longrightarrow} B^{(\mulax)}_{\ttSp} (4)
\eq
Now compare with the final summand in the identity in the Theorem.

In this case note that the subset does not span a submodule. The
quotient in the definition of the module
with respect to the $\gamma$ order 
has a slightly different
effect here (just as in the corresponding classical $P_n$ problem
\cite{Martin94}).
It will be convenient to articulate the argument using an example. 

An example of a basis element of $\specht^{n+1}_{\lmu}$ in case 4 is:
$q=$
\raisebox{-.31in}{\includegraphics[width=1in]{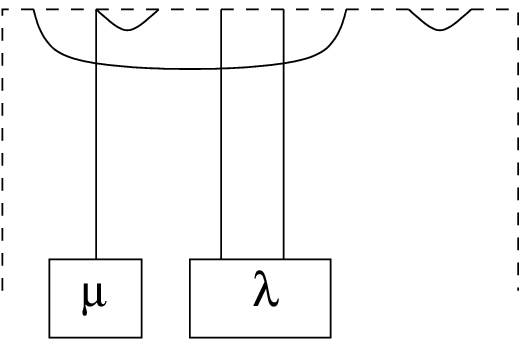}}
.
We consider the action of the $1_1 \otimes P^l_n$ subalgebra as
indicated
by the following schematic:
\raisebox{-.31in}{\includegraphics[width=1in]{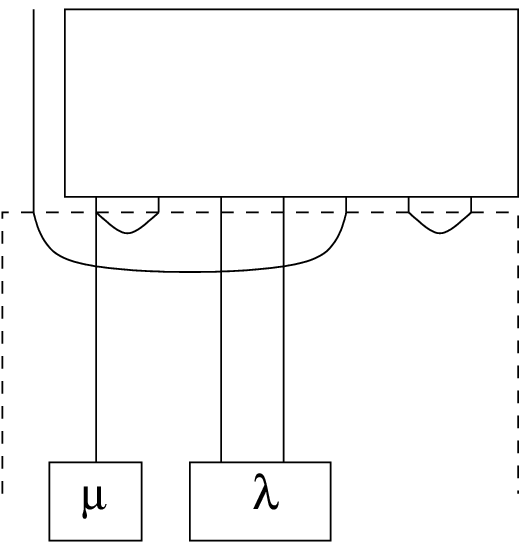}}
.$\;$ 
For example we have, for any completion of the acting partition as
indicated here
\raisebox{-.31in}{\includegraphics[width=1in]{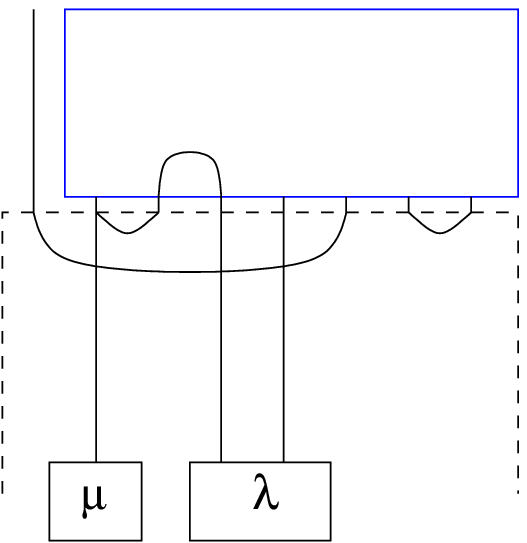}}
$ \equiv 0 $ by the $I^{< \mm}$ quotient.
On the other hand the following has completions for which
\raisebox{-.31in}{\includegraphics[width=1in]{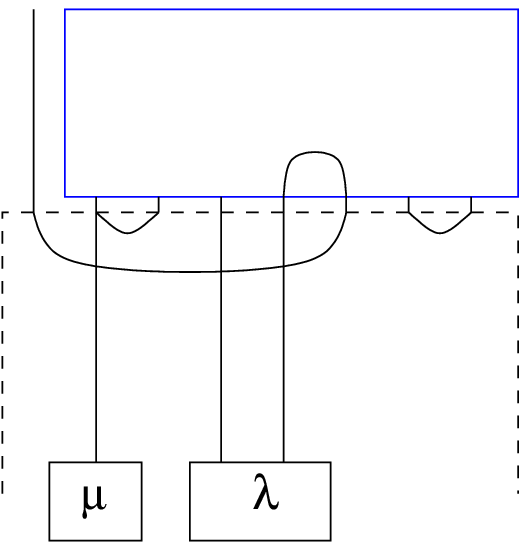}}
$ \not\equiv 0 $, such as
\raisebox{-.31in}{\includegraphics[width=1in]{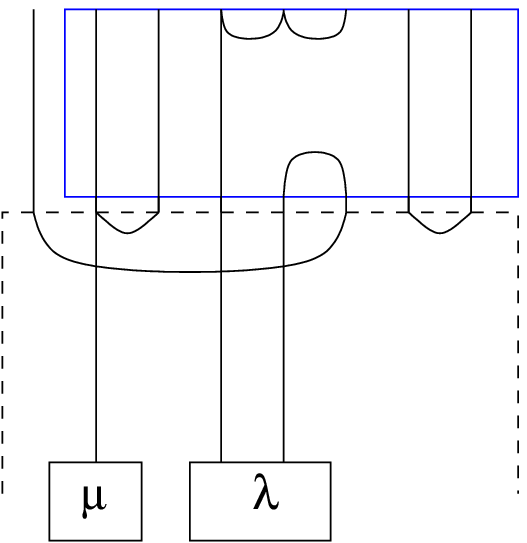}}
.$\;$
This is not 0, but lies in $B^{\lmu}_{\ttSp}(1)$. 
But if $q$ were a basis element of $\specht^n_{(\lambda+e_i ,\mu)}$
under the isomorphism (\ref{eq:basiso4}) then such an element of $P^l_n$ 
would act as 0.
In order to get our filtration here, therefore, we will quotient by a
submodule containing $B^{\lmu}_{\ttSp}(1)$.
Finally consider actions like this:
\raisebox{-.31in}{\includegraphics[width=1in]{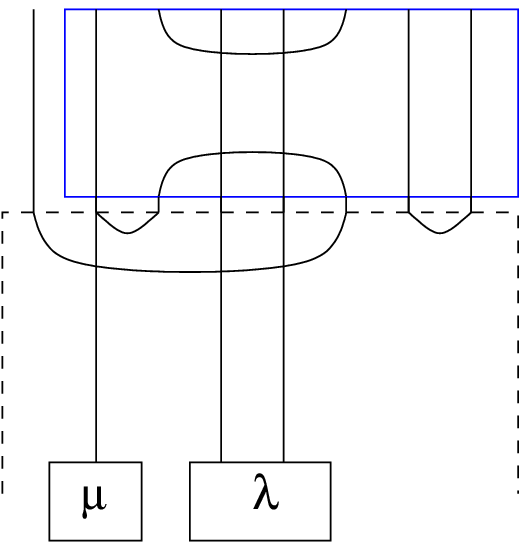}}
.$\;$
Again this is non-zero, lying in  $B^{\lmu}_{\ttSp}(2)$, but should be
0 in this factor of the restriction.
To see that this 
requirement is satisfied  
note
the position of $B(4)$ in
\redx{(\ref{eq:resses2}). }
\qed



\begin{figure}
\[
\includegraphics[width=13.643cm]{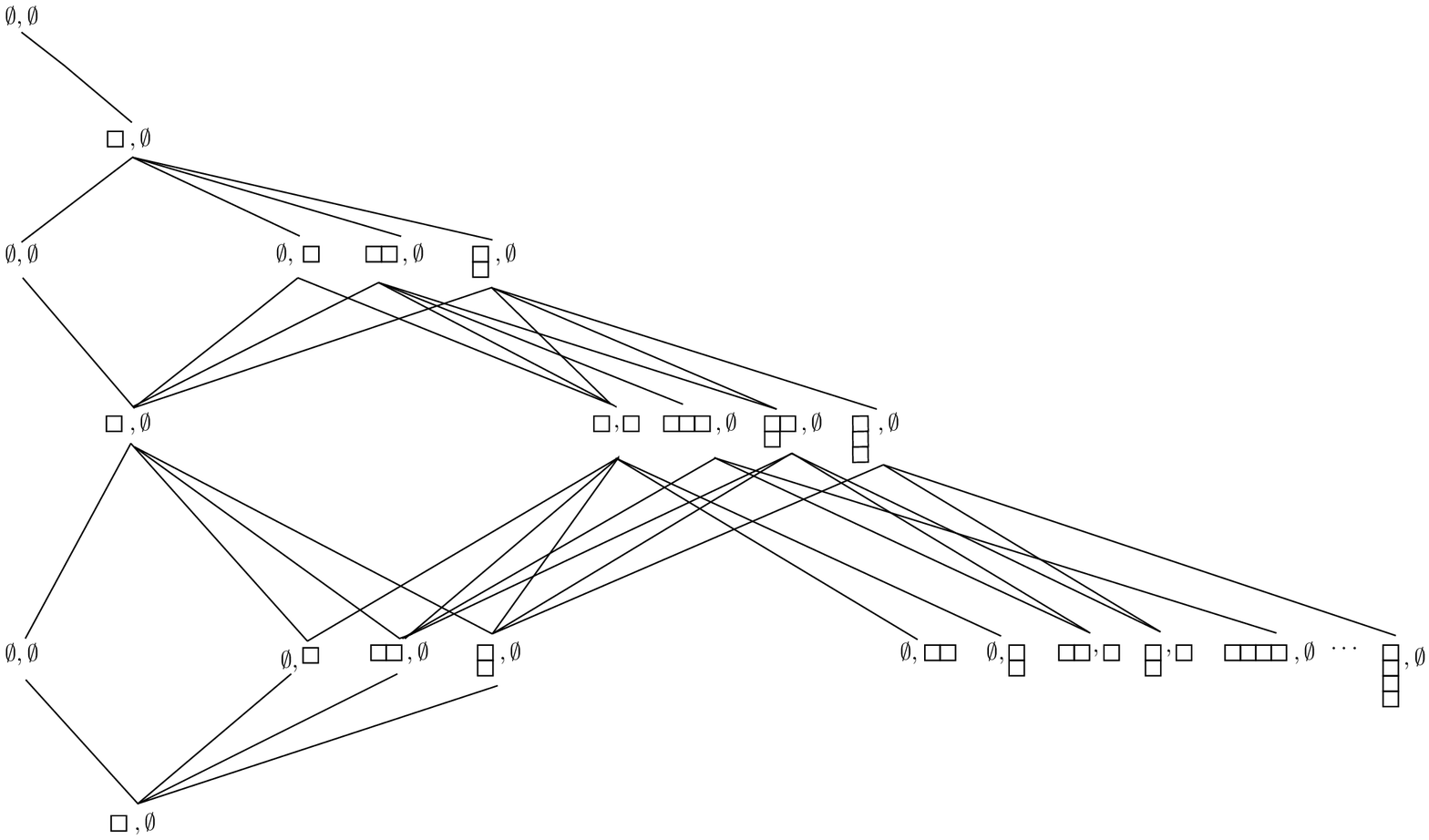} 
\] 
\caption{Standard $E_n$-module restriction diagram
  (complete up to $n=3$, partial up to $n=5$).
\label{fig:Xres22}}
\end{figure}

\begin{figure}
\[ 
\includegraphics[width=9.643cm]{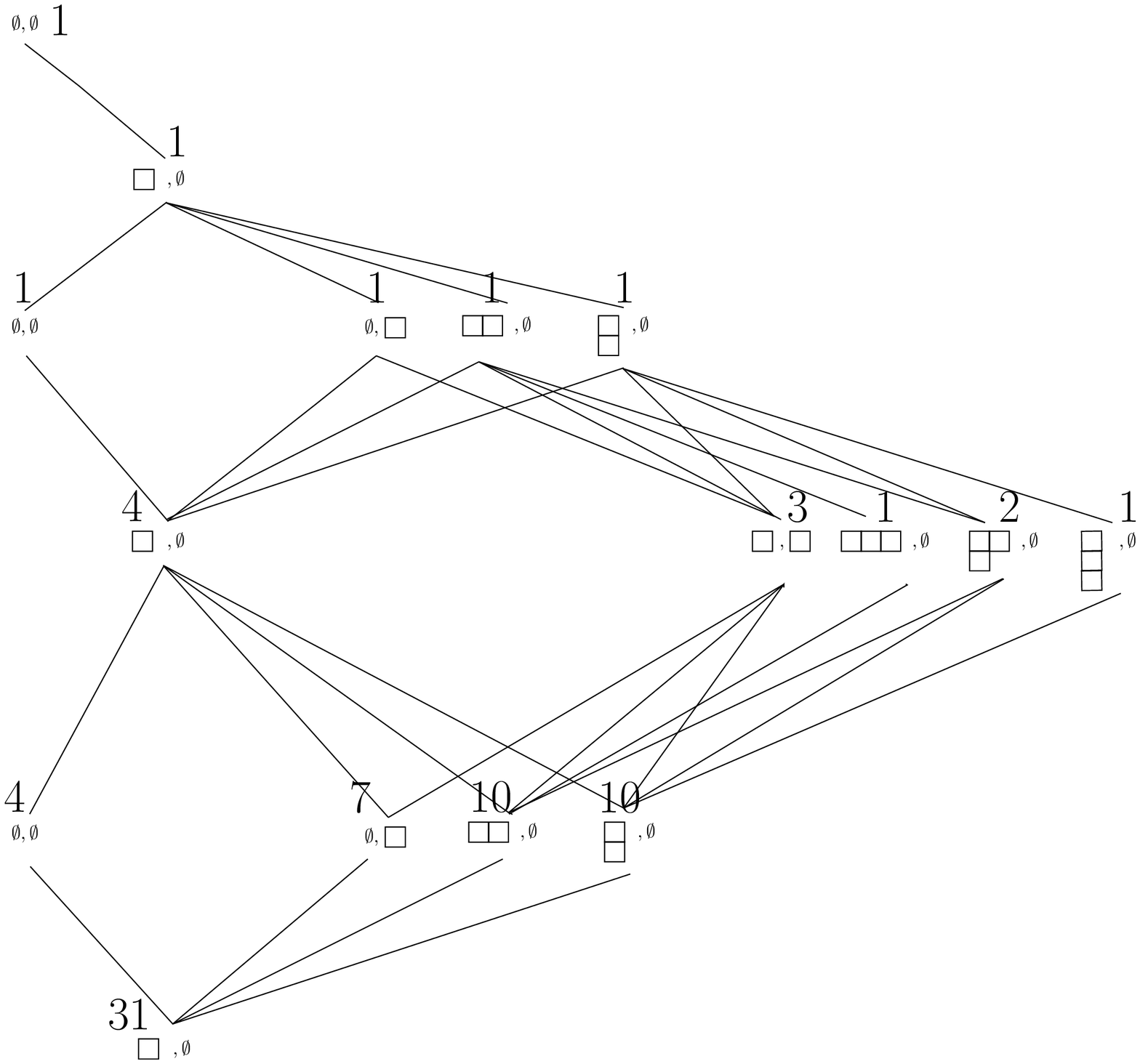}
\hspace{1.4in}
\]
\caption{Standard module restriction diagram with dimensions.
\label{fig:Xres222}}
\end{figure}


\newcommand{\lu}{\lambda+e_i}
\newcommand{\ld}{\lambda-e_i}
\newcommand{\muu}{\mu+e_j}
\newcommand{\mud}{\mu-e_j}

\mdef
Here is a schematic illustrating restriction in general position:
\[
\xymatrix{
\lu,\mud         &  \lu,\mu & \lu,\muu  \\
\lambda,\mud     & \lambda,\mu  \ar[u] \ar[d] \ar[ul] \ar[dr] 
                                         & \lambda,\muu \\
\ld,\mud         &  \ld,\mu & \ld,\muu
}
\]


\mdef It follows from (\ref{th:brat1}) 
that the Bratteli diagram takes the form in
Fig.\ref{fig:Xres22}.
It also follows that the `even Bell numbers'
(cf. e.g. \cite{AMM1})
have an intriguing
expression as a sum of squares --- a generalised Robinson--Schensted
correspondence \cite{Knuth98,MartinRollet}.
See Fig.\ref{fig:Xres222}.
\redxx{...}



\subsection{Examples: applications of modularity}

\mex \label{exa:3nonsplit}
By Theorem~\ref{th:brat1} we have the following short exact sequences  of $P^2_3$-modules
\beq \label{eq:ex1}
0 \to
  \specht^{3}_{((1),0)} +    \specht^{3}_{((1),(1))}
\to
\rest_3(\specht^{4}_{((2),0)})
\to
     \specht^{3}_{((3),0)}
       \oplus
     \specht^{3}_{((2,1),0)}
    \to 0
\eq
\begin{equation} \label{nsses}
0\to \specht^{3}_{((1),(0))}\to \rest_3(\specht^{4}_{((0),(1))})\to \specht^{3}_{((1),(1))}\to0.
\end{equation}
Let us consider the case $k=\C$ and $\delta=1$.
One can see from Fig.\ref{fig:Xres222} that $P^2_2$ is semisimple in this case,
so that $\specht^2_{((2),0)}$ is simple-projective.
It follows that $\rest_3(\specht^{4}_{((2),0)})$
is projective.
To see this note the following.


{\mlem{
In case $l=2$ we have the identification of functors
  $\ind -   \;
     \equiv \; \rest \, G_W -$.
  }}

\proof
We may essentially use the usual $P_n$ argument as for example in
\cite{Martin94} or \cite{Martin96}.
The extra requirement is to check that the $l=2$ constraint is
preserved by these manipulations.
For this the key point is the
alternative characterisation of $P^2$ as the subset of partitions of
even order (as in (\ref{de:even1})).
It follows that $P^2_{n+1} W \cong P^2(n+2,n)$ as $k$-space.
Furthermore the even property is invariant under the disk isomorphism
\cite[Appendix]{Martin2000}
$P(n,m) \to P(n-1,m+1)$, so
$P^2(n,m) \to P^2(n-1,m+1)$.
It follows that
\[
{}_{P_{n+1}} P^2_{n+2} W {}_{P_n} \cong \; {}_{P_{n+1}} P^2_{n+1} {}_{P_n}
\]
The result now follows by unpacking the various definitions. 
\qed


\mdef \label{de:indproj1}
Since 
$G_W \specht^{n}_{\lmu} = \specht^{n+2}_{\lmu}$, by (\ref{lem:Gspecht}), 
and induction preserves projectivity we see
that 
$\rest\, \specht^{n+2}_{\lmu} 
   =\rest\, G_W \specht^{n}_{\lmu} =\ind \,\specht^{n}_{\lmu}  $ 
is projective when 
$\specht^{n}_{\lmu}$ is projective.

Now by (\ref{th:brat1}) and (\ref{de:Llmu}) we see that
$\rest \,\specht_{\mulax}$, when projective, contains at least
each $P_{\lambda+e_i, \mu}$ (the indecomposable
projective cover of $L_{\lambda+e_i,\mu}$).
So in our case
\redx{\eqref{eq:ex1}}
the restriction contains 
$P_{((3),0)}$ and  $P_{((2,1),0)}$.


\mdef
Continuing with (\ref{exa:3nonsplit}), 
elementary linear algebra shows that when $\delta=1$
\beq \label{eq:la1}
\specht_{((1),0)} =L_{((1),0)} + L_{((1),(1))} .
\eq
The argument is as follows. 
Firstly, inspection of the gram matrix  shows that the socle has dimension 3.
Secondly the ideal generated by $a^{((1),(1))}$ acts as 0 on
$\specht_{\lambda,0}$
for $\lambda\vdash 3$,
but the subspace of $\specht_{((1),0)}$ on which the ideal generated
by $a^{((1),(1))}$ acts as 0 is
easily seen to be empty:
\newcommand{\aalpha}{\alpha_1}
\newcommand{\abeta}{\alpha_2}
\newcommand{\agamma}{\alpha_3}
\newcommand{\adelta}{\alpha_4}
\[
\rb{-.1in}{\ig[width=1cm]{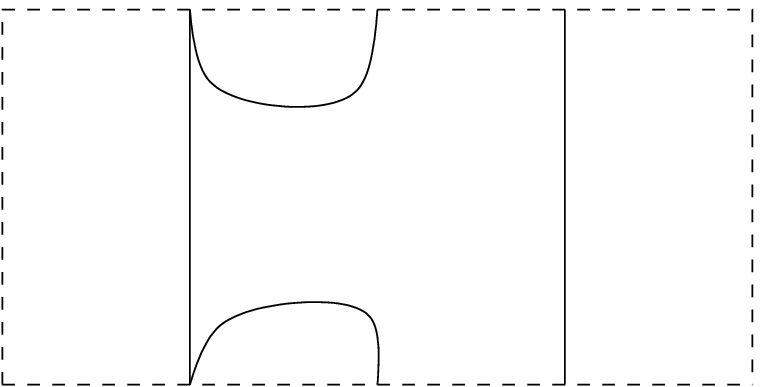}}
\left(
 \aalpha \rb{-.1in}{\ig[width=1cm]{xfig/S10b2.eps}}
+ \abeta \rb{-.1in}{\ig[width=1cm]{xfig/S10b3.eps}}
+\agamma \rb{-.1in}{\ig[width=1cm]{xfig/S10b4.eps}}
+\adelta \rb{-.1in}{\ig[width=1cm]{xfig/S10b1.eps}}
\right)
= \agamma \rb{-.1in}{\ig[width=1cm]{xfig/S10b4.eps}}
   + (\aalpha+\abeta+\adelta)  \rb{-.1in}{\ig[width=1cm]{xfig/S10b1.eps}}
\]
so the subspace has $\agamma=0$ (and $\aalpha=\abeta=0$ by symmetry, and
hence $\adelta=0$).

Meanwhile, we see from (\ref{lem:Stanss}) 
that all the $\specht_{\lambda,\mu}$ except $\specht_{((1),0)}$ are simple.
Thus in fact $\specht_{((3),0)}$ and $\specht_{((21),0)}$ are
simple-projective and 
the big sequence in (\ref{eq:ex1}) does split.
But then by (\ref{eq:BH1}) and (\ref{eq:la1}) the restriction includes
\beq \label{eq:P11}
P_{((1),(1))} = \specht_{((1),0)}+\specht_{((1),(1))} , 
\eq
so 
$\specht_{((1),0)}+\specht_{((1),(1))}$  in (\ref{eq:ex1}) does not split. 

\medskip



\ignore{{
\mex
By Theorem~\ref{th:brat1} we have the following short exact sequence  of $P^2_3$-modules
\begin{equation}\label{nsses}
0\to \specht^{3}_{((1),(0))}\to \rest_3(\specht^{4}_{((0),(1))})\to \specht^{3}_{((1),(1))}\to0.
\end{equation}
Now   let $k=\mathbb{C}$ and  $\delta=1$.
}}
By (\ref{de:indproj1}) we have that $\rest_3(\specht^{4}_{((0),(1))}) $ is
also projective.
Indeed it contains $P_{((1),(1))}$.
So by (\ref{eq:P11})  $\rest_3(\specht^{4}_{((0),(1))}) = P_{((1),(1))}$
and (\ref{nsses}) is non-split.

\ignore{{
From (\ref{}) we have 
$\rad(\specht^{3}_{((1),(0))})\simeq \specht^{3}_{((1),(1))}$ as 
a left $P^2_3(1)$-module.
The Brauer-Humphreys reciprocity implise that $$(P^3_{((1),(1))}: \specht^{3}_{((1),(0))})=[ \specht^{3}_{((1),(0))}: \specht^{3}_{((1),(1))}]\neq0$$ and $$(P^3_{((1),(1))}: \specht^{3}_{((1),(1))})\neq0.$$
  Hence  $\specht^{3}_{((1),(1))}$ and  $\specht^{3}_{((1),(0))}$ are filtration factors of $P^3_{((1),(1))}$.
The functor $\ind_n(-)$ preserves the projectivity and $\rest_3(\mathcal{S}^{4}_{((0),(1))})=\ind_2(\specht^{2}_{((0),(1))})$, this means $\rest_3(\specht^{4}_{((0),(1))})$ is projective as a left $P^2_3(1)$-module. Therefore  $P^3_{((1),(1))}$ is isomorphic to $\rest_3(\specht^{4}_{((0),(1))})$ as a left $P^2_3(1)$-module. We know $P^3_{((1),(1))}$ is indecomposable, so the sequence \ref{nsses} does not split.
}}



\ignore{{
\bigskip

\mex \red{ fix me!!! Chwas this cannot be true!} 
\textup{Let $k=\mathbb{C}$. Then $P^2_4(1)$ and $P^2_3(1)$ are not
  semisimple, and  have the following short exact sequence  of
  $P^2_3(1)$-modules 
$$
0 \to \mathcal{S}^{3}_{(\phi,(1))}\to
    \rest_3(\mathcal{S}^{4}_{(\phi,(2))})\to
    \mathcal{S}^{3}_{((1),(1))}\oplus\mathcal{S}^{3}_{(\phi,(3))}
    \oplus\mathcal{S}^{3}_{(\phi,(1))}
    \to 0
$$
which is non-split since
(some routine linear algebra shows that)
    $\rest_3(\mathcal{S}^{4}_{(\phi,(2))})$ has
no proper submodule containing
$
\mathcal{S}^{3}_{((1),(1))}\oplus\mathcal{S}^{3}_{(\phi,(3))}\oplus\mathcal{S}^{3}_{(\phi,(1))}
.
$
}
}}


\subsection{General $l$}
\newcommand{\ula}{\underline{\lambda}} 

\newcommand{\addo}{add}
\newcommand{\remo}{rem}

Let
$
\underline\lambda=
(\lambda^1,\lambda^2,\cdots,\lambda^l)\vdash\textup{\textbf{m}}\in\gamma^{l,n+1}
$
be a multi-partition.
For $1\leq j\leq l$ we define
$\addo_j(\underline{\lambda})$ to be the set of all multi-partitions
obtained from $\underline{\lambda}$ by adding an addable box to the
Young diagram of $\lambda^{j}$, and
$\remo_j(\underline{\lambda})$ to
be the set of all multi-partitions obtained from
$\underline{\lambda}$ by removing  a removable box from the Young
diagram of $\lambda^{j}$. 

 
%

{\mth{
Let $\rest_n:P^l_{n+1}-mod\to P^l_n-mod$ denote the natural
restriction corresponding to the inclusion
$P^l_n \hookrightarrow P^l_{n+1}$ given by
$d\mapsto 1\otimes d$. For
$
\underline\lambda=(\lambda^1,\lambda^2,\cdots,\lambda^l)
   \vdash\textup{\textbf{m}}\in\gamma^{l,n+1}
$
the module $\rest_n(\mathcal{S}^{n+1}_{\underline\lambda})$
has a filtration by
 standard/Specht modules.  
Specifically there is  a short exact sequence of $ P^l_n$-modules
\beq \label{eq:resses}
0 \rightarrow k A_{\ttSp} \rightarrow
 \rest_n(\mathcal{S}^{n+1}_{\underline\lambda})
    \rightarrow \bigplus_{{\underline\mu}\in \addo_{l-1}(\underline\lambda)}
\mathcal{S}^{n}_{\underline\mu} \rightarrow 0 
\eq
Where $k A_{\ttSp}= {\bigplus}_{\substack{{\underline\mu\in \addo_{i}(\remo_{i+1}(\underline\lambda))}\\ {1\leq i\leq l-1}}}                         \mathcal{S}^{n}_{\underline\mu}
\; + \;
\bigplus_{{\underline\mu}\in \addo_l(\remo_1(\underline\lambda))}\mathcal{S}^{n}_{\underline\mu}
\; + \;
\bigplus_{{\underline\mu}\in \remo_1(\underline\lambda)}\mathcal{S}^{n}_{\underline\mu}.$
\ignore{
\begin{align*}
& \rest_n(\mathcal{S}^{n+1}_{\underline\lambda})
=
\\
&   {\bigplus}_{\substack{{\underline\mu\in \addo_{i}(\remo_{i+1}(\underline\lambda))}\\ {1\leq i\leq l-1}}}                         \mathcal{S}^{n}_{\underline\mu}
\; +\;
\bigplus_{{\underline\mu}\in \addo_{l-1}(\underline\lambda)}
\mathcal{S}^{n}_{\underline\mu}
\; + \;
\bigplus_{{\underline\mu}\in \addo_l(\remo_1(\underline\lambda))}\mathcal{S}^{n}_{\underline\mu}
\; + \;
\bigplus_{{\underline\mu}\in \remo_1(\underline\lambda)}\mathcal{S}^{n}_{\underline\mu}
\end{align*}
}}}


\begin{proof}
 We  organise the set
$B^{\underline{\lambda}}_{\ttSp}$ of basis diagrams of
$\mathcal{S}^{n+1}_{\underline\lambda}$  into $l+2$
subsets denoted
$B^{\underline{\lambda}}_{\textup{sp}}(i)$,
$i=1,2,\cdots,l,l',0$, 
containing diagrams in which the first vertex is:
\begin{enumerate}
\item
  a ``propagating singleton''. This subset is denoted by
  $B^{\underline{\lambda}}_{\textup{sp}}(1)$.   
\item part of a co-$i$ propagating part, for $i=2,\cdots,l$.   
  This subset is denoted by $B^{\underline{\lambda}}_{\textup{sp}}(i)$.
\item part of a  co-1 propagating part
  $ {p}$ with $|\textup{cora}(p)|>l$.
  This subset is denoted by $B^{\underline{\lambda}}_{\textup{sp}}(l^\prime)$.
\item
  part of a non-propagating part. This subset is denoted by
  $B^{\underline{\lambda}}_{\textup{sp}}(0)$.
\end{enumerate}
  

\noindent
We proceed organisationally as for Theorem~\ref{th:brat1}.

 Case~1. This is similar to
Theorem~\ref{th:brat1}. Here the factor associated to $\lambda^1$ is the restriction of the Specht module $S_{|\lambda^1|}$ to $S_{|\lambda^1|-1}$, and this gives the submodule $\bigplus_{{\underline\mu}\in \remo_1(\underline\lambda)}\mathcal{S}^{n}_{\underline\mu}$ in (\ref{eq:resses}).

Case 2.
Here $b \in B^{\ula}_{sp}(i) $ passes on restriction to one fewer
co-$i$ line and one more co-$(i-1)$ line (a mild generalisation of the
$l=2$ case). Hence we have restriction of the factor $S_{|\lambda^i|}$ and induction in the factor $S_{|\lambda^{i-1}|}$. This case gives the submodule  
 $ {\bigplus}_{ {\underline\mu\in \addo_{i-1}(\remo_{i}(\underline\lambda))}} \mathcal{S}^{n}_{\underline\mu}$ of  $\rest_n(\mathcal{S}^{n+1}_{\underline\lambda})$ in  (\ref{eq:resses}),  for each   $i=2,\cdots,l$.

Case 3. Here we remove from co-1 and add to co-$l$. Hence we have induction in the factor $S_{|\lambda^{l}|}$ and restriction in the factor $S_{|\lambda^{1}|}$. This gives the submodule $\bigplus_{{\underline\mu}\in \addo_l(\remo_1(\underline\lambda))}\mathcal{S}^{n}_{\underline\mu}$ of  $\rest_n(\mathcal{S}^{n+1}_{\underline\lambda})$ in  (\ref{eq:resses}).

Case 4. Here we add to co-$(l-1)$. Hence we have induction in the factor $S_{|\lambda^{l-1}|}$. By the argument we have in the proof of  Theorem~\ref{th:brat1} Case~4 the factor  $\bigplus_{{\underline\mu}\in \addo_{l-1}(\underline\lambda)} \mathcal{S}^{n}_{\underline\mu}$ is not necessarily a submodule of $\rest_n(\mathcal{S}^{n+1}_{\underline\lambda})$. This leads to position of $\bigplus_{{\underline\mu}\in \addo_{l-1}(\underline\lambda)} \mathcal{S}^{n}_{\underline\mu}$ in  (\ref{eq:resses}).\end{proof}


\ignore{{
Bibliography [this is tempeorary and just to know what I have cited]

1. Klucznik, Michael, and Steffen König.
Characteristic tilting modules over quasi-hereditary algebras. Univ., 1999. 
}}





\ignore{{
\subsection{Source 
  partition categories }

In the construction of the partition category one chooses the two sets
of vertices 
in a morphism $p \in \Part_{n,m}$ to be disjoint
--- i.e.  $\underline{n}$ and $\underline{m'}$. 
However the composition does not require this restriction. 
Indeed for $S$ any set (called the  `source' set)
one may define a category $\Part_S$
by  
considering this set to be appended to every $\underline{n}$ and
$\underline{m'}$.
(This allows statistical mechanical transfer matrix computations with
a magnetic field, correlation functions,  
and other physical features, via a `source term' --- see
\cite{Martin2000}.) 

\mdef
Let $S$ be a set disjoint from any $\underline{n}$  or  $\underline{n}'$,
and $s=|S|$. 
Define
\[
\Part_{n,m}^S = \Part_{\underline{n}\cup\underline{m'}\cup S}
\]
The set $\Part_{n,m}^S$
is  in bijection with the subset $\Part_{n+s,m+s}^s$ of
$\Part_{n+s,m+s}$ of elements in which the pair $j+n,(j+m)'$ are in the same
part whenever $j \geq 1$. 

We have immediately
\mth{
For each $s \in \N_0$, 
the partition category composition closes on the collection
of subsets of form  $k\Part_{n+s,m+s}^s$. 
This defines a category 
$$
\Part_s = (\N_0 , k\Part_{n+s,m+s}^s , *).
$$
\Qed
}


\mdef
Here we will mainly focus on the category $\Part_1$.
In this case we call the unique vertex in $S$ the {\em source vertex}
(following \cite[\S 1.1]{Martin95-2000}).
We write $\Part_n^s = \Part_n^s(\delta)$ 
for the $k$-algebra that is the $n$-th end-set in
$\Part_s$. 



\mdef Note that there is an injection of 
$\Part_{n,m} \hookrightarrow \Part_{n+s,m+s}^s$ given by 
$p \mapsto p \otimes 1_s$.
In this way we realise $\Part$ as a subcategory of $\Part_s$.

}}



\section{The fusion $\chiy$-functor}
\ignore{{
\subsection{The original $\chiy$-functor} \label{ss:chiy1}

Here we recall the construction of a functor between $\Part_n^1(\delta) $-mod
and  $\Part_n(\delta-1) $-mod that gives a Morita equivalence.
Later a new but 
directly analogous functor will be used to study the \rB\ algebra.

\mdef
Still following \cite{Martin95-2000}, for each $n \in \N_0$ define 
\[
\chi_{n+1} = \prod_{i=1}^n ( 1 - A^{i\; n+1} )  \;\;\; \in \Part_{n+1}
\]
Note that this is a product of commuting idempotents. 
Note also that $\chi_{n+1} \in \Part_n^1$. 

We have 
\mpr{ {\rm \cite[Pr.4]{Martin2000}}
The map
$ \;
\psi: \chi_{n+1} \Part_n^1(\delta) \chi_{n+1}  \rightarrow P_n(\delta-1)
$
given on elements of form $ \chi_{n+1} d \chi_{n+1}$ 
where $d$ is a diagram by  $ \chi_{n+1} d \chi_{n+1} \mapsto d$ 
is a $k$-algebra isomorphism. \Qed
}

More generally, $\{ \chi_{n+1} \}_{n\in\N_0}$ constitutes a marked idempotent
in each end-set. Thus we have an idempotent subcategory of $\Part_1$:
\[
\Part_{\chi} = (\N_0 , \chi_{n+1} k \Part_{n+s,m+s}^s \chi_{m+1}, * )
\]
and hence
\mpr{CLAIM:
$\;
\Part_{\chi}^{\delta} \cong \Part^{\delta-1}
$
is a category isomorphism.  
}
\proof{needed?}

\mdef
Finally we note that $\chi_{n+1} \Part_n^1(\delta) \chi_{n+1}$ and 
$ \Part_n^1(\delta)$ are Morita equivalent (see \cite{Martin95-2000}).
In particular, then, 
the representation theory of $\Part_1$ (all $\delta$) is determined by
that of $\Part$ (all $\delta$). 
(The complex representation theory of $\Part$ is well understood, for example,
so this device completes the study of $\Part^1$ representation theory in this
case,
up to some combinatorics.)

Our job here is to use this device in the \rB\ setting.

\subsection{Marked partitions}

\mdef \label{para:marked partition}
One way of 
diagramatically encoding partitions in $\Part_1$ is to draw the vertices
that are in the same part as the source vertex in a different colour to the
rest. In this case it is not necessary to draw edges between them or
the source vertex, or indeed draw the source vertex at all.

With this twist in mind, we can ``restrict'' the entire programme
in \S\ref{ss:chiy1}
above to \rook\ partitions (i.e. to the \rB\ category
--- the subcategory where only partitions into parts with at most two elements appear).
We do this next.

In this setting one can consider the case in which the only singleton
parts are those connected to the source (i.e. all other parts are pairs). 
In this case, then, the
different colour is not needed.


\mdef
The product $*$ closes in $\Part_1$ on the span of partitions in which
parts not involving the source vertex are pair parts.
One can see that (apart from the way the second parameter is
manifested) these constructions are equivalent.

\subsection{The $\chiy$-functor for $\PPart_n$}

\mdef In the marked partition realisation of 
$\Part_1$ given in (\ref{para:marked partition})
we see that $\chi_{n+1}$ lies in the \rB\  subcategory. 
Thus we can do representation theory with $\chi_{n+1}$ in the same way as before. 
To express this in the flat deformation setting 
(i.e. where there is no $n+1$-th vertex)
we need the image of
$A^{i \; i+1}$: 
\[
A^{i \; i+1} \mapsto \; \frac{1}{\delta'}  A^{i.} \;  
\]
(Note that nothing like this can work in the full partition algebra
setting. It relies on restriction to the \rB\ case.)
So
\[
\chi_{n+1} \; \mapsto \; \prod_{i=1}^n (1 - \frac{1}{\delta'} A^{i.})
\]
Hereafter we may write $\chiy$ for this idempotent where no ambiguity
arises. 

}}
\newcommand{\epi}{e_{\pi}}

Suppose $n$ even and define
\[
\epi = a^{\otimes n/2} = A^{12} A^{34} ... A^{n-1 \; n}
\]
{\mpr{
$\epi E_n \epi \cong P_{n/2}$
}}
\proof Consider the picture:
\[
\includegraphics[width=4.3cm]{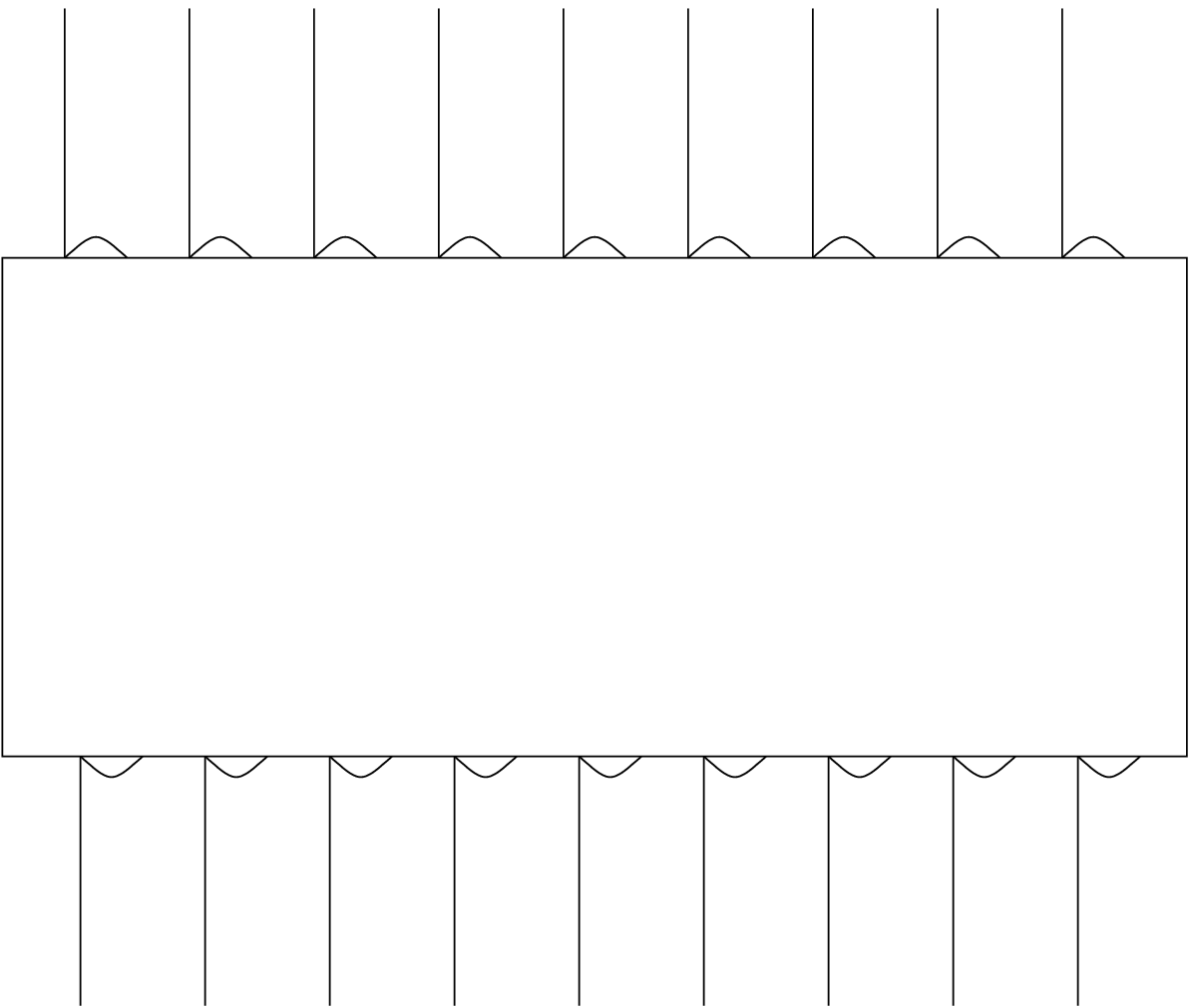} 
\]
It will be evident that inserting $E_n$ into the box defines a
map to $P_{n/2}$. It is straightforward to check that the map is
surjective. 
\qed

It follows that $P_{n/2}$-mod fully embeds in $E_n$-mod, and hence that
$E_n$ is non-semisimple whenever $P_{n/2}$ is non-semisimple. 
The structure of $P_n$ is given in \cite{Martin96}.
Thus we have determined another (small but interesting)
part of the structure of $E_n$. 
In particular we may deduce  that non-negative integer $\delta$ values are
non-semisimple for sufficienly large $n$.
(Confer for example
\cite{MazorchukMartin10},
\cite{GanyushkinMazorchuk09}.)




\medskip \bigskip

\noindent
    {\bf Acknowledgments.} 
    CA thanks the KRG for scholarship funding
    and then the University of Sulaimani for support. 
PM thanks EPSRC for funding under grant EP/I038683/1.
VM thanks also EPSRC for partial support during the visit to
Leeds during which this project was started.
VM is partially supported by the Swedish Research Council and
G\"oran Gustafsson Foundation.

\clearpage \appendix \section*{Appendix} 
\section{Partition category composition} \label{ss:pdefa}
\ignore{{  
\newcommand{\Tri}{\mbox{Tri}}  
\newcommand{\Dgroup}{{\mathcal D}} 
\newcommand{\Aregular}{$A$-regular}
\newcommand{\nngraph}{\mbox{-graph}}
\newcommand{\hash}{\#}
\newcommand{\ex}{\sigma}  
\newcommand{\chiy}{\chi}  
\newcommand{\prop}{propagating}
\newcommand{\propno}{\hash^p}
\newcommand{\MM}{{\mathcal M}} 
\newcommand{\chiyn}{\chiy_{}}
\newcommand{\monotone}{monotone}
\newcommand{\catP}{{\mathcal P}}
\newcommand{\catT}{{\mathcal T}}
\newcommand{\JB}{J}  
\newcommand{\J}{J}  


\renewcommand{\nome}[2]{}
\newcommand{\href}[2]{#2}
\newcommand{\highdetail}[1]{}
\newcommand{\complete}[1]{\check{#1}}
\newcommand{\nullseq}{{\mathfrak n}}
\newcommand{\RHS}{RHS}
\newcommand{\RTP}{RTS}

\newcommand{\mdeff}[1]{{\sc #1}}
\newcommand{\mdefff}[2]{{\sc #2}}
\newcommand{\mdeft}[4]{{ #2} & #3 $\;$ #4}


\newcommand{\Trace}{\mbox{Tr}} 
\newcommand{\Power}{P}  
\newcommand{\Part}{P}  
\newcommand{\Pset}{{\mathsf P}} 
\newcommand{\Pair}{J}  
\newcommand{\pcirc}{\circ}  
\newcommand{\dcirc}{\oslash}  
\newcommand{\ER}{E}  
\newcommand{\udot}{u}   
\newcommand{\bdot}{u^*} 
\newcommand{\RB}{R}   
\newcommand{\PPart}{\RB} 
\newcommand{\PPPart}{\RB'} 
\newcommand{\PParto}{\RB^{o}} 
\newcommand{\PParte}{\RB^{e}} 
\newcommand{\pPart}{B} 

\newcommand{\CC}{{\mathcal C}} 
\newcommand{\rook}{partial}
\newcommand{\Rook}{Partial}

\newcommand{\specht}{{\mathcal S}}  

\newcommand{\muth}{{\textsc{Theorem}.}\ } 
\newcommand{\mulem}{{\textsc{Lemma}.}\ }

\newcommand{\myceil}[1]{\left \lceil #1 \right \rceil }
\newcommand{\floor}[1]{\left \lfloor #1 \right \rfloor }
\newcommand{\flip}{\star}  
\newcommand{\mE}{\mathtt{E}} 
\newcommand{\mP}{\mathtt{P}} 
\newcommand{\mD}{\mathtt{D}} 
}}

\newcommand{\un}[1]{\underline{#1}} 
\newcommand{\nN}{{\mathcal N}}

For $S$ a set 
let $\mP_S$ denote the set of set partitions of $S$.
Let 
$\underline{n} := \{ 1,2,..,n \}$ 
and 
$\underline{n'} := \{ 1',2',..,n' \}$ 
and so on; and 
$P_{n,m} \; := \; \mP_{\un{n} \cup \un{m'}}$.

Given a graph $g=(V,E)$ then $\pi(g) \in \mP_V$ denotes the partition
according to connected components of $g$. 
Given a graph $g=(V,E)$ with $V \supseteq \un{n}\cup\un{m'}$
then $\pi_{nm}(g) \in P_{n,m}$ denotes the partition according to
connected components of $g$.

Given graphs $g=(V,E)$, $g'=(V',E')$ then graph
\[
g.g' = (V \cup V', E \sqcup E')
\]
Given a graph $g$ with $V=\un{n}\cup\un{m'}$ and a graph $g'$ with
$V'= \un{m}\cup\un{l'}$ then 
\beq \label{eq:gbarg}
g|g' = g_+ . g'_-
\eq
where graph $g_+ =(V_+,E_+)$ 
is $g$ with vertices $i'$ ($i\in \un{m}$) replaced by $i''$;
and graph $g'_-$ is $g'$ with vertices $i \in \un{m}$ replaced by $i''$
(so $V_+ \cap V'_- = \un{m''}$).

A connected component of $g|g'$ is {\em internal} if it intersects
neither $\un{n}$ or $\un{l'}$. Define $c(g|g')$ as the number of
internal components. 
Define 
\[
*:P_{n,m} \times P_{m,l} \rightarrow kP_{n,l}
\]
as follows. For $p \in P_{n,m}$ pick any graph with $\pi(g)=p$;
and similarly for $p' \in P_{m,l}$. 
Then
\[
p*p' = \delta^{c(g|g')} \pi_{nl}(g|g')
\]

{\mth{ {\rm \cite{Martin94}}
The composition $*$ is well-defined and, extended $k$-linearly, 
makes $kP_{n,n}$ an associative unital algebra;
and $\left( kP_{n,m} \right)_{n,m \in \N_0}$ a $k$-linear category.   
\qed
}}






\bibliographystyle{amsplain}
\bibliography{bib/new31,bib/bb,bib/bbc,
local}

\providecommand{\bysame}{\leavevmode\hbox to3em{\hrulefill}\thinspace}
\providecommand{\MR}{\relax\ifhmode\unskip\space\fi MR }
\providecommand{\MRhref}[2]{%
  \href{http://www.ams.org/mathscinet-getitem?mr=#1}{#2}
}
\providecommand{\href}[2]{#2}
\begin{thebibliography}{10}

\bibitem{AMM1}
C~Ahmed, P~P Martin, and V~Mazorchuk, \emph{On the number of principal ideals
  in d-tonal partition monoids}, arXiv:1503.06718 [math.CO] (2015).

\bibitem{Assem06}
I~Assem, D~Simson, and A~Skowronski, \emph{Elements of the representation
  theory of associative algebras}, Cambridge UP, 2006.

\bibitem{Benson95}
D~J Benson, \emph{Representations and cohomology {I}}, Cambridge, 1995.

\bibitem{boltje2012twisted}
R~Boltje and S~Danz, \emph{Twisted split category algebras as quasi-hereditary
  algebras}, Archiv der Mathematik \textbf{99} (2012), no.~6, 589--600.

\bibitem{Brauer37}
R~Brauer, \emph{On algebras which are connected with the semi--simple
  continuous groups}, Annals of Mathematics \textbf{38} (1937), 854--872.

\bibitem{Brauer39}
\bysame, \emph{On modular and p-adic representations of algebras}, Proc Nat
  Acad Sci USA \textbf{25} (1939), 252--258.

\bibitem{CPS}
E~Cline, B~Parshall, and L~Scott, J. reine angew Math. \textbf{391} (1988), 85.

\bibitem{CoxMartinParkerXi06}
A~G Cox, P~P Martin, A~E Parker, and C~C Xi, \emph{Representation theory of
  towers of recollement: theory, notes and examples}, J Algebra \textbf{302}
  (2006), 340--360, DOI 10.1016 online (math.RT/0411395).

\bibitem{CoxDeVisscherMartin0609}
A~G Cox, M~De Visscher, and P~P Martin, \emph{A geometric characterisation of
  the blocks of the {B}rauer algebra}, JLMS \textbf{80} (2009), 471--494,
  (math.RT/0612584).

\bibitem{CurtisReiner62}
C~W Curtis and I~Reiner, \emph{Representation theory of finite groups and
  associative algebras}, Wiley Interscience, New York, 1962.

\bibitem{CurtisReiner90}
\bysame, \emph{Methods of representation theory with applications to finite
  groups and orders}, vol.~1, Wiley, New York, 1990.

\bibitem{DlabRingel}
V~Dlab and C~M Ringel, Compositio Mathematica \textbf{70} (1989), 155--175.

\bibitem{EhrigStroppel13}
M~Ehrig and C~Stroppel, \emph{{Nazarov-{W}enzl Algebras, Coideal Subalgebras
  and Categorified Skew {H}owe Duality}},  (2013), arXiv:1310.1972.

\bibitem{EhrigStroppel14b}
\bysame, \emph{{Schur-{W}eyl Duality for the {B}rauer Algebra and the
  Ortho-symplectic {L}ie Superalgebra }},  (2014), arXiv:1412.7853.

\bibitem{enyang2017cellular}
J~Enyang and F~M Goodman, \emph{Cellular bases for algebras with a jones basic
  construction}, Algebras and Representation Theory \textbf{20} (2017), no.~1,
  71--121.

\bibitem{ErdmannSaenz03}
K~Erdmann and C~Saenz, \emph{On standardly stratified algebras}, Communications
  in Algebra \textbf{31} (2003), no.~7, 3429--3446.

\bibitem{GanyushkinMazorchuk09}
O~Ganyushkin and V~Mazorchuk, \emph{Classical finite transformation semigroups:
  an introduction}, Springer-Verlag, 2009.

\bibitem{GoodmanWallach98}
R~Goodman and N~R Wallach, \emph{Representations and invariants of the
  classical groups}, Cambridge, 1998.

\bibitem{GrahamLehrer96}
J~J Graham and G~I Lehrer, \emph{Cellular algebras}, Invent. Math. \textbf{123}
  (1996), 1--34.

\bibitem{Green80}
J~A Green, \emph{Polynomial representations of ${GL}_n$}, Springer-Verlag,
  Berlin, 1980.

\bibitem{green2004standard}
R~M Green, \emph{Standard modules for tabular algebras}, Algebras and
  representation theory \textbf{7} (2004), no.~4, 419--440.

\bibitem{James}
G~D James, \emph{The representation theory of the symmetric groups}, Lecture
  Notes in Mathematics 682, Springer, 1978.

\bibitem{JamesKerber}
G~D James and A~Kerber, \emph{The representation theory of the symmetric
  group}, Addison-Wesley, London, 1981.

\bibitem{Jantzen}
J~C Jantzen, \emph{Representations of algebraic groups}, Academic Press, 1987.

\bibitem{Jones94}
V~F~R Jones, \emph{The {P}otts model and the symmetric group}, in {S}ubfactors:
  {P}roceedings of the {T}aniguchi {S}ymposium on {O}perator {A}lgebras, World
  Scientific, Singapore, 1994.

\bibitem{KS1999}
M~Klucznik and S~Koenig, \emph{\emph{Characteristic tilting modules over
  quasi-hereditary algebras}}, Bar Ilan / Bielefeld Univ. preprint (1999).

\bibitem{Knuth98}
D~E Knuth, \emph{Sorting and searching}, 2 ed., The Art of Computer
  Programming, vol.~3, Addison Wesley, 1998.

\bibitem{Kosudax}
M~Kosuda, \emph{Party algebra of type b and construction of its irreducible
  representations}, Formal Power Series and Algebraic Combinatorics 2005
  Proceedings, ed. L Carini, H Barcello, J-Y Thibon, manuscript online (2005).

\bibitem{Kosuda06}
\bysame, \emph{Irreducible representations of the party algebra}, Osaka J Math
  \textbf{43} (2006), 431--474.

\bibitem{Kosuda08}
\bysame, \emph{Characterization for the modular party algebra}, Journal of Knot
  Theory and Its Ramifications \textbf{17} (2008), 939--960.

\bibitem{lam2012lectures}
T~Y Lam, \emph{Lectures on modules and rings}, vol. 189, Springer Science \&
  Business Media, 2012.

\bibitem{LehrerZhang14}
G~Lehrer and R~Zhang, \emph{The first fundamental theorem of invariant theory
  for the orthosymplectic supergroup}, arxiv:1401.7395.

\bibitem{Martin91}
P~P Martin, \emph{Potts models and related problems in statistical mechanics},
  World Scientific, Singapore, 1991.

\bibitem{Martin94}
\bysame, \emph{Temperley--{L}ieb algebras for non--planar statistical mechanics
  --- the partition algebra construction}, Journal of Knot Theory and its
  Ramifications \textbf{3} (1994), no.~1, 51--82.

\bibitem{Martin96}
\bysame, \emph{The structure of the partition algebras}, J Algebra \textbf{183}
  (1996), 319--358.

\bibitem{Martin2000}
\bysame, \emph{The partition algebra and the {P}otts model transfer matrix
  spectrum in high dimensions}, J Phys A \textbf{32} (2000), 3669--3695.

\bibitem{Martin08a}
\bysame, \emph{On diagram categories, representation theory and statistical
  mechanics}, AMS Contemp Math \textbf{456} (2008), 99--136.

\bibitem{Martin0915}
\bysame, \emph{The decomposition matrices of the {B}rauer algebra over the
  complex field}, Trans. A.M.S. \textbf{367} (2015), 1797--1825,
  (http://arxiv.org/abs/0908.1500).

\bibitem{MazorchukMartin10}
P~P Martin and V~Mazorchuk, \emph{On the representation theory of partial
  brauer algebras}, The Quarterly Journal of Mathematics \textbf{65} (2013),
  no.~1, 225--247.

\bibitem{MartinRollet}
P~P Martin and G~Rollet, \emph{The {P}otts model representation and a
  {R}obinson--{S}chensted correspondence for the partition algebra}, Compositio
  Math \textbf{112} (1998), 237--254, I have a copy.

\bibitem{MartinSaleur94b}
P~P Martin and H~Saleur, \emph{Algebras in higher dimensional statistical
  mechanics --- the exceptional partition algebras}, Lett. Math. Phys.
  \textbf{30} (1994), 179--185, (hep-th/9302095).

\bibitem{MartinWoodcock98}
P~P Martin and D~Woodcock, \emph{The partition algebras and a new deformation
  of the {S}chur algebras}, J Algebra \textbf{203} (1998), 91--124.

\bibitem{Martin}
S~Martin, \emph{Schur algebras and representation theory}, Cambridge University
  Press, 1993.

\bibitem{Orellana07}
R~Orellana, \emph{On partition algebras for complex reflection groups}, J
  Algebra \textbf{313} (2007), 590--616.

\bibitem{Peel75}
M~H Peel, \emph{Specht modules and symmetric groups}, J Algebra \textbf{36}
  (1975), 88--97.

\bibitem{Tanabe97}
K~Tanabe, \emph{On the centralizer algebra of the unitary reflection group
  $g(m,p,n)$}, Nagoya Math J \textbf{148} (1997), 113--126.

\bibitem{TemperleyLieb71}
H~N~V Temperley and E~H Lieb, \emph{Relations between percolation and colouring
  problems and other graph theoretical problems associated with regular planar
  lattices: some exact results for the percolation problem}, Proceedings of the
  Royal Society A \textbf{322} (1971), 251--280.

\end{thebibliography}


\printindex

\end{document}